\documentclass[11pt]{article}
\usepackage[english]{babel}
\usepackage[latin1]{inputenc}
\usepackage{amsmath}
\usepackage{amssymb}
\usepackage{amsfonts}
\usepackage{amsthm}
\usepackage{upgreek}
\usepackage{stmaryrd}

\newcommand{\n}{\mathfrak{n}}

\newcommand{\w}{\mathfrak{w}}
\newcommand{\x}{{\rm x}}
\newcommand{\C}{\mathbb{C}}
\newcommand{\N}{\mathbb{N}}

\newcommand{\R}{\mathbb{R}}

\newcommand{\boA}{\mathcal{A}}
\newcommand{\boC}{\mathcal{C}}

\newcommand{\boK}{\mathcal{K}}
\newcommand{\boL}{\mathcal{L}}
\newcommand{\boM}{\mathcal{M}}
\newcommand{\boN}{\mathcal{N}}
\newcommand{\boO}{\mathcal{O}}

\newcommand{\boT}{\mathcal{T}}
\newcommand{\boU}{\mathcal{U}}

\newcommand{\boW}{\mathcal{W}}
\newcommand{\boX}{\mathcal{X}}
\newcommand{\boZ}{\mathcal{Z}}

\newcommand{\eps}{\varepsilon}

\newcommand{\sign}{{\rm sign}}

\newtheorem{case}{Case}
\newtheorem{claim}{Claim}

\newtheorem{lemma}{Lemma}
\newtheorem{prop}{Proposition}
\newtheorem{step}{Step}
\newtheorem{theorem}{Theorem}
\theoremstyle{definition}
\newtheorem*{merci}{Acknowledgements}
\newtheorem{remark}{Remark}

\begin{document}

\title{On the Korteweg-de Vries long-wave approximation of the Gross-Pitaevskii equation II}
\author{
\renewcommand{\thefootnote}{\arabic{footnote}}
Fabrice B\'ethuel \footnotemark[1], Philippe Gravejat \footnotemark[2], Jean-Claude Saut \footnotemark[3], Didier Smets \footnotemark[4]}
\footnotetext[1]{UPMC, Universit\'e Paris 06, UMR 7598, Laboratoire Jacques-Louis Lions, F-75005, Paris, France. E-mail: bethuel@ann.jussieu.fr}
\footnotetext[2]{Centre de Recherche en Math\'ematiques de la D\'ecision, Universit\'e Paris Dauphine, Place du Mar\'echal De Lattre De Tassigny, 75775 Paris Cedex 16, France, and \'Ecole Normale Sup\'erieure, DMA, UMR 8553, F-75005, Paris, France. E-mail: gravejat@ceremade.dauphine.fr}
\footnotetext[3]{Laboratoire de Math\'ematiques, Universit\'e Paris Sud and CNRS UMR 8628, B\^atiment 425, 91405 Orsay Cedex, France. E-mail: Jean-Claude.Saut@math.u-psud.fr}
\footnotetext[4]{UPMC, Universit\'e Paris 06, UMR 7598, Laboratoire Jacques-Louis Lions, F-75005, Paris, France, and \'Ecole Normale Sup\'erieure, DMA, UMR 8553, F-75005, Paris, France. E-mail: smets@ann.jussieu.fr}
\maketitle

\begin{abstract}
In this paper, we proceed along our analysis of the Korteweg-de Vries approximation of the Gross-Pitaevskii equation initiated in \cite{BeGrSaS2}. At the long-wave limit, we establish that solutions of small amplitude to the one-dimensional Gross-Pitaevskii equation split into two waves with opposite constant speeds $\pm \sqrt{2}$, each of which are solutions to a Korteweg-de Vries equation. We also compute an estimate of the error term which is somewhat optimal as long as travelling waves are considered. At the cost of higher regularity of the initial data, this improves our previous estimate in \cite{BeGrSaS2}.
\end{abstract}

\section{Introduction}

\subsection{Statement of the results}

In this paper, we proceed along our study initiated in \cite{BeGrSaS2} of the one-dimensional Gross-Pitaevskii equation
\renewcommand{\theequation}{GP}
\begin{equation}
\label{GP}
i \partial_t \Psi + \partial_\x^2 \Psi = \Psi (|\Psi|^2 - 1) \ {\rm on} \ \R \times \R,
\end{equation}
supplemented with the boundary condition at infinity
$$|\Psi(\x, t)| \to 1, \ {\rm as} \ |\x| \to + \infty.$$
This boundary condition is suggested by the formal conservation of the Ginzburg-Landau energy
$$E(\Psi) = \frac{1}{2} \int_{\R} |\partial_\x \Psi|^2 + \frac{1}{4} \int_{\R} (1 - |\Psi|^2)^2.$$
In this paper, we will only consider finite energy solutions to \eqref{GP}.

The Gross-Pitaevskii equation is integrable by means of the inverse scattering method, and it has been formally analyzed within this framework in \cite{ShabZak2}, and rigorously in \cite{GeraZha1}. Concerning the Cauchy problem, it can be shown (see \cite{Zhidkov1,Gerard2,BeGrSaS2}) that \eqref{GP} is globally well-posed in the spaces
$$X^k(\R) = \big\{ u \in L^1_{\rm loc}(\R, \C), \ {\rm s.t.} \ 1 - |u| ^2 \in L^2(\R) \ {\rm and} \ \partial_\x u \in H^{k - 1}(\R) \big\},$$
for any $k \geq 1$. More precisely, we have

\begin{prop}[\cite{BeGrSaS2}]
\label{thm:existe}
Let $k \in \N^*$ and $\Psi_0 \in X^k(\R)$. Then, there exists a unique solution $\Psi(\cdot, t)$ in $\boC^0(\R, X^k(\R))$
\footnote{Here, the space $X^k(\R)$ is endowed with the distance
$$d_{X^k, A}(u, v) = \| u - v \|_{L^\infty([- A, A])} + \| \partial_\x u - \partial_\x v \|_{H^{k - 1}(\R)} + \| |u| - |v| \|_{L^2(\R)},$$
for some given $A > 0$ (see e.g. \cite{Gerard2,BeGrSaS1} for more details).}
to \eqref{GP} with initial data $\Psi_0$. Furthermore, the energy $E$ is conserved along the flow.
\end{prop}

If $u$ belongs to $X^1(\R)$ and satisfies
$$E(u) < \frac{2 \sqrt{2}}{3},$$
then it does not vanish, and we may write $u= |u| \exp i \theta$, where $\theta$ is continuous (see e.g. \cite{BetGrSa2}). Here, we will focus on solutions with small energy, so that in view of the conservation of the energy, we may write
$$\Psi(\cdot, t) = \varrho(\cdot, t) \exp i \varphi(\cdot, t).$$
More precisely, we will consider initial data which are small long-wave perturbations of the constant one, namely
$$\left\{ \begin{array}{ll}
\varrho(\x, 0) = \Big( 1 - \frac{\varepsilon^2}{6} N_\varepsilon^0(\varepsilon \x) \Big)^\frac{1}{2},\\
\varphi (\x, 0) = \frac{\varepsilon}{6 \sqrt{2}} \Theta_\varepsilon^0(\varepsilon \x), \end{array} \right.$$
where $0 < \varepsilon < 1$ is a small parameter, and $N_\varepsilon^0$ and $W_\varepsilon^0 = \partial_x \Theta_\varepsilon^0$ are uniformly bounded in some Sobolev spaces $H^k(\R)$ for sufficiently large $k$. We will add two additional assumptions on $\Theta_\varepsilon^0$ and $N_\varepsilon^0$. We will assume that
\renewcommand{\theequation}{\arabic{equation}}
\setcounter{equation}{0}
\begin{equation}
\label{H1}
\| N_\varepsilon^0 \|_{\boM(\R)}+ \| \partial_x \Theta_\varepsilon^0 \|_{\boM(\R)} < + \infty,
\end{equation}
with a uniform bound in $\varepsilon$. Here, $\| \cdot \|_{\boM(\R)}$ denotes the norm defined on $L_{\rm loc}^1(\R)$ by
\begin{equation}
\label{massnorm}
\| f \|_{\boM(\R)} = \underset{(a, b) \in \R^2}{\sup} \bigg| \int_a^b f(x) dx \bigg|,
\end{equation}
so that \eqref{H1} implies in particular that $\Theta_\eps^0$ is uniformly bounded in $L^\infty(\R)$. In the appendix, we will introduce a notion of mass for $\Psi$ closely related to the $\boM$-norm of $1 - |\Psi|^2$, and prove its conservation by the Gross-Pitaevskii flow.

We next introduce the slow coordinates
$$x^- = \varepsilon (\x + \sqrt{2} t), \ x^+ = \varepsilon (\x - \sqrt{2} t), \ {\rm and} \ \tau = \frac{\varepsilon^3}{2 \sqrt{2}} t.$$
The definition of the new coordinates $x^-$ and $x^+$ corresponds to reference frames travelling to the left and to the right respectively with speed $\sqrt{2}$ in the original coordinates $(\x, t)$. We define accordingly the rescaled functions $N_\varepsilon^\pm$ and $\Theta_\varepsilon^\pm$ as follows
\begin{equation}
\label{slow-var}
\begin{split}
N_\varepsilon^\pm(x^\pm, \tau) & = \frac{6}{\varepsilon^2} \eta(\x, t) = \frac{6}{\varepsilon^2} \eta \Big( \frac{x^\pm}{\varepsilon} \pm \frac{4 \tau}{\varepsilon^3}, \frac{2 \sqrt{2} \tau}{\varepsilon^3} \Big),\\
\Theta_\varepsilon^\pm(x^\pm, \tau) & = \frac{6 \sqrt{2}}{\varepsilon} \varphi(\x, t) =
\frac{6 \sqrt{2}}{\varepsilon} \varphi \Big( \frac{x^\pm}{\varepsilon} \pm \frac{4
\tau}{\varepsilon^3}, \frac{2 \sqrt{2} \tau}{\varepsilon^3} \Big),
\end{split}
\end{equation}
where $\eta = 1 - \varrho^2$. Setting
\begin{equation}
\label{eq:uv}
\begin{split}
U_\varepsilon^-(x^-, \tau) = \frac{1}{2} \Big( N_\varepsilon^-(x^-, \tau) + \partial_{x^-} \Theta_\varepsilon^-(x^-, \tau) \Big),\\
U_\varepsilon^+(x^+, \tau) = \frac{1}{2} \Big( N_\varepsilon^+(x^+, \tau) - \partial_{x^+} \Theta_\varepsilon^+(x^+, \tau) \Big),
\end{split}
\end{equation}
our main result is

\begin{theorem}
\label{cochondore}
Let $k \geq 0$ and $\varepsilon > 0$ be given. Assume that the initial data $\Psi_0$ belongs to $X^{k + 6}(\R)$ and satisfies the assumption
\begin{equation}
\label{grinzing1}
\| N_\varepsilon^0 \|_{\boM(\R)} + \|\partial_x \Theta_\varepsilon^0 \|_{\boM(\R)} + \| N_\varepsilon^0 \|_{H^{k + 5}(\R)} + \varepsilon \| \partial_x^{k + 6} N_\varepsilon^0 \|_{L^2(\R)} + \|\partial_x \Theta_\varepsilon^0 \|_{H^{k + 5}(\R)} \leq K_0.
\end{equation}
Let $\boU^-$ and $\boU^+$ denote the solutions to the Korteweg-de Vries equations
\renewcommand{\theequation}{KdV}
\begin{equation}
\label{KdV}
\partial_\tau \boU^- + \partial_{x^-}^3 \boU^- + \boU^- \partial_{x^-} \boU^- = 0,
\end{equation}
and
\renewcommand{\theequation}{\arabic{equation}}
\setcounter{equation}{5}
\begin{equation}
\label{KdV-}
\partial_\tau \boU^ + - \partial_{x^+}^3 \boU^+ - \boU^+ \partial_{x^+} \boU^+ = 0,
\end{equation}
with the same initial value
\footnote{Since their initial data depend on $\varepsilon$, $\boU^-$ and $\boU^+$ do as well. We voluntarily hide this dependence in the notations in order to stress the fact that the equations they satisfy are independent of $\varepsilon$.}
as $U_\varepsilon^-$, respectively $U_\varepsilon^+$. Then, there exist positive constants $\varepsilon_1$ and $K_1$, depending only on $k$ and $K_0$, such that
\begin{equation}
\label{ineq1}
\| U_\varepsilon^-(\cdot, \tau) - \boU^-(\cdot, \tau) \|_{H^k(\R)} + \|U_\varepsilon^+(\cdot, \tau) - \boU^+(\cdot, \tau) \|_{H^k(\R)} \leq K_1 \varepsilon^2 \exp K_1 |\tau|,
\end{equation}
for any $\tau \in \R$ provided $\varepsilon \leq \varepsilon_1$.
\end{theorem}

\begin{remark}
In the original time variable, the Korteweg-de Vries approximation is valid on a time interval $t \in [0, T_\eps]$ with
$$T_\eps = o \bigg( \frac{|\log(\eps)|}{\eps^3} \bigg).$$
Moreover, the approximation error remains of order $\boO(\eps^2)$ on a time interval $t \in [0, T_\eps']$ with $T_\eps' = \boO(\eps^{- 3})$.
\end{remark}

In order to explain the statements of Theorem \ref{cochondore}, it is presumably useful to recast them in the context of known results about the long-wave limit of the Gross-Pitaevskii equation. First, we rewrite \eqref{GP} in the slow coordinates $(y, s) = (\varepsilon \x, \varepsilon t)$ and set
$$\left\{ \begin{array}{ll}
\varrho(\x, t) = \Big( 1 - \frac{\varepsilon^2}{6} n_\varepsilon(\varepsilon \x, \varepsilon t) \Big)^\frac{1}{2},\\
\partial_{\x} \varphi (\x, t) = \frac{\varepsilon^2}{6 \sqrt{2}} w_\varepsilon(\varepsilon \x, \varepsilon t).
\end{array} \right.$$
In this setting, \eqref{GP} translates into the system for $n_\varepsilon$ and $w_\varepsilon$,
\begin{equation}
\label{eq:dynaslow}
\left\{ \begin{array}{ll}
\partial_s n_\varepsilon - \sqrt{2} \partial_y w_\varepsilon = - \frac{\varepsilon^2}{3 \sqrt{2}} \partial_y (n_\varepsilon w_\varepsilon),\\
\partial_s w_\varepsilon - \sqrt{2} \partial_y n_\varepsilon = - 3 \sqrt{2} \varepsilon^2 \partial_y \Big( \frac{\partial_y^2 n_\varepsilon}{6 - \varepsilon^2 n_\varepsilon} + \frac{w_\varepsilon^2}{36} + \frac{\varepsilon^2}{2} \frac{(\partial_y n_\varepsilon)^2}{(6 - \varepsilon^2 n_\varepsilon)^2} \Big).
\end{array} \right.
\end{equation}
It has been shown in \cite{BetDaSm1} that for suitably small data and times, this system is well-approximated by the linear wave equation. More precisely, assume that $s \geq 2$ and
$$K(s) \varepsilon^2 \| (N_\varepsilon^0, W_\varepsilon^0) \|_{H^{s + 1}(\R) \times H^s(\R)} \leq 1,$$
where $K(s)$ refers to some positive constant depending only on $s$. Let $(\n, \w)$ denote the solution of the free wave equation
\begin{equation}
\label{eq:wamu}
\left\{ \begin{array}{ll}
\partial_s \n - \sqrt{2} \partial_y \w = 0,\\
\partial_s \w - \sqrt{2} \partial_y \n = 0,
\end{array} \right.
\end{equation}
with initial data $(N_\varepsilon^0, W_\varepsilon^0)$. Then, for any $0 \leq t \leq T_\varepsilon$, we have
\begin{align*}
& \|(n_\varepsilon, w_\varepsilon)(\cdot, \varepsilon t) - (\n, \w)(\cdot, \varepsilon t)\|_{H^{s - 2}(\R) \times H^{s - 2}(\R)}\\
\leq K(s) \varepsilon^3 t \Big( \|(N^0_\varepsilon, & W^0_\varepsilon) \|_{H^{s + 1}(\R) \times H^s(\R)} + \| (N^0_\varepsilon, W^0_\varepsilon) \|^2_{H^{s + 1}(\R) \times H^s(\R)} \Big),
\end{align*}
where $T_\varepsilon = \big( K(s) \varepsilon^3 \| (N_\varepsilon^0, W_\varepsilon^0) \|_{H^{s + 1}(\R) \times H^s(\R)} \big)^{- 1}$. In particular, when $ \| (N_\varepsilon^0, W_\varepsilon^0) \|_{H^{s + 1}(\R) \times H^s(\R)}$ remains uniformly bounded, then $T_\varepsilon = \boO \big( \varepsilon^{- 3} \big)$, and the wave equation is a good approximation for times of order $o \big( \varepsilon^{- 3} \big)$. The general solution to \eqref{eq:wamu} may be written as
$$(\n, \w) = (\n^+, \w^+) + (\n^-, \w^-),$$
where the functions $(\n^\pm, \w^\pm)$ are solutions to \eqref{eq:wamu} given by the d'Alembert formulae,
\begin{align*}
\big( \n^+(y, s), \w^+(y, s) \big) = \big( N^+(y - \sqrt{2} s), W^+(y - \sqrt{2} s) \big),\\
\big( \n^-(y, s), \w^-(y, s) \big) = \big( N^-(y + \sqrt{2} s), W^-(y + \sqrt{2} s) \big),
\end{align*}
where the profiles $N^\pm$ and $W^\pm$ are real-valued functions on $\R$. Solutions may therefore be split into right and left going waves of speed $\sqrt{2}$. Since the functions $(\n^\pm, \w^\pm)$ are solutions to \eqref{eq:wamu}, it follows that
$$\partial_y \big( N^+ + W^+ \big) = 0, \ {\rm and} \ \partial_y \big( N^- - W^- \big) = 0,$$
so that, if the functions decay to zero at infinity, then
$$N^\pm = \mp W^\pm = \frac{N_\varepsilon^0 \mp W_\varepsilon^0}{2}.$$

Theorem \ref{cochondore} extends our earlier results in \cite{BeGrSaS2} (see also \cite{ChirRou2} for an alternative approach and an extension to the higher dimensional case). It shows that the Korteweg-de Vries equation provides the appropriate approximation for time scales of order $\boO \big( \varepsilon^{- 3} \big)$. The definition of the new coordinate $x^+$, respectively $ x^-$, corresponds to a reference frame travelling to the left, respectively to the right, with speed $\sqrt{2}$ in the original coordinates $(\x, t)$. In the frame corresponding to $x^-$, the wave $(\n^-, \w^-)$, originally travelling to the left, is now stationary, whereas the right going wave $(\n^+, \w^+)$ now has a speed equal to $8 \varepsilon^{- 2}$. The coordinate $x^-$ is therefore particularly appropriate for the study of waves travelling to the left, whereas the coordinate $x^+$ is appropriate for the study of waves travelling to the right. In \cite{BeGrSaS2}, we imposed some additional assumptions which implied in particular the smallness of $U_\varepsilon^+$, so that it was only the study of waves going to the left which was addressed. This approach simplifies somehow the analysis.

In this paper, we remove this smallness assumption, at the cost however of new assumption \eqref{H1}, and we analyze both waves at the same time. Finally, we would like to emphasize also that comparing Theorem \ref{cochondore} with Theorem 1.4 in \cite{BeGrSaS2}, the error term now involves $\varepsilon^2$ instead of $\varepsilon$. As explained in \cite{BeGrSaS2}, the $\varepsilon^2$ is somewhat optimal, as the specific examples provided by travelling waves show. This improvement is related to the fact that we use higher order derivatives. As a matter of fact, the same improvement holds in the setting of Theorem 1.4 of \cite{BeGrSaS2}. More precisely, we have

\begin{theorem}
\label{cochon2}
Let $\varepsilon > 0$ and $k \geq 0$ be given. Assume that the initial data $\Psi_0$ belongs to $X^{k + 6}(\R)$ and satisfies
\begin{equation}
\label{grinzing2}
\| N_\varepsilon^0 \|_{H^{k + 5}(\R)} + \varepsilon \| \partial_x^{k + 6} N_\varepsilon^0 \|_{L^2(\R)} + \|\partial_x \Theta_\varepsilon^0\|_{H^{k + 5}(\R)} \leq K_0.
\end{equation}
Let $\boN^\pm$ and $\boW^\pm$ denote the solutions to the Korteweg-de Vries equation
\footnote{As well as the functions $\boU^\pm$, the functions $\boN^\pm$ and $\boW^\pm$ depend on $\varepsilon$. We again hide this dependence in the notations in order to stress the fact that the equations they satisfy are independent of $\varepsilon$.}
$$\partial_\tau \boU^\pm \mp \partial_x^3 \boU^\pm \mp \boU^\pm \partial_x \boU^\pm = 0,$$
with initial data $N_\varepsilon^0$, respectively $\partial_x \Theta_\varepsilon^0$. There exists positive constants $\varepsilon_2$ and $K_2$, depending possibly on $K_0$ and $k$, such that
\begin{equation}
\begin{split}
\label{ineq2}
\| \boN^\pm & (\cdot, \tau) - N_\varepsilon^\pm(\cdot, \tau ) \|_{H^k(\R)} + \| \boW^\pm(\cdot, \tau) - \partial_x \Theta_\varepsilon^\pm(\cdot, \tau ) \|_{H^k(\R)}\\
& \leq K_2 \big( \varepsilon^2 + \| N_\varepsilon^0 \pm \partial_x \Theta_\varepsilon^0 \|_{H^k(\R)} \big) \exp K_2 |\tau|,
\end{split}
\end{equation}
for any $\tau \in \R$, provided $\varepsilon \leq \varepsilon_2$.
\end{theorem}

\begin{remark}
If the term $\| N_\varepsilon^0 \pm \partial_x \Theta_\varepsilon^0 \|_{H^k(\R)}$ is small, then the Korteweg-de Vries approximation is valid on a time interval (in the original time variable) $t \in [0, T_\eps]$ with
$$T_\eps = o \bigg( \min \bigg\{ \frac{|\log(\eps)|}{\eps^3}, \frac{|\log(\|N_\eps^0 \pm \partial_x \Theta_\eps^0 \|_{H^k(\R)})|}{\eps^3} \bigg\} \bigg).$$
In particular, if $\| N_\eps^0 \pm \partial_x \Theta_\eps^0 \|_{H^k(\R)} \leq C \eps^\alpha$, with $\alpha > 0$, then the approximation is valid on a time interval $t \in [0, T_\eps']$ with $T_\eps' = o(\eps^{- 3} |\log(\eps)|)$. Moreover, if $\| N_\eps^0 \pm \partial_x \Theta_\eps^0 \|_{H^k(\R)}$ is of order $\boO(\eps^2)$, then the approximation error remains of order $\boO(\eps^2)$ on a time interval $t \in [0, T_\eps'']$ with $T_\eps'' = \boO(\eps^{- 3})$.
\end{remark}

The main difference between Theorems \ref{cochondore} and \ref{cochon2} is that the second one involves the functions $N_\varepsilon^\pm$ and $\partial_x \Theta_\varepsilon^\pm$ instead of $U_\eps^\pm$. The error term in \eqref{ineq2} involves the quantity $\| N_\varepsilon^0 \pm \partial_x \Theta_\varepsilon^0 \|_{H^k(\R)}$, which is small if the wave travelling in the other direction is small. Let us also emphasize that, in contrast with the assumptions of Theorem \ref{cochondore}, the assumptions of Theorem \ref{cochon2} do not involve any assumption on the $\boM$-norm.

Finally, it is worthwhile to notice that similar issues have been addressed and solved in the case of the long-wave limit of the water wave system. As a matter of fact, system \eqref{eq:dynaslow} bears some resemblance with a Boussinesq system. In a seminal work \cite{Craig1}, Craig proved the first rigorous convergence result towards the Korteweg-de Vries equation, under assumptions which are similar in spirit to the ones in Theorem \ref{cochon2}, focusing on a wave travelling in a single direction. Schneider and Wayne \cite{SchnWay1} completed the analysis and were able to handle both left and right-going waves at the same time, a result similar in spirit to Theorem \ref{cochondore}. Bona, Colin and Lannes provided a sharp error estimate in \cite{BonCoLa1} (see also \cite{Wright1}). The asymptotics were fully justified in \cite{AlSaLan1}, as well as in the higher dimensional case. In order to control the interactions between the two waves, sometimes called secular growth, these authors introduce additional assumptions on the initial data, of different nature but in the spirit, similar to our introduction of the $\boM$-norm which is more natural in our context.

We also emphasize that in the regime under study, the solutions of the Gross-Pitaevskii
equation have their modulus close to one, so that by the Madelung transform, one
is reduced to analyze a dispersive perturbation of a hyperbolic system. In this direction, Ben Youssef and Colin considered in \cite{BeYoCol1} similar limits for general hyperbolic systems perturbed by linear dispersive terms.

With respect to those works, the main difficulty regarding system \eqref{eq:dynaslow} is related to the fact that the dispersion terms are nonlinear. This difficulty is overcome using the integrability of the Gross-Pitaevskii and Korteweg-de Vries equations which allow to derive suitable bounds in high regularity spaces.

\subsection{Some elements in the proofs}

The proofs are somewhat parallel with the proofs in \cite{BeGrSaS2}, so that we will try to emphasize the new ideas and ingredients. Concerning Theorem \ref{cochondore}, the left and right going waves $U_\varepsilon^-$ and $U_\varepsilon^+$ play the same role in estimate \eqref{ineq1}, so that we may focus for instance on the estimates on $U_\varepsilon^-$. In order to simplify our notation, we set $N_\varepsilon = N_\varepsilon^-$, $\Theta_\varepsilon = \Theta_\varepsilon^-$, $U_\varepsilon = U_\varepsilon^-$, $\boU = \boU^-$, and
$$x = x^- = \varepsilon (\x + \sqrt{2} t).$$
We also introduce the new notation
\begin{equation}
\label{defV}
V_\varepsilon(x,\tau) = \frac{1}{2} \Big( N_\varepsilon(x, \tau) - \partial_x \Theta_\varepsilon(x, \tau) \Big) = U_\varepsilon^+ \Big( x - \frac{8}{\varepsilon^2}\tau, \tau \Big),
\end{equation}
and compute the relevant equations for $U_\varepsilon$ and $V_\varepsilon$,
\begin{equation}
\label{slow1}
\partial_\tau U_\varepsilon + \partial_x^3 U_\varepsilon + U_\varepsilon \partial_x U_\varepsilon = f_\varepsilon - \varepsilon^2 r_\varepsilon,
\end{equation}
where
\begin{equation}
\label{fF}
f_\varepsilon \equiv \partial_x \Big( \frac{1}{6} V_\varepsilon^2 - \partial_x^2 V_\varepsilon + \frac{1}{3} U_\varepsilon V_\varepsilon \Big) \equiv \partial_x F_\varepsilon,
\end{equation}
and
\begin{equation}
\label{slow2}
\partial_\tau V_\varepsilon + \frac{8}{\varepsilon^2} \partial_x V_\varepsilon = g_\varepsilon + \varepsilon^2 r_\varepsilon,
\end{equation}
where
\begin{equation}
\label{gG}
g_\varepsilon \equiv \partial_x \Big( \partial_x^2 N_\varepsilon + \frac{1}{2} V_\varepsilon^2 - \frac{1}{6} U_\varepsilon^2 - \frac{1}{3} U_\varepsilon V_\varepsilon \Big) \equiv \partial_x G_\varepsilon.
\end{equation}
The remainder term $r_\varepsilon$ is given by the formula
\begin{equation}
\label{rR}
r_\varepsilon = \partial_x \Big( \frac{N_\varepsilon \partial_x^2 N_\varepsilon}{6 (1 - \frac{\varepsilon^2}{6} N_\varepsilon)} + \frac{(\partial_x N_\varepsilon)^2 }{12 (1 - \frac{\varepsilon^2}{6} N_\varepsilon)^2} \Big) \equiv \partial_x R_\varepsilon.
\end{equation}
On the left-hand side of \eqref{slow1}, we recognize the \eqref{KdV} equation, whereas on the left-hand side of \eqref{slow2}, we recognize the transport operator with constant speed $- 8 \varepsilon^{- 2}$, which stems from the fact that we are working in moving frames with speed $\pm 4 \varepsilon^{- 2}$. It remains to establish that the terms on the right-hand side of \eqref{slow1} behave as error terms. The first step is to establish that all quantities are uniformly bounded on finite time intervals.

\begin{prop}
\label{Bounded}
Let $k \in \N$. Given any $\varepsilon > 0$ sufficiently small, assume that the initial data $\Psi_0(\cdot) = \Psi(\cdot, 0)$ belongs to $X^{k + 1}(\R)$ and satisfies the assumption
\begin{equation}
\label{prater}
\| N_\varepsilon^0 \|_{H^k(\R)} + \varepsilon \| \partial_x^{k + 1} N_\varepsilon^0 \|_{L^2(\R)}+ \|\partial_x \Theta_\varepsilon^0\|_{H^k(\R)} \leq K_0,
\end{equation}
where $K_0$ is some given positive constant. Then, there exists a positive constant $K$ depending only on $K_0$ and $k$, such that
\begin{equation}
\label{granderoue}
\| N_\varepsilon(\cdot, \tau) \|_{H^k(\R)} + \varepsilon \| \partial_x^{k + 1} N_\varepsilon(\cdot, \tau) \|_{L^2(\R)}+ \|\partial_x \Theta_\varepsilon(\cdot, \tau) \|_{H^k(\R)} \leq K \exp K |\tau|,
\end{equation}
for any $\tau \in \R$. In particular, we have
\begin{equation}
\label{frites}
\| U_\varepsilon(\cdot, \tau) \|_{H^k(\R)} + \| V_\varepsilon(\cdot, \tau) \|_{H^k(\R)} \leq K \exp K |\tau|.
\end{equation}
\end{prop}

\begin{remark}
Here and in the sequel, when we write $\eps$ sufficiently small, we mean that $0 < \eps \leq \eps_0$, where $\eps_0$ is some constant which depends only on $K_0$, but not on the order of differentiation $k$. In the course of our proofs, the constant $\eps_0$ is determined so that, when assumption \eqref{prater} holds, the energy of $\Psi$ is sufficiently small in order that \eqref{half} holds.
\end{remark}

It follows from Proposition \ref{Bounded} that the quantity $r_\varepsilon$ remains bounded on finite time intervals, so that the error term $\varepsilon^2 r_\varepsilon$ is of order $\varepsilon^2$ as desired in Theorem \ref{cochondore}. In contrast, the term $f_\varepsilon$ describes the interaction of the two waves and is therefore more delicate to handle. Since $V_\varepsilon$ is not supposed to be small in Theorem \ref{cochondore} (in contrast with Theorem \ref{cochon2}), $f_\varepsilon$ is not small in general. However, due to the dynamics, the interaction turns out to be of lower order. Indeed, in view of \eqref{slow2}, at leading order, the function $V_\varepsilon$ is shifted to the left with speed $8 \varepsilon^{- 2}$. Since the definition of $f_\varepsilon$ strongly depends on the function $V_\varepsilon$, a related property also holds for $f_\varepsilon$, so that the average interaction turns out to be small. To provide a rigorous justification of this last claim, we need however to localize the functions $U_\varepsilon$ and $V_\varepsilon$. This leads us to use the norm $\| \cdot \|_{\boM(\R)}$.

As in \cite{BeGrSaS2}, our proofs rely on energy methods. We introduce the difference $Z_\varepsilon \equiv U_\varepsilon - \boU$, which satisfies the equation
\begin{equation}
\label{dede}
\partial_\tau Z_\varepsilon + \partial_x^3 Z_\varepsilon + \boU \partial_x Z_\varepsilon + Z_\varepsilon \partial_x \boU + Z_\varepsilon \partial_x Z_\varepsilon = f_\varepsilon - \varepsilon^2 r_{\varepsilon}.
\end{equation}
In order to compute the $L^2$-norm of $\partial_x^k Z_\varepsilon$, we apply the differential operator $\partial_x^k$ to \eqref{dede}, multiply the resulting equation by $\partial_x^k Z_\varepsilon$ and integrate on $\R$ to obtain
\begin{align*}
\frac{1}{2} \partial_\tau \| \partial_x^k Z_\varepsilon \|_{L^2(\R)}^2 = & - \int_\R \partial_x^{k+1} \big( \boU Z_\varepsilon \big) \partial_x^k Z_\varepsilon - \frac{1}{2} \int_\R \partial_x^{k+1} \big( Z_\varepsilon^2 \big) \partial_x^k Z_\varepsilon\\
& + \int_\R \partial_x^k f_\varepsilon \partial_x^k Z_\varepsilon - \varepsilon^2 \int_\R \partial_x^k r_\varepsilon \partial_x^k Z_\varepsilon.
\end{align*}
Setting
\begin{equation}
\label{beauzire}
\boZ_\varepsilon^k(\tau) \equiv \int_0^\tau \int_\R \big( \partial_x^k Z_\varepsilon \big)^2,
\end{equation}
and integrating in time, we are led to the differential equation
\begin{equation}
\label{gege}
\begin{split}
\frac{1}{2} \partial_\tau \boZ_\varepsilon^k(\tau) = - & \int_0^\tau \int_\R \partial_x^{k+1} \big( \boU Z_\varepsilon \big) \partial_x^k Z_\varepsilon - \frac{1}{2} \int_0^\tau \int_\R \partial_x^{k+1} \big( Z_\varepsilon^2 \big) \partial_x^k Z_\varepsilon\\
+ & \int_0^\tau \int_\R \partial_x^k f_\varepsilon \partial_x^k Z_\varepsilon - \varepsilon^2 \int_0^\tau \int_\R \partial_x^k r_\varepsilon \partial_x^k Z_\varepsilon.
\end{split}
\end{equation}
The proof of Theorem \ref{cochondore} then follows applying the Gronwall lemma to \eqref{gege} provided we are first able to bound suitably all the terms on the right-hand side of \eqref{gege}. The first, second and fourth terms can be handled thanks to Proposition \ref{Bounded}. Indeed, we will show in Section \ref{quatre} that these terms can be bounded as follows
\begin{align}
\label{bound1}
&\bigg| \int_0^\tau \int_\R \partial_x^{k+1} \big( \boU Z_\varepsilon \big) \partial_x^k Z_\varepsilon \bigg| 
\leq K \bigg( \bigg| \int_0^\tau \int_\R \big( \partial_x^k Z_\varepsilon \big)^2 \bigg| + \bigg| \int_0^\tau \exp K |s| \ \| Z_\varepsilon(\cdot, s) \|_{H^k(\R)} ds \bigg| \bigg),\\ \label{bound2}
&\bigg| \int_0^\tau \int_\R \partial_x^{k+1} \big( Z_\varepsilon^2 \big) \partial_x^k Z_\varepsilon \bigg|
\leq K \bigg( \bigg| \int_0^\tau \int_\R \big( \partial_x^k Z_\varepsilon \big)^2 \bigg| + \bigg| \int_0^\tau \exp K |s| \ \| Z_\varepsilon(\cdot, s) \|_{H^k(\R)} ds \bigg| \bigg),\\\label{bound3}
&\bigg| \int_0^\tau \int_\R \partial_x^k r_\varepsilon \partial_x^k Z_\varepsilon \bigg| \leq K \bigg| \int_0^\tau \exp K |s| \ \| Z_\varepsilon(\cdot, s) \|_{H^k(\R)} ds \bigg|.
\end{align}
For the third term, we will prove

\begin{prop}
\label{Estim-f}
Let $\varepsilon > 0$ be given sufficiently small. Given any $k \in \N$, assume that the initial datum $\Psi_0(\cdot) = \Psi(\cdot, 0)$ belongs to $X^{k + 6}(\R)$ and satisfies \eqref{grinzing1} for some positive constant $K_0$. Then, there exists a positive constant $K$ depending only on $K_0$ and $k$, such that
\begin{equation}
\label{estim-f}
\begin{split}
& \bigg| \int_0^\tau \int_\R \partial_x^k f_\varepsilon \partial_x^k Z_\varepsilon \bigg|\\
\leq K \varepsilon^2 \bigg( \big( \varepsilon^2 + \| \partial_x^k Z_\varepsilon(\cdot, \tau) \|_{L^2(\R)} \big) & \exp K |\tau| + \bigg|\int_0^\tau \exp K |s| \ \| Z_\varepsilon(\cdot, s) \|_{H^k(\R)} ds \bigg| \bigg),
\end{split}
\end{equation}
for any $\tau \in \R$.
\end{prop}

Combining \eqref{gege} with bounds \eqref{bound1}, \eqref{bound2} and \eqref{estim-f}, we will obtain
$$\partial_\tau \boZ_\varepsilon^k(\tau) \leq K \Big( \sign(\tau) \boZ_\varepsilon^k(\tau) + \varepsilon^4 \exp K |\tau| \bigg).$$
The proof of Theorem \ref{cochondore} will then follow applying the Gronwall lemma.

We next say a few words about the proof of Proposition \ref{Estim-f}. For sake of clarity, we assume here that $k = 0$. For given $\tau \in \R$, we then have
\begin{equation}
\label{huitre}
\begin{split}
\int_0^\tau \int_\R f_\varepsilon Z_\varepsilon & = \int_0^\tau \int _\R \Big( \frac{1}{6} \partial_x V_\varepsilon^2 - \partial_x^3 V_\varepsilon + \frac{1}{3} \partial_x \big( U_\varepsilon V_\varepsilon \big) \Big) Z_\varepsilon\\
& = \int_0^\tau \int_\R \partial_x V_\varepsilon \Big( \frac{1}{3} V_\varepsilon Z_\varepsilon - \partial_x^2 Z_\varepsilon + \frac{1}{3} U_\varepsilon Z_\varepsilon \Big) + \frac{1}{3} \int_0^\tau \int_\R \partial_x U_\varepsilon V_\varepsilon Z_\varepsilon.
\end{split}
\end{equation}
For the first integral on the right-hand side, we use transport equation \eqref{slow2} and write
$$\partial_x V_\varepsilon = \frac{\varepsilon^2}{8} \Big( - \partial_\tau V_\varepsilon + g_\varepsilon + \varepsilon^2 r_\varepsilon \Big).$$
This change makes apparent an $\varepsilon^2$ factor, and then we continue the computation using integration by parts for the time and space variables as well as the bounds provided by Proposition \ref{Bounded}. It remains to handle the second integral on the right-hand side of \eqref{huitre}, which does not involve as before the spatial derivative $\partial_x V_\varepsilon$. We somewhat artificially introduce such a derivative considering an antiderivative $\Upsilon_\varepsilon$ of $V_\varepsilon$ defined by
\begin{equation}
\label{defupsilon}
\Upsilon_\varepsilon(x, \tau)=\int_{-R}^x V_\varepsilon(y, \tau) dy,
\end{equation}
where the positive number $R$ will be suitably chosen in the course of the computations. We obtain 
$$\int_0^\tau \int_\R \partial_x U_\varepsilon V_\varepsilon Z_\varepsilon = \int_0^\tau \int_\R \partial_x U_\varepsilon \partial_x \Upsilon_\varepsilon Z_\varepsilon.$$
The function $\Upsilon_\varepsilon$ satisfies a transport equation, namely
\begin{equation}
\label{equpsilon}
\partial_\tau \Upsilon_\varepsilon + \frac{8}{\varepsilon^2} \partial_x \Upsilon_\varepsilon = G_\varepsilon + \varepsilon^2 R_\varepsilon - G_\varepsilon(- R) - \varepsilon^2 R_\varepsilon(- R) + \frac{8}{\varepsilon^2} V_\varepsilon(- R),
\end{equation}
so that it is possible to implement a similar argument as above replacing $\partial_x \Upsilon_\varepsilon$ by the expression provided by \eqref{equpsilon}. In order to perform our computation, it turns out that we only need to obtain some control in time of the $L^\infty$-norm of $\Upsilon_\varepsilon$, which essentially amounts to controlling the $\boM$-norm of $V_\eps.$ At initial time, this corresponds to assumption \eqref{H1} on $N_\varepsilon^0$ and $\partial_x \Theta_\varepsilon^0$. For later time, we have

\begin{lemma}
\label{Transport}
Let $0 \leq E_0 < \frac{2 \sqrt{2}}{3}$ and $\Psi_0 \in \boM(\R) \cap X^4(\R)$ be given such that $E(\Psi_0) \leq E_0$. Then, there exists a positive constant $K$, depending only on $E_0$, such that
\begin{equation}
\label{transm}
\begin{split}
\| \eta(\cdot, t) \|_{\boM(\R)} + & \sqrt{2} \| \partial_\x \varphi(\cdot, t) \|_{\boM(\R)} \leq K \bigg( \| \eta^0 \|_{\boM(\R)} + \sqrt{2} \| \partial_\x \varphi^0 \|_{\boM(\R)}\\
+ & \bigg| \int_0^t \Big( \| \partial_\x^2 \eta(\cdot, s) \|_{L^\infty(\R)} + \| \eta(\cdot, s) \|_{W^{2, \infty}(\R)}^2 + \| \partial_\x \varphi(\cdot, s) \|_{L^\infty(\R)}^2 \Big) ds \bigg| \bigg),
\end{split}
\end{equation}
for any $t \in \R$.
\end{lemma}

In the slow coordinate $t$, Lemma \ref{Transport} provides a control on the norm $\| \Upsilon_\varepsilon \|_{L^\infty(\R)}$, which is independent of $\varepsilon$ and grows linearly in time.

\begin{lemma}
\label{Wanted}
Let $\varepsilon > 0$ be sufficiently small. Assume that the initial datum $\Psi_0(\cdot) = \Psi(\cdot, 0)$ belongs to $X^4(\R)$ and satisfies assumption \eqref{prater} for $k = 3$ and some positive constant $K_0$. Then, there exists a positive constant $K$, which does not depend on $\varepsilon$ nor $\tau$, such that
\begin{equation}
\label{massive}
\| N_\varepsilon(\cdot, \tau) \|_{\boM(\R)} + \| \partial_x \Theta_\varepsilon(\cdot, \tau) \|_{\boM(\R)} \leq K \Big( \| N_\varepsilon^0 \|_{\boM(\R)} + \| \partial_x \Theta_\varepsilon^0 \|_{\boM(\R)} + |\tau| \Big),
\end{equation}
for any $\tau \in \R$. In particular, we have
\begin{equation}
\label{attack}
\| U_\varepsilon(\cdot, \tau) \|_{\boM(\R)} + \| V_\varepsilon(\cdot, \tau) \|_{\boM(\R)} \leq K \Big( \| U_\varepsilon^0 \|_{\boM(\R)} + \| V_\varepsilon^0 \|_{\boM(\R)} + |\tau| \Big).
\end{equation}
\end{lemma}

The proof of estimate \eqref{estim-f} follows combining the bounds of Proposition \ref{Bounded} and Lemma \ref{Wanted}.

We next give some elements of the proof of Theorem \ref{cochon2}. Notice first that the functions $N_\varepsilon^\pm$ and $\partial_x \Theta_\varepsilon^\pm$ associated to the coordinates $x^+$ and $x^-$ play the same role in Theorem \ref{cochon2}, so that we may focus again for the estimates on $N_\varepsilon \equiv N_\varepsilon^-$ and $\partial_x \Theta_\varepsilon \equiv \partial_x \Theta_\varepsilon^-$. Recall that the main differences with respect to Theorem \ref{cochondore} are that Theorem \ref{cochon2} addresses the functions $N_\varepsilon$ and $\partial_x \Theta_\varepsilon$ instead of $U_\varepsilon$, and does not involve any assumption on $\boM$-norms. The proof of Theorem \ref{cochon2} is however parallel to the one of Theorem \ref{cochondore}. We write for the function $N_\varepsilon$,
\begin{equation}
\label{trianguler1}
\| N_\varepsilon - \boN \|_{H^k(\R)} \leq \| V_\varepsilon \|_{H^k(\R)} + \| U_\varepsilon - \boU \|_{H^k(\R)} + \| \boU - \boN \|_{H^k(\R)},
\end{equation}
since by definition, $V_\varepsilon = N_\varepsilon - U_\varepsilon$. Similarly, the functions $\partial_x \Theta_\varepsilon$ and $\boW$ satisfy
\begin{equation}
\label{trianguler2}
\| \partial_x \Theta_\varepsilon - \boW \|_{H^k(\R)} \leq \| V_\varepsilon \|_{H^k(\R)} + \| U_\varepsilon - \boU \|_{H^k(\R)} + \| \boU - \boW \|_{H^k(\R)}, 
\end{equation}
so that the proof of \eqref{ineq2} reduces to bound each of the terms on the right-hand side of \eqref{trianguler1} and \eqref{trianguler2}. For the $H^k$-norm of $V_\varepsilon$, we invoke again energy estimates, based now on equation \eqref{slow2}. This yields

\begin{prop}
\label{Estim-V}
Let $\varepsilon > 0$ be sufficiently small. Given any $k \in \N$, assume that the initial datum $\Psi_0(\cdot) = \Psi(\cdot, 0)$ belongs to $X^{k + 6}(\R)$ and satisfies \eqref{grinzing2} for some positive constant $K_0$. Then, there exists a positive constant $K$ depending only on $K_0$ and $k$, such that
\begin{equation}
\label{estim-V}
\| V_\varepsilon(\cdot, \tau) \|_{H^k(\R)} \leq K \big( \| V_\varepsilon^0 \|_{H^k(\R)} + \varepsilon^2 \big) \exp K |\tau|,
\end{equation}
for any $\tau \in \R$.
\end{prop}

Similarly, concerning the differences between the solutions to \eqref{KdV}, $\boU$ and $\boN$ on one hand, and $\boU$ and $\boW$ on the other, we invoke a general stability result for the Korteweg-de Vries equation, which we recall for the sake of completeness. The proof follows from standard $H^k$-energy methods for the difference, using the conserved quantities of \eqref{KdV} in order to bound uniformly the $H^{k+1}$-norms coming from the quadratic terms.

\begin{lemma}
\label{Initier}
Let $k \in \N$ be given. Consider two functions $F^0$ and $G^0$ in $H^{k + 2}(\R)$ and denote $F$ and $G$ the solutions to \eqref{KdV} with initial data $F^0$, respectively $G^0$. Then, there exists a constant $K$, depending only on $k$ and the $H^{k + 2}$-norms of $F^0$ and $G^0$, such that
\begin{equation}
\label{init}
\| F(\cdot, \tau) - G(\cdot, \tau) \|_{H^k(\R)} \leq K \| F^0 - G^0 \|_{H^k(\R)} \exp K |\tau|,
\end{equation}
for any $\tau \in \R$.
\end{lemma}

In view of Lemma \ref{Initier}, and the fact that
$$\boU^0 - \boW^0 = \boN^0 - \boU^0 = V_\varepsilon^0,$$
we are led to
\begin{equation}
\label{initfin}
\| \boU(\cdot, \tau) - \boN(\cdot, \tau) \|_{H^k(\R)} + \| \boU(\cdot, \tau) - \boW(\cdot, \tau) \|_{H^k(\R)} \leq K \| V_\varepsilon^0 \|_{H^k(\R)} \exp K |\tau|,
\end{equation}
for any $\tau \in \R$. Going back to \eqref{trianguler2}, it remains to estimate the term $\| U_\varepsilon - \boU \|_{H^k(\R)} = \| Z_\varepsilon \|_{H^k(\R)}$. An upper bound for this term was given in Theorem \ref{cochondore} using bounds on the $\boM$-norm of $N_\eps$ and $\partial_x \Theta_\eps.$ Here, we use instead the estimates of Proposition \ref{Estim-V} to bound the interaction terms.

\begin{prop}
\label{Estim-f2}
Let $\varepsilon > 0$ be given sufficiently small. Given any $k \in \N$, assume that the initial datum $\Psi_0(\cdot) = \Psi(\cdot, 0)$ belongs to $X^{k + 6}(\R)$ and satisfies \eqref{grinzing2} for some positive constant $K_0$. Then, there exists a positive constant $K$ depending only on $K_0$ and $k$, such that
\begin{equation}
\label{estim-f2}
\begin{split}
\bigg| \int_0^\tau \int_\R \partial_x^k f_\varepsilon \partial_x^k Z_\varepsilon \bigg| \leq & K \varepsilon^2 \Big( \varepsilon^2 + \| V_\varepsilon^0 \|_{H^k(\R)} + \| \partial_x^k Z_\varepsilon(\cdot, \tau) \|_{L^2(\R)} \Big) \exp K |\tau|\\
+ & K \Big( \varepsilon^2 + \| V_\varepsilon^0 \|_{H^k(\R)} \Big) \bigg| \int_0^\tau \exp K |s| \ \| Z_\varepsilon(\cdot, s) \|_{H^k(\R)} ds \bigg|,
\end{split}
\end{equation}
for any $\tau \in \R$.
\end{prop}

We then adapt the arguments of the proof of \eqref{ineq1} in order to obtain that
\begin{equation}
\label{ineq12}
\| U_\varepsilon(\cdot, \tau) - \boU(\cdot, \tau) \|_{H^k(\R)} \leq K \big( \varepsilon^2 + \| V_\varepsilon^0 \|_{H^k(\R)} \big) \exp K |\tau|,
\end{equation}
for any $\tau \in \R$. Combining with inequalities \eqref{trianguler1} and \eqref{trianguler2}, and bounds \eqref{estim-V} and \eqref{initfin}, this will complete the proof of Theorem \ref{cochon2}.

\subsection{Outline of the paper}

This paper is organized as follows. In the next section, we provide the proofs to Proposition \ref{Bounded}, and Lemmas \ref{Transport} and \ref{Wanted}. In Section \ref{trois}, we prove Proposition \ref{Estim-f}. The proof of Theorem \ref{cochondore} is completed in Section \ref{quatre}, whereas we derive Proposition \ref{Estim-V}, Lemma \ref{Initier}, Proposition \ref{Estim-f2}, and finally Theorem \ref{cochon2} in Section \ref{cinq}. In a separate appendix, we extend the arguments of the proof of Lemma \ref{Transport} to provide a rigorous framework for the notion of mass and establish its conservation.

\numberwithin{cor}{section}
\numberwithin{equation}{section}
\numberwithin{lemma}{section}
\numberwithin{prop}{section}
\numberwithin{remark}{section}
\numberwithin{theorem}{section}
\section{Bounds for the rescaled functions}
\label{sobolev}

In this section, we establish a certain number of bounds for the rescaled functions $N_\varepsilon$, $\partial_x \Theta_\varepsilon$, $U_\varepsilon$ and $V_\varepsilon$, which are useful in the course of the proofs of Theorems \ref{cochondore} and \ref{cochon2}. We first derive the Sobolev bounds of Proposition \ref{Bounded}, and then we compute estimate \eqref{massive} of Lemma \ref{Wanted}.

\subsection{Sobolev bounds}
\label{Sobounds}

Two different arguments are under hand to derive the Sobolev bounds given by inequality \eqref{granderoue}. The first one relies on the integrability properties of \eqref{GP}. It is proved in \cite{BeGrSaS2} that the quantities $E_1(\Psi) \equiv E(\Psi)$,
$$E_2(\Psi) \equiv \frac{1}{2} \int_{\R} |\partial_\x^2 \Psi|^2 - \frac{3}{2} \int_{\R} \eta |\partial_\x \Psi|^2 + \frac{1}{4} \int_{\R} (\partial_\x \eta)^2 - \frac{1}{4} \int_\R \eta^3,$$
\begin{align*}
E_3(\Psi) \equiv & \frac{1}{2} \int_{\R} |\partial_\x^3 \Psi|^2 + \frac{1}{4} \int_{\R} |\partial_\x^2 \eta|^2 + \frac{5}{4} \int_{\R} |\partial_\x \Psi|^4 + \frac{5}{2} \int_{\R} \partial_\x^2 \eta |\partial_\x \Psi|^2 - \frac{5}{2} \int_{\R} \eta |\partial_\x^2 \Psi|^2\\
& - \frac{5}{4} \int_\R \eta (\partial_\x \eta)^2 + \frac{15}{4} \int_\R \eta^2 |\partial_\x \Psi|^2 + \frac{5}{16} \int_\R \eta^4,
\end{align*}
and
\begin{align*}
E_4(\Psi) & \equiv \frac{1}{2} \int_{\R} |\partial_\x^4 \Psi|^2 + \frac{1}{4} \int_{\R} |\partial_\x^3 \eta|^2 - \frac{7}{4} \int_\R \eta (\partial_\x^2 \eta)^2 - \frac{7}{2} \int_\R \eta |\partial_\x^3 \Psi|^2 + \frac{35}{8} \int_{\R} \eta^2 (\partial_\x \eta)^2\\
+ & \frac{35}{4} \int_\R \eta^2 |\partial_\x^2 \Psi|^2 - \frac{35}{4} \int_{\R} (\partial_\x \eta)^2 |\partial_\x \Psi|^2 - \frac{7}{2} \int_\R |\partial_\x \Psi|^2 |\partial_\x^2 \Psi|^2 - 7 \int_\R \partial_\x^2 \eta \langle \partial_\x \Psi, \partial_\x^3 \Psi\rangle\\
- & 7 \int_\R |\partial_\x \Psi|^2 \langle \partial_\x \Psi, \partial_\x^3 \Psi \rangle - \frac{35}{2} \int_\R \eta \partial_\x^2 \eta |\partial_\x \Psi|^2 - \frac{35}{4} \int_\R \eta^3 |\partial_\x \Psi|^2 - \frac{35}{4} \int_\R \eta |\partial_\x \Psi|^4 - \frac{7}{16} \int_\R \eta^5,
\end{align*}
are conserved along the Gross-Pitaevskii flow, provided that the initial datum $\Psi_0$ belongs to $X^4(\R)$. Notice that, for $1\leq k\leq 4$, the quantities $E_k$ defined above give a control on the $L^2$-norms of the functions $\partial_\x^k \Psi$ and $\partial_\x^{k - 1} \eta$. As a matter of fact, invoking the Sobolev embedding theorem, one can establish that there exists some universal constant $K$ such that
$$E_k(\psi) \leq K \Big( \| \psi \|_{H^k(\R)} + \| 1 - |\psi|^2 \|_{H^{k - 1}(\R)} \Big)^{k + 1},$$
for any function $\psi \in X^k(\R)$. Similarly, there exists a positive constant $K(E_1(\psi), \ldots, E_{k - 1}(\psi))$, depending only on the quantities $E_1(\psi)$, $\ldots$ and $E_{k - 1}(\psi)$, such that
$$\| \psi \|_{H^k(\R)}^2 + \| 1 - |\psi|^2 \|_{H^{k - 1}(\R)}^2 \leq K(E_1(\psi), \ldots, E_{k - 1}(\psi)) \ E_k(\psi).$$
In particular, the conservation of the quantities $E_k(\Psi)$ along the Gross-Pitaevskii flow provides bounds on the Sobolev norms of the functions $\Psi$ and $\eta$, which only depends on the $X^k$-norms of the initial datum $\Psi_0$. In the rescaled variables, we obtain

\begin{prop}[\cite{BeGrSaS2}]
\label{Boundthemall}
Let $0 \leq k \leq 3$ and $\varepsilon > 0$ be given sufficiently small. Assume that the initial datum $\Psi_0(\cdot) = \Psi(\cdot, 0)$ belongs to $X^{k + 1}(\R)$ and satisfies assumption \eqref{prater} for some positive constant $K_0$. Then, there exists a positive constant $K$, which does not depend on $\varepsilon$ nor $\tau$, such that
\begin{equation}
\label{praterbis}
\| N_\varepsilon(\cdot, \tau) \|_{H^k(\R)} + \varepsilon \| \partial_x^{k + 1} N_\varepsilon(\cdot, \tau) \|_{L^2(\R)} + \|\partial_x \Theta_\varepsilon(\cdot, \tau) \|_{H^k(\R)} \leq K,
\end{equation}
for any $\tau \in \R$.
\end{prop}

Inequality \eqref{praterbis} presents the advantage to be uniform in time. We will invoke this property to derive Lemma \ref{Wanted}. We believe that inequality \eqref{praterbis} remains valid for higher order Sobolev spaces. However, it seems rather involved to prove this claim since this requires to compute, or at least to describe precisely, the higher order invariants of \eqref{GP} (see \cite{BeGrSaS2} for more details).

In order to compute higher Sobolev bounds, we rely on the energy estimates derived in \cite{BetDaSm1} in the context of the wave limit of the Gross-Pitaevskii equation mentioned above in the introduction. More precisely, given any $k \in \N$ and any positive number $E_0 < \frac{2 \sqrt{2}}{3}$, we consider some initial datum $\Psi_0 \in X^k(\R)$ such that $E(\Psi_0) \leq E_0$. Then, there exists a positive constant $m_0$, depending only on $E_0$, such that
\begin{equation}
\label{half}
m_0 \leq |\Psi(\x, t)| \leq \frac{1}{m_0},
\end{equation}
for any $(\x, t) \in \R^2$ (see \cite{BetGrSa2} for more details). Under this additional assumption, the computation of energy estimates for system \eqref{eq:dynaslow} achieved in Proposition 1 of \cite{BetDaSm1}, provides the tame estimates
\begin{equation}
\label{hippocampe}
\partial_t \Gamma_\varepsilon^k(t) \leq K(k, m_0) \varepsilon^2 \big( 1 + \varepsilon^2 \| n_\varepsilon(\cdot, t) \|_{L^\infty} \big) \big( \| \partial_x n_\varepsilon(\cdot, t) \|_{L^\infty} + \| \partial_x y_\varepsilon(\cdot, t) \|_{L^\infty} \big) \Big( \Gamma_\varepsilon^k(t) + \frac{\Gamma_\varepsilon^0(t)}{8} \Big),
\end{equation}
where $K(k, m_0)$ is some positive constant which does not depend on $\varepsilon$ nor $t$, and where the function $y_\varepsilon$ is defined by
$$y_\varepsilon = w_\varepsilon + \frac{i \varepsilon}{\sqrt{2}} \frac{\partial_x n_\varepsilon}{m_\varepsilon}, \ {\rm where} \ m_\varepsilon = 1 - \frac{\varepsilon^2}{6} n_\varepsilon,$$
while the notation $\Gamma_\varepsilon^k(t)$ refers to the functional
\begin{equation}
\label{gammvert}
\Gamma_\varepsilon^k(t) \equiv \| \partial_x^k n_\varepsilon(\cdot, t) \|_{L^2(\R)}^2 + \| m_\varepsilon(\cdot, t)^\frac{1}{2} \partial_x^k y_\varepsilon(\cdot, t) \|_{L^2(\R)}^2.
\end{equation}
In view of \eqref{half}, and since we have
$$m_\varepsilon(x, t) = \Big| \Psi \Big( \frac{x}{\varepsilon}, \frac{t}{\varepsilon} \Big) \Big|,$$
the quantity $\Gamma_\varepsilon^k$ controls the $H^k$-norms of $n_\varepsilon$ and $y_\varepsilon$. Going back to the original setting, we have in particular,
$$\Gamma_\varepsilon^0(t) = \frac{144}{\varepsilon^3} E \Big( \Psi \Big( \cdot, \frac{t}{\varepsilon} \Big) \Big),$$
so that, by the conservation of the energy,
\begin{equation}
\label{gammaun}
\Gamma_\varepsilon^0(t) = \frac{144}{\varepsilon^3} E(\Psi_0).
\end{equation}

\begin{proof}[Proof of Proposition \ref{Bounded}]
In the case $k < 4$, Proposition \ref{Bounded} is a direct consequence of Proposition \ref{Boundthemall}. In the case $k \geq 4$, the proof of Proposition \ref{Bounded} is obtained applying the Gronwall lemma to inequality \eqref{hippocampe}, using identity \eqref{gammaun} and bounds \eqref{praterbis}. We conclude rescaling the derived inequality in the variables $N_\varepsilon$ and $\partial_x \Theta_\varepsilon$.

More precisely, invoking assumption \eqref{prater} for $k = 0$, we first compute in the rescaled setting,
\begin{equation}
\label{petitapetit}
E(\Psi_0) = \frac{\varepsilon^3}{144} \int_\R \bigg( M_\varepsilon^0 (\partial_x \Theta_\varepsilon^0)^2 + (N_\varepsilon^0)^2 + \frac{\varepsilon^2 (\partial_x N_\varepsilon^0)^2}{2 M_\varepsilon^0} \bigg) \leq K \varepsilon^3,
\end{equation}
where $M_\varepsilon^0 \equiv 1 - \frac{\varepsilon^2}{6} N_\varepsilon^0$ and $K$ is a positive constant depending only on $K_0$. In particular, given any $\varepsilon$ sufficiently small, assumption \eqref{half} holds for $m_0=\frac{1}{2}$, so that in view of \eqref{hippocampe} and \eqref{gammaun}, we are led to
\begin{equation}
\label{gammamia}
\Gamma_\varepsilon^k(t) \leq \Gamma_\varepsilon^k(0) \exp \bigg| \int_0^t A_\varepsilon(s) ds \bigg| + \frac{18}{\varepsilon^3} E(\Psi_0) \bigg( \exp \bigg| \int_0^t A_\varepsilon(s) ds \bigg| - 1 \bigg),
\end{equation}
where
\begin{equation}
\label{ahahah}
A_\varepsilon(t) \equiv K(k, m_0) \varepsilon^2 \big( 1 + \varepsilon \| n_\varepsilon(\cdot, t) \|_{L^\infty(\R)} \big) \big( \| \partial_x n_\varepsilon(\cdot, t) \|_{L^\infty(\R)} + \| \partial_x y_\varepsilon(\cdot, t) \|_{L^\infty(\R)} \big).
\end{equation}
We now rescale inequality \eqref{gammamia} in the variables $N_\varepsilon$ and $\Theta_\varepsilon$ using the fact that
$$N_\varepsilon(x, \tau) = n_\varepsilon \Big( x - \frac{4 \tau}{\varepsilon^2}, \frac{2 \sqrt{2}}{\varepsilon^2} \tau \Big), \ {\rm and} \ \partial_x \Theta_\varepsilon(x, \tau) = w_\varepsilon \Big( x - \frac{4 \tau}{\varepsilon^2}, \frac{2 \sqrt{2}}{\varepsilon^2} \tau \Big).$$
In this scaling, definition \eqref{gammvert} may be written as
\begin{equation}
\label{gammbleu}
\upgamma_\varepsilon^k(\tau) \equiv \Gamma_\varepsilon^k \Big( \frac{2 \sqrt{2} \tau}{\varepsilon^2} \Big) = \| \partial_x^k N_\varepsilon(\cdot, \tau) \|_{L^2(\R)}^2 + \| M_\varepsilon(\cdot, \tau)^\frac{1}{2} \partial_x^k Y_\varepsilon(\cdot, \tau) \|_{L^2(\R)}^2,
\end{equation}
where $Y_\varepsilon \equiv \partial_x \Theta_\varepsilon + \frac{i \varepsilon}{\sqrt{2}} \frac{\partial_x N_\varepsilon}{M_\varepsilon}$ and $M_\varepsilon \equiv 1 - \frac{\varepsilon^2}{6} N_\varepsilon$. Similarly, definition \eqref{ahahah} may be rescaled as
$$\boA_\varepsilon(\cdot, \tau) \equiv \frac{2 \sqrt{2}}{\varepsilon^2} A_\varepsilon \Big( \frac{2 \sqrt{2} \tau}{\varepsilon^2} \Big) = K(k, m_0) \big( 1 + \varepsilon \| N_\varepsilon(\cdot, \tau) \|_{L^\infty} \big) \big( \| \partial_x N_\varepsilon(\cdot, \tau) \|_{L^\infty} + \| \partial_x Y_\varepsilon(\cdot, \tau) \|_{L^\infty} \big),$$
so that inequality \eqref{gammamia} becomes
$$\upgamma_\varepsilon^k(\tau) \leq \upgamma_\varepsilon^k(0) \exp \bigg| \int_0^\tau \boA_\varepsilon(s) ds \bigg| + \frac{18}{\varepsilon^3} E(\Psi_0) \bigg( \exp \bigg| \int_0^\tau \boA_\varepsilon(s) ds \bigg| - 1 \bigg).$$
Invoking assumption \eqref{prater} for $k = 2$, and the Sobolev embedding theorem, Proposition \ref{Boundthemall} provides
$$|\boA_\varepsilon(\cdot, \tau)| \leq K,$$
for any $\tau \in \R$, where $K$ is some positive constant depending on $k$ and $K_0$. Hence, by \eqref{petitapetit},
$$\upgamma_\varepsilon^k(\tau) \leq \big( \upgamma_\varepsilon^k(0) + K \big) \exp K |\tau|.$$
The proof of \eqref{granderoue} then follows from definition \eqref{gammbleu} and inequalities \eqref{half}. Inequality \eqref{frites} is a direct consequence of definitions \eqref{eq:uv} and \eqref{defV}.
\end{proof}

\subsection{Bounds in the space $\boM(\R)$}
\label{Mbounds}

We turn first to the proof of Lemma \ref{Wanted}. As mentioned in the introduction, inequality \eqref{massive} is a rescaled version of inequality \eqref{transm}, using the Sobolev estimates of Proposition \ref{Boundthemall} to bound the integral on the right-hand side of \eqref{transm}.

\begin{proof}[Proof of Lemma \ref{Wanted}]
In view of definitions \eqref{slow-var}, inequality \eqref{transm} may be recast in the variables $N_\varepsilon$ and $\partial_x \Theta_\varepsilon$ as
\begin{align*}
\| N_\varepsilon(\cdot, \tau) \|_{\boM(\R)} + & \| \partial_x \Theta_\varepsilon(\cdot, \tau) \|_{\boM(\R)} \leq \| N_\varepsilon^0 \|_{\boM(\R)} + \| \partial_x \Theta_\varepsilon^0 \|_{\boM(\R)}\\
+ & K \bigg| \int_0^\tau \Big( \| \partial_x^2 N_\varepsilon(\cdot, s) \|_{L^\infty(\R)} + \| N_\varepsilon(\cdot, s) \|_{W^{2, \infty}(\R)}^2 + \| \partial_x \Theta_\varepsilon(\cdot, s) \|_{L^\infty(\R)}^2 \Big) ds \bigg|,
\end{align*}
so that by the Sobolev embedding theorem,
\begin{equation}
\label{thiercelin}
\begin{split}
\| N_\varepsilon(\cdot, \tau) \|_{\boM(\R)} + & \| \partial_x \Theta_\varepsilon(\cdot, \tau) \|_{\boM(\R)} \leq \| N_\varepsilon^0 \|_{\boM(\R)} + \| \partial_x \Theta_\varepsilon^0 \|_{\boM(\R)}\\
+ & K \bigg| \int_0^\tau \Big( \| \partial_x^2 N_\varepsilon(\cdot, s) \|_{H^1(\R)} + \| N_\varepsilon(\cdot, s) \|_{H^3(\R)}^2 + \| \partial_x \Theta_\varepsilon(\cdot, s) \|_{H^1(\R)}^2 \Big) ds \bigg|.
\end{split}
\end{equation}
Invoking assumption \eqref{prater} for $k = 3$, we deduce from Proposition \ref{Boundthemall} that the integrand on the right-hand side of inequality \eqref{thiercelin} is bounded by some constant depending only on $K_0$. This completes the proof of \eqref{massive}. Inequality \eqref{attack} is then a direct consequence of definitions \eqref{eq:uv} and \eqref{defV}.
\end{proof}

We next prove Lemma \ref{Transport}.

\begin{proof}[Proof of Lemma \ref{Transport}]
The proof of Lemma \ref{Transport} relies on the conservative form of the system of equations satisfied by $\eta$ and $\partial_\x \varphi$ which may be written as
\begin{equation}
\label{sysetaphi}
\left\{ \begin{array}{ll}
\partial_t \eta - 2 \partial_\x \big( \partial_\x \varphi \big) = - 2 \partial_\x \big( \eta \partial_\x \varphi),\\
\partial_t \big( \partial_\x \varphi \big) - \partial_\x \eta = - \frac{1}{2} \partial_\x^3 \eta - \partial_\x \Big( |\partial_\x \varphi|^2 + \frac{\eta \partial_\x^2 \eta}{2 ( 1 -\eta)} + \frac{(\partial_\x \eta)^2}{4(1 - \eta)^2} \Big).
\end{array} \right.
\end{equation}
As already mentioned in the introduction, we recognize on the left-hand side of \eqref{sysetaphi}, a transport operator $\boT$, given by
$$\boT = (\partial_t - 2 \partial_\x, \partial_t - \partial_\x).$$
Introducing the variables
\begin{equation}
\label{duvet}
\begin{split}
u = \frac{1}{2} \Big( \eta + \sqrt{2} \partial_\x \varphi \Big), \ {\rm and} \ v = \frac{1}{2} \Big( \eta - \sqrt{2} \partial_\x \varphi \Big),
\end{split}
\end{equation}
we diagonalize the transport operator $\boT$, so that \eqref{sysetaphi} becomes
\begin{equation}
\label{sysuvet}
\left\{ \begin{array}{ll}
\partial_t u - \sqrt{2} \partial_\x u = - \partial_\x \big( \eta \partial_\x \varphi) - \frac{\sqrt{2}}{4} \partial_\x^3 \eta - \sqrt{2} \partial_\x \Big( \frac{1}{2} |\partial_\x \varphi|^2 + \frac{\eta \partial_\x^2 \eta}{4 ( 1 -\eta)} + \frac{(\partial_\x \eta)^2}{8(1 - \eta)^2} \Big),\\
\partial_t v + \sqrt{2} \partial_\x v = - \partial_\x \big( \eta \partial_\x \varphi) + \frac{\sqrt{2}}{4} \partial_\x^3 \eta + \sqrt{2} \partial_\x \Big( \frac{1}{2} |\partial_\x \varphi|^2 + \frac{\eta \partial_\x^2 \eta}{4 ( 1 -\eta)} + \frac{(\partial_\x \eta)^2}{8(1 - \eta)^2} \Big).
\end{array} \right.
\end{equation}
Given any real numbers $a < b$, we deduce from the first equation of \eqref{sysuvet} that
$$\partial_t \bigg( \int_a^b u \big( x -\sqrt{2} t, t \big) dx \bigg) = \bigg[ - \eta \partial_\x \varphi - \frac{\sqrt{2}}{4} \partial_\x^2 \eta - \sqrt{2} \Big( \frac{1}{2} |\partial_\x \varphi|^2 + \frac{\eta \partial_\x^2 \eta}{4 ( 1 -\eta)} + \frac{(\partial_\x \eta)^2}{8(1 - \eta)^2} \Big) \bigg]_{a - \sqrt{2} t}^{b - \sqrt{2} t}.$$
At this stage, it is worthwhile to recall that, since $E(\Psi_0) \leq E_0$, inequalities \eqref{half} hold for some positive number $m_0$, depending only on $E_0$,
so that
$$\bigg| \partial_t \bigg( \int_a^b u \big( x -\sqrt{2} t, t \big) dx \bigg) \bigg| \leq K \Big( \| \partial_\x^2 \eta(\cdot, t) \|_{L^\infty(\R)} + \| \eta(\cdot, t) \|_{W^{2, \infty}(\R)}^2 + \| \partial_\x \varphi(\cdot, t) \|_{L^\infty(\R)}^2 \Big),$$
where $K$ is some positive constant, depending only on $E_0$. Performing an integration in time, this may be recast as
\begin{equation}
\label{riou}
\begin{split}
\bigg| \int_a^b & u \big( x -\sqrt{2} t, t \big) dx \bigg| \leq \bigg| \int_a^b u \big( x, 0 \big) dx \bigg|\\
& + K \bigg| \int_0^t \Big( \| \partial_\x^2 \eta(\cdot, s) \|_{L^\infty(\R)} + \| \eta(\cdot, s) \|_{W^{2, \infty}(\R)}^2 + \| \partial_\x \varphi(\cdot, s) \|_{L^\infty(\R)}^2 \Big) ds \bigg|.
\end{split}
\end{equation}
Applying the change of coordinates $y = x - \sqrt{2} t$ to the left-hand side of \eqref{riou}, and invoking definition \eqref{massnorm}, we are led to
$$\| u( \cdot, t) \|_{\boM(\R)} \leq \| u( \cdot, 0) \|_{\boM(\R)} + K \bigg| \int_0^t \Big( \| \partial_\x^2 \eta(\cdot, s) \|_{L^\infty(\R)} + \| \eta(\cdot, s) \|_{W^{2, \infty}(\R)}^2 + \| \partial_\x \varphi(\cdot, s) \|_{L^\infty(\R)}^2 \Big) ds \bigg|.$$
The same proof applies to the second equation in \eqref{sysuvet} and this provides the same inequality for the function $v$. Combining with definitions \eqref{duvet}, inequality \eqref{transm} follows.
\end{proof}

\section{Estimates for the interaction terms}
\label{trois}

This section is devoted to the proof of Proposition \ref{Estim-f} which provides estimates of the interaction term $\int_0^\tau \int_\R \partial_x^k f_\varepsilon \partial_x^k Z_\varepsilon$ in \eqref{gege}.

\subsection{Proof of Proposition \ref{Estim-f}}

Given any integer $k$ and any number $\tau$, the interaction term $\int_0^\tau \int_\R \partial_x^k f_\varepsilon \partial_x^k Z_\varepsilon$ may be written in view of definition \eqref{fF},
$$\int_0^\tau \int_\R \partial_x^k f_\varepsilon \partial_x^k Z_\varepsilon = \int_0^\tau \int_\R \partial_x^{k+1} \Big( \frac{1}{6} V_\varepsilon^2 - \partial_x^2 V_\varepsilon + \frac{1}{3} U_\varepsilon V_\varepsilon \Big) \partial_x^k Z_\varepsilon,$$
so that by the Leibniz formula,
\begin{equation}
\label{mandanda}
\int_0^\tau \int_\R \partial_x^k f_\varepsilon \partial_x^k Z_\varepsilon = I_\varepsilon(\tau) + J_\varepsilon(\tau),
\end{equation}
where we denote
$$I_\varepsilon(\tau) \equiv \int_0^\tau \int_\R \bigg( \frac{1}{6} \partial_x^{k + 1} \big( V_\varepsilon^2 \big) - \partial_x^{k +3} V_\varepsilon + \frac{1}{3} \sum_{j = 0}^k \binom{k + 1}{j} \partial_x^j U_\varepsilon \partial_x^{k + 1 - j} V_\varepsilon \bigg) \partial_x^k Z_\varepsilon,$$
and
\begin{equation}
\label{zubar}
J_\varepsilon(\tau) \equiv \frac{1}{3} \int_0^\tau \int_\R \partial_x^{k + 1} U_\varepsilon V_\varepsilon \partial_x^k Z_\varepsilon.
\end{equation}
In view of \eqref{mandanda}, the proof of Proposition \ref{Estim-f} reduces to bound each term in the integrals $I_\varepsilon(\tau)$ and $J_\varepsilon(\tau)$ combining the estimates of Proposition \ref{Bounded} and Lemma \ref{Wanted} with some H\"older inequality and Sobolev embedding theorems. In particular, we will repetitively invoke the following bounds of the nonlinear functions $f_\varepsilon$, $g_\varepsilon$, $r_\varepsilon$, $F_\varepsilon$, $G_\varepsilon$ and $R_\varepsilon$, as well as the following estimates of the time derivative $\partial_\tau U_\varepsilon$, and of the solution $\boU$ to \eqref{KdV} with initial datum $U_\varepsilon^0$. We state these bounds in a series of lemmas whose proofs are each the object of a separate subsection here after the completion of the proof of Proposition \ref{Estim-f}.

\begin{lemma}
\label{Nonbounded}
Let $\varepsilon > 0$ be given sufficiently small. Given any $k \in \N$, assume that the initial datum $\Psi_0(\cdot) = \Psi(\cdot, 0)$ belongs to $X^{k + 6}(\R)$ and satisfies \eqref{grinzing2} for some positive constant $K_0$. Then, there exists a positive constant $K$ depending only on $K_0$ and $k$, such that
\begin{equation}
\label{allfolks1}
\begin{split}
& \| \partial_\tau U_\varepsilon(\cdot, \tau) \|_{H^{k + 2}(\R)} + \| f_\varepsilon(\cdot, \tau) \|_{H^{k + 2}(\R)} + \| g_\varepsilon(\cdot, \tau) \|_{H^{k + 2}(\R)} + \| r_\varepsilon(\cdot, \tau) \|_{H^{k + 2}(\R)}\\
& + \| F_\varepsilon(\cdot, \tau) \|_{H^{k + 3}(\R)} + \| G_\varepsilon(\cdot, \tau) \|_{H^{k + 3}(\R)} + \| R_\varepsilon(\cdot, \tau) \|_{H^{k + 3}(\R)} \leq K \exp K |\tau|,
\end{split}
\end{equation}
for any $\tau \in \R$. Similarly, there exists a positive constant $K$ depending only on $K_0$ such that
\begin{equation}
\label{allfolks2}
\| \partial_\tau \boU(\cdot, \tau) \|_{H^{k + 2}(\R)} + \| \boU(\cdot, \tau) \|_{H^{k + 5}(\R)} \leq K,
\end{equation}
for any $\tau \in \R$. In particular, we have
\begin{equation}
\label{allfolks3}
\| \partial_\tau Z_\varepsilon(\cdot, \tau) \|_{H^{k + 2}(\R)} + \| Z_\varepsilon(\cdot, \tau) \|_{H^{k + 5}(\R)} \leq K \exp K |\tau|,
\end{equation}
for any $\tau \in \R$.
\end{lemma}

Notice that the estimates of Lemma \ref{Nonbounded} do not contain any $\varepsilon^2$ factor. In order to make apparent such a factor, we rely on equations \eqref{slow2} and \eqref{equpsilon}. Indeed, in view of \eqref{slow2} and \eqref{equpsilon}, the space derivative of the functions $V_\varepsilon$ and $\Upsilon_\varepsilon$ may be replaced by the time derivative of $V_\varepsilon$ and $\Upsilon_\varepsilon$ respectively, up to some remainder terms, the whole being multiplied by the desired $\varepsilon^2$ factor.

Concerning the integral $I_\varepsilon$, all of its terms involve space derivatives of $V_\varepsilon$, so that we may invoke the argument above to gain a $\varepsilon^2$ factor. More precisely, we obtain

\begin{lemma}
\label{Itisbounded}
Let $\varepsilon > 0$ be given sufficiently small. Given any $k \in \N$, assume that the initial datum $\Psi_0(\cdot) = \Psi(\cdot, 0)$ belongs to $X^{k + 6}(\R)$ and satisfies \eqref{grinzing2} for some positive constant $K_0$. Then, there exists a positive constant $K$ depending only on $K_0$ and $k$, such that
\begin{equation}
\label{ohele}
\begin{split}
& \bigg| I_\varepsilon(\tau) - \frac{\varepsilon^2}{24} \sum_{j = 1}^k \binom{k}{j - 1} \int_0^\tau \int_\R \partial_x^k f_\varepsilon \partial_x^{k - j} V_\varepsilon \partial_x^j U_\varepsilon \bigg|\\
\leq K \varepsilon^2 \bigg( \big( \varepsilon^2 + \| \partial_x^k & Z_\varepsilon(\cdot, \tau) \|_{L^2(\R)} \big) \exp K |\tau| + \bigg| \int_0^\tau \exp K |s| \ \| Z_\varepsilon(\cdot, s) \|_{H^k(\R)} ds \bigg| \bigg),
\end{split}
\end{equation}
for any $\tau \in \R$.
\end{lemma}

The sum in the first line of \eqref{ohele}, namely
\begin{equation}
\label{sommedure}
\sum_{j = 1}^k \binom{k}{j - 1} \int_0^\tau \int_\R \partial_x^k f_\varepsilon \partial_x^{k - j} V_\varepsilon \partial_x^j U_\varepsilon
\end{equation}
corresponds to some remainder terms mentioned above in the computation of the space derivatives of $V_\varepsilon$, related to the computation of the term
$$\sum_{j = 0}^k \binom{k + 1}{j} \partial_x^j U_\varepsilon \partial_x^{k + 1 - j} V_\varepsilon \partial_x^k Z_\varepsilon,$$
appearing in $I_\eps$ . The estimates of the sum \eqref{sommedure} is more involved, so that we postpone its analysis in Lemma \ref{Ture} below.

In contrast, the integral $J_\varepsilon$ does not contain any derivative of $V_\varepsilon$. Our argument does not rely anymore on \eqref{slow2}, but instead we introduce the function $\Upsilon_\varepsilon$ in the right-hand side of \eqref{zubar} for some suitably chosen number $R$. We then invoke \eqref{equpsilon} to gain some $\varepsilon^2$ factor. This provides

\begin{lemma}
\label{Jesuisborne}
Let $\varepsilon > 0$ be given sufficiently small. Given any $k \in \N$, we assume that the initial datum $\Psi_0(\cdot) = \Psi(\cdot, 0)$ belongs to $X^{k + 6}(\R)$ and satisfies \eqref{grinzing1} for some positive constant $K_0$. Then, given any $\tau \in \R$, there exist a positive number $R_1$, depending on $\varepsilon$ and $V_\varepsilon$, and a positive constant $K$, depending only on $K_0$ and $k$, such that
\begin{equation}
\label{oheme}
\begin{split}
& \bigg| J_\varepsilon(\tau) - \frac{\varepsilon^2}{24} \int_0^\tau \int_\R \partial_x^k f_\varepsilon \Upsilon_\varepsilon \partial_x^{k + 1} U_\varepsilon \bigg|\\
\leq K \varepsilon^2 \bigg( \big( \varepsilon^2 + \| \partial_x^k & Z_\varepsilon(\cdot, \tau) \|_{L^2(\R)} \big) \exp K |\tau| + \bigg| \int_0^\tau \exp K |s| \ \| Z_\varepsilon(\cdot, s) \|_{H^k(\R)} ds \bigg| \bigg),
\end{split}
\end{equation}
for any choice of the number $R$ of definition \eqref{defupsilon} in $(R_1, + \infty)$.
\end{lemma}

Once again, the integral on the left-hand side of \eqref{oheme} corresponds to some of the extra remainder terms in \eqref{equpsilon}. The estimates of this integral are also more involved, so that we postpone its analysis in Lemma \ref{Ture} below.

In order to achieve the proof of Proposition \ref{Estim-f}, it remains to estimate the interactions terms
\begin{equation}
\label{pourri}
K_\varepsilon(\tau) \equiv \frac{\varepsilon^2}{24} \int_0^\tau \int_\R \bigg( \partial_x^k f_\varepsilon \Upsilon_\varepsilon \partial_x^{k + 1} U_\varepsilon + \sum_{j = 1}^k \binom{k}{j - 1} \partial_x^k f_\varepsilon \partial_x^{k - j} V_\varepsilon \partial_x^j U_\varepsilon \bigg),
\end{equation}
which appear on the left-hand side of \eqref{ohele} and \eqref{oheme}. Using the bounds of Proposition \ref{Bounded} and Lemmas \ref{Wanted} and \ref{Nonbounded}, we compute the estimate
$$|K_\varepsilon(\tau)| \leq K \varepsilon^2 \exp K |\tau|,$$
for any $\tau \in \R$. However, this bound is not sufficient to complete the proof of Theorem \ref{cochondore}, since it would provide an $\varepsilon$ factor instead of an $\varepsilon^2$ factor in \eqref{ineq1}. In order to gain some further $\varepsilon$ factor, we iterate the argument, and replace once more the space derivatives of $V_\varepsilon$ and $\Upsilon_\varepsilon$, which appear in the expression of $K_\varepsilon(\tau)$, by the time derivatives of $V_\varepsilon$, respectively $\Upsilon_\varepsilon$, plus some additional remainder terms. We obtain

\begin{lemma}
\label{Ture}
Let $\varepsilon > 0$ be given sufficiently small. Given any $k \in \N$, assume that the initial datum $\Psi_0(\cdot) = \Psi(\cdot, 0)$ belongs to $X^{k + 6}(\R)$ and satisfies \eqref{grinzing1} for some positive constant $K_0$. Then, given any $\tau \in \R$, there exist a positive constant $K$, depending only on $K_0$ and $k$, such that
\begin{equation}
\label{moisi}
| K_\varepsilon(\tau)| \leq K \varepsilon^4 \exp K |\tau|,
\end{equation}
for any choice of the number $R$ of definition \eqref{defupsilon} in $(R_1, + \infty)$, where $R_1$ denotes the positive number given by Lemma \ref{Jesuisborne}.
\end{lemma}

Proposition \ref{Estim-f} follows combining Lemmas \ref{Itisbounded}, \ref{Jesuisborne} and \ref{Ture}.

\begin{proof}[Proof of Proposition \ref{Estim-f} completed]
Given any $\tau \in \R$, we first fix the number $R$ such that \eqref{oheme} and \eqref{moisi} hold. Then, in view of \eqref{mandanda}, \eqref{ohele}, \eqref{oheme} and \eqref{pourri}, we have
\begin{align*}
& \bigg| \int_0^\tau \int_\R \partial_x^k f_\varepsilon \partial_x^k Z_\varepsilon - K_\varepsilon(\tau) \bigg|\\
\leq K \varepsilon^2 \bigg( \big( \varepsilon^2 + \| \partial_x^k Z_\varepsilon(\cdot, \tau) & \|_{L^2(\R)} \big) \exp K |\tau| + \bigg| \int_0^\tau \exp K |s| \ \| Z_\varepsilon(\cdot, s) \|_{H^k(\R)} ds \bigg| \bigg),
\end{align*}
so that \eqref{estim-f} is a direct consequence of \eqref{moisi}. This concludes the proof of Proposition \ref{Estim-f}.
\end{proof}

The rest of this section is devoted to the proofs of Lemmas \ref{Nonbounded}, \ref{Itisbounded}, \ref{Jesuisborne} and \ref{Ture}.

\subsection{Proof of Lemma \ref{Nonbounded}}

Concerning the nonlinear functions $f_\varepsilon$, $F_\varepsilon$, $g_\varepsilon$ and $G_\varepsilon$, it follows from definitions \eqref{fF} and \eqref{gG}, the Leibniz formula, the H\"older inequality and the Sobolev embedding theorem
\footnote{For $k > \frac{1}{2}$, one could invoke instead the fact that $H^k(\R)$ is a Banach algebra, but we write the proofs so that they work for $k \geq 0$.}
that
\begin{equation}
\label{gourcuff}
\begin{split}
& \| f_\varepsilon \|_{H^k(\R)} + \| g_\varepsilon \|_{H^k(\R)} + \| F_\varepsilon \|_{H^{k + 1}(\R)} + \| G_\varepsilon \|_{H^{k + 1}(\R)}\\
\leq K \Big( \| U_\varepsilon & \|_{H^{k + 3}(\R)} + \| V_\varepsilon \|_{H^{k + 3}(\R)} + \big( \| U_\varepsilon \|_{H^{k + 1}(\R)} + \| V_\varepsilon \|_{H^{k + 1}(\R)} \big)^2 \Big),
\end{split}
\end{equation}
where $K$ is some positive constant depending only on $k$. Similarly for $r_\varepsilon$ and $R_\varepsilon$, we have in view of $\eqref{rR}$,
\begin{equation}
\label{chamakh}
\begin{split}
\| r_\varepsilon \|_{H^k(\R)} & + \| R_\varepsilon \|_{H^{k + 1}(\R)} \leq K \| N_\varepsilon \|_{H^{k + 1}(\R)}^{k + 2} \| N_\varepsilon \|_{H^{k + 3}(\R)}\\
& \leq K \Big( \| U_\varepsilon \|_{H^{k + 1}(\R)} + \| V_\varepsilon \|_{H^{k + 1}(\R)} \Big)^{k + 2} \Big( \| U_\varepsilon \|_{H^{k + 3}(\R)} + \| V_\varepsilon \|_{H^{k + 3}(\R)} \Big).
\end{split}
\end{equation}
Notice that for \eqref{chamakh}, we also invoke the bound
$$m_0 \leq 1 - \frac{\varepsilon^2}{6} N_\varepsilon \leq \frac{1}{m_0},$$
which is a consequence of \eqref{half}. Combining \eqref{gourcuff} and \eqref{chamakh} with the bounds of Proposition \ref{Bounded}, we obtain \eqref{allfolks1}, except for the time derivative $\partial_\tau U_\varepsilon$. For this function, we have in view of \eqref{slow1},
$$\| \partial_\tau U_\varepsilon \|_{H^k(\R)} \leq \| U_\varepsilon \|_{H^{k + 3}(\R)} + \| U_\varepsilon \|_{H^{k + 1}(\R)}^2 + \| f_\varepsilon \|_{H^k(\R)} + \varepsilon^2 \| r_\varepsilon \|_{H^k(\R)},$$
so that the bound for $\partial_\tau U_\varepsilon$ follows from the previous bounds on $f_\varepsilon$ and $r_\varepsilon$ combined with the bound on $U_\varepsilon$ of Proposition \ref{Bounded}.

For the uniform bound \eqref{allfolks2} on $\boU$, we invoke the integrability properties of the Korteweg-de Vries equation. As a matter of fact, equation \eqref{KdV} owns an infinite number of invariant quantities which control the $H^k$-norms of the solutions (see e.g. \cite{GarKrMi1}). Therefore, the $H^k$-norm of a solution $\boU$ at time $\tau$ is controlled by the $H^k$-norm of its initial datum $\boU^0 = U_\varepsilon^0$. More precisely, there exists a positive constant $K$, depending only on $\| U_\varepsilon^0 \|_{H^k(\R)}$, such that
$$\| \boU(\cdot, t) \|_{H^k(\R)} \leq K,$$
for any $t \in \R$ (see \cite{BonaSmi1}). This control is uniform with respect to $\varepsilon$ provided that $\| U_\varepsilon^0 \|_{H^k(\R)}$ is bounded independently of $\varepsilon$. Since $\boU$ solves the Korteweg-de Vries equation, the $H^{k - 3}$-norm of $\partial_\tau \boU$ is then uniformly bounded with respect to $\varepsilon$. Thus, under assumption \eqref{grinzing1}, we obtain inequality \eqref{allfolks2}. Estimate \eqref{allfolks3} follows combining the definition of $Z_\varepsilon$ with \eqref{allfolks1} and \eqref{allfolks2}.

\subsection{Proof of Lemma \ref{Itisbounded}}

Given any $\tau \in \R$, we split the expression of $I_\varepsilon(\tau)$ into three terms
\begin{equation}
\label{lloris}
I_\varepsilon(\tau) = I_1(\tau) + I_2(\tau) + I_3(\tau),
\end{equation}
where
\begin{equation}
\label{I1}
I_1(\tau) \equiv \frac{1}{6} \int_0^\tau \int_\R \partial_x^k Z_\varepsilon \partial_x^{k + 1} \big( V_\varepsilon^2 \big),
\end{equation}
$$I_2(\tau) \equiv - \int_0^\tau \int_\R \partial_x^k Z_\varepsilon \partial_x^{k + 3} V_\varepsilon,$$
and
$$I_3(\tau) \equiv \frac{1}{3} \sum_{j = 0}^k \binom{k + 1}{j} \int_0^\tau \int_\R \partial_x^k Z_\varepsilon \partial_x^j U_\varepsilon \partial_x^{k + 1 - j} V_\varepsilon.$$
We now compute estimates of each term $I_k(\tau)$. For the first one, we have

\begin{step}
\label{Estim-I1}
Under the assumptions of Lemma \ref{Jesuisborne}, there exists a positive constant $K$ depending only on $K_0$ and $k$, such that
\begin{equation}
\label{estim-I1}
\begin{split}
& \bigg| I_1(\tau) - \frac{\varepsilon^2}{48} \int_0^\tau \int_\R \partial_x^k \big( V_\varepsilon^2 \big) \partial_x^k f_\varepsilon \bigg|\\
\leq K \varepsilon^2 \bigg( \big( \varepsilon^2 + \| \partial_x^k Z_\varepsilon(\cdot, \tau) & \|_{L^2(\R)} \big) \exp K |\tau| + \bigg| \int_0^\tau \exp K |s| \ \| Z_\varepsilon(\cdot, s) \|_{H^k(\R)} ds \bigg| \bigg).
\end{split}
\end{equation}
\end{step}

In order to obtain the $\varepsilon^2$ factor on the right-hand side of \eqref{estim-I1}, we write the transport equation satisfied by $\partial_x^k (V_\varepsilon^2)$, namely
\begin{equation}
\label{eq:Vcarre}
\partial_x \Big( \partial_x^k \big( V_\varepsilon^2 \big) \Big) = - \frac{\varepsilon^2}{8} \partial_\tau \Big( \partial_x^k \big( V_\varepsilon^2 \big) \Big) + \frac{\varepsilon^2}{4} \partial_x^k \big( g_\varepsilon V_\varepsilon \big) + \frac{\varepsilon^4}{4} \partial_x^k \big( r_\varepsilon V_\varepsilon \big),
\end{equation}
and replace $\partial_x^{k + 1} \big( V_\varepsilon^2 \big)$ in \eqref{I1} by its expression provided by \eqref{eq:Vcarre}. Integrating by parts in time and using the fact that $Z_\varepsilon^0 = 0$, the integral $I_1(\tau)$ becomes
\begin{equation}
\label{boumsong}
\begin{split}
I_1(\tau) = & - \frac{\varepsilon^2}{48} \int_\R \partial_x^k Z_\varepsilon(x, \tau) \partial_x^k \big( V_\varepsilon^2 \big) (x, \tau) dx + \frac{\varepsilon^2}{48} \int_0^\tau \int_\R \partial_\tau \partial_x^k Z_\varepsilon \partial_x^k \big( V_\varepsilon^2 \big)\\
& + \frac{\varepsilon^2}{24} \int_0^\tau \int_\R \partial_x^k Z_\varepsilon \partial_x^k \big( g_\varepsilon V_\varepsilon \big) + \frac{\varepsilon^4}{24} \int_0^\tau \int_\R \partial_x^k Z_\varepsilon \partial_x^k \big( r_\varepsilon V_\varepsilon \big).
\end{split}
\end{equation}
We now apply the Leibniz formula, the H\"older inequality and the Sobolev embedding theorem to obtain for the first term on the right-hand side of \eqref{boumsong},
$$\bigg| \int_\R \partial_x^k Z_\varepsilon(x, \tau) \partial_x^k \big( V_\varepsilon^2 \big)(x, \tau) dx \bigg| \leq K \| \partial_x^k Z_\varepsilon(\cdot, \tau) \|_{L^2(\R)} \| V_\varepsilon(\cdot, \tau) \|_{H^k(\R)} \| V_\varepsilon(\cdot, \tau) \|_{H^{k + 1}(\R)}.$$
Hence, we deduce from \eqref{frites} that there exists a positive constant $K$ depending only on $K_0$ and $k$, such that
\begin{equation}
\label{cris}
\bigg| \int_\R \partial_x^k Z_\varepsilon(x, \tau) \partial_x^k \big( V_\varepsilon^2 \big)(x, \tau) dx \bigg| \leq K \| \partial_x^k Z_\varepsilon(\cdot, \tau) \|_{L^2(\R)} \exp K |\tau|.
\end{equation}
Similarly, in view of Lemma \ref{Nonbounded}, we have for the third and fourth terms
\begin{equation}
\label{mensah}
\begin{split}
& \bigg| \int_0^\tau \int_\R \partial_x^k Z_\varepsilon \partial_x^k \big( g_\varepsilon V_\varepsilon \big) \bigg| + \bigg| \int_0^\tau \int_\R \partial_x^k Z_\varepsilon \partial_x^k \big( r_\varepsilon V_\varepsilon \big) \bigg|\\
\leq & K \bigg| \int_0^\tau \| \partial_x^k Z_\varepsilon \|_{L^2(\R)} \| V_\varepsilon \|_{H^{k + 1}(\R)} \big( \| r_\varepsilon \|_{H^k(\R)} + \| g_\varepsilon \|_{H^k(\R)}\big) \bigg|\\
\leq & K \bigg| \int_0^\tau \exp K |s| \ \| \partial_x^k Z_\varepsilon(\cdot, s) \|_{L^2(\R)} ds \bigg|.
\end{split}
\end{equation}
The analysis of the second term on the right-hand side of \eqref{boumsong} is more involved since estimate \eqref{estim-I1} requires to gain some further $\varepsilon^2$ factor. We first replace $\partial_\tau \partial_x^k Z_\varepsilon$ by its expression given by \eqref{dede}. We obtain
\begin{equation}
\label{toulalan}
\int_0^\tau \int_\R \partial_\tau \partial_x^k Z_\varepsilon \partial_x^k \big( V_\varepsilon^2 \big) = I_1^1(\tau) + I_1^2(\tau) + \int_0^\tau \int_\R \partial_x^k \big( V_\varepsilon^2 \big) \partial_x^k f_\varepsilon,
\end{equation}
where
$$I_1^1(\tau) \equiv - \int_0^\tau \int_\R \partial_x^k \big( V_\varepsilon^2 \big) \Big( \partial_x^{k + 3} Z_\varepsilon + \partial_x^{k + 1} \big( \boU Z_\varepsilon \big) + \frac{1}{2} \partial_x^{k + 1} \big( Z_\varepsilon^2 \big) \Big),$$
and
$$I_1^2(\tau) \equiv - \varepsilon^2 \int_0^\tau \int_\R \partial_x^k \big( V_\varepsilon^2 \big) \partial_x^k r_{\varepsilon}.$$
We then estimate each of the above integral $I_1^j(\tau)$. Integrating by parts in space the first one, we have
$$I_1^1(\tau) = \int_0^\tau \int_\R \Big( \partial_x^{k + 3} \big( V_\varepsilon^2 \big) \partial_x^k Z_\varepsilon + \partial_x^{k + 1} \big( V_\varepsilon^2 \big) \partial_x^k \big( \boU Z_\varepsilon \big) + \frac{1}{2} \partial_x^{k + 1} \big( V_\varepsilon^2 \big) \partial_x^k \big( Z_\varepsilon^2 \big) \Big),$$
so that
\begin{align*}
|I_1^1(\tau)| \leq \bigg| \int_0^\tau \| V_\varepsilon \|_{H^{k + 3}(\R)} \| Z_\varepsilon \|_{H^k(\R)} \Big( \| V_\varepsilon \|_{H^{k + 3}(\R)} + \| V_\varepsilon \|_{H^{k + 1}(\R)} \big( \| \boU \|_{H^k(\R)} + \| Z_\varepsilon \|_{H^k(\R)} \big) \Big) \bigg|.
\end{align*}
In view of the bounds of Proposition \ref{Bounded} and Lemma \ref{Nonbounded}, we are led to
\begin{equation}
\label{juninho}
|I_1^1(\tau)| \leq K \bigg| \int_0^\tau \exp K |s| \ \| Z_\varepsilon(\cdot, s) \|_{H^k(\R)} ds \bigg|.
\end{equation}
Similarly, in view of the bounds of Proposition \ref{Bounded} and Lemma \ref{Nonbounded}, we compute
\begin{equation}
\label{ederson}
|I_1^2(\tau)| \leq K \varepsilon^2 \bigg| \int_0^\tau \| V_\varepsilon \|_{H^k(\R)} \| V_\varepsilon \|_{H^{k + 1}(\R)} \| r_\varepsilon \|_{H^k(\R)} \bigg| \leq K \varepsilon^2 \exp K |\tau|.
\end{equation}
Combining \eqref{toulalan} with \eqref{juninho} and \eqref{ederson}, we are led to
\begin{align*}
& \bigg| \int_0^\tau \int_\R \partial_\tau \partial_x^k Z_\varepsilon \partial_x^k \big( V_\varepsilon^2 \big) - \int_0^\tau \int_\R \partial_x^k \big( V_\varepsilon^2 \big) \partial_x^k f_\varepsilon \bigg|
\\ \leq K \bigg( & \varepsilon^2 \exp K |\tau| + \bigg| \int_0^\tau \exp K |s| \ \| Z_\varepsilon(\cdot, s) \|_{H^k(\R)} ds \bigg| \bigg).
\end{align*}
In view of \eqref{boumsong}, \eqref{cris} and \eqref{mensah}, this completes the proof of Step \ref{Estim-I1}.

We now turn to the integral $I_2(\tau)$ for which we have

\begin{step}
\label{Estim-I2}
Under the assumptions of Lemma \ref{Jesuisborne}, there exists a positive constant $K$ depending only on $K_0$ and $k$, such that
\begin{equation}
\label{estim-I2}
\begin{split}
& \bigg| I_2(\tau) + \frac{\varepsilon^2}{8} \int_0^\tau \int_\R \partial_x^{k + 2} V_\varepsilon \partial_x^k f_\varepsilon \bigg|\\
\leq K \varepsilon^2 \bigg( \big( \varepsilon^2 + \| \partial_x^k Z_\varepsilon(\cdot, \tau) & \|_{L^2(\R)} \big) \exp K |\tau| + \bigg| \int_0^\tau \exp K |s| \ \| Z_\varepsilon(\cdot, s) \|_{H^k(\R)} ds \bigg| \bigg).
\end{split}
\end{equation}
\end{step}

The proof is similar to the proof of Step \ref{Estim-I1}. In view of \eqref{slow2}, we have integrating by parts in time and using the fact that $Z_\varepsilon^0 = 0$,
\begin{equation}
\label{benzema}
\begin{split}
I_2(\tau) & = \frac{\varepsilon^2}{8} \int_\R \partial_x^k Z_\varepsilon(x, \tau) \partial_x^{k + 2} V_\varepsilon(x, \tau) dx - \frac{\varepsilon^2}{8} \int_0^\tau \int_\R \partial_\tau \partial_x^k Z_\varepsilon \partial_x^{k + 2} V_\varepsilon\\
& - \frac{\varepsilon^2}{8} \int_0^\tau \int_\R \partial_x^{k + 2} g_\varepsilon \partial_x^k Z_\varepsilon - \frac{\varepsilon^4}{8} \int_0^\tau \int_\R \partial_x^{k + 2} r_\varepsilon \partial_x^k Z_\varepsilon.
\end{split}
\end{equation}
Similarly to \eqref{cris}, we compute for the first term on the right-hand side of \eqref{benzema},
\begin{equation}
\label{govou}
\bigg| \int_\R \partial_x^k Z_\varepsilon(x, \tau) \partial_x^{k + 2} V_\varepsilon(x, \tau) dx \bigg| \leq K \| \partial_x^k Z_\varepsilon(\cdot, \tau) \|_{L^2(\R)} \exp K |\tau|,
\end{equation}
while for the third and fourth terms, we have similarly to \eqref{mensah},
\begin{equation}
\label{fred}
\bigg| \int_0^\tau \int_\R \partial_x^{k + 2} g_\varepsilon \partial_x^k Z_\varepsilon \bigg| + \bigg| \int_0^\tau \int_\R \partial_x^{k + 2} r_\varepsilon \partial_x^k Z_\varepsilon \bigg| \leq K \int_0^\tau \exp K |s| \ \| \partial_x^k Z_\varepsilon(\cdot, s) \|_{L^2(\R)} ds.
\end{equation}
Concerning the second term on the right-hand side of \eqref{benzema}, we replace as above $\partial_\tau \partial_x^k Z_\varepsilon$ by its expression given by \eqref{dede}. This leads to
\begin{equation}
\label{delgado}
\int_0^\tau \int_\R \partial_\tau \partial_x^k Z_\varepsilon \partial_x^{k + 2} V_\varepsilon = I_2^1(\tau) + I_2^2(\tau) + \int_0^\tau \int_\R \partial_x^{k + 2} V_\varepsilon \partial_x^k f_\varepsilon,
\end{equation}
where
$$I_2^1(\tau) \equiv - \int_0^\tau \int_\R \partial_x^{k + 2} V_\varepsilon \Big( \partial_x^{k + 3} Z_\varepsilon + \partial_x^{k + 1} \big( \boU Z_\varepsilon \big) + \frac{1}{2} \partial_x^{k + 1} \big( Z_\varepsilon^2 \big) \Big),$$
and
$$I_2^2(\tau) \equiv - \varepsilon^2 \int_0^\tau \int_\R \partial_x^{k + 2} V_\varepsilon \partial_x^k r_{\varepsilon}.$$
Following the lines of the proofs of \eqref{juninho} and \eqref{ederson}, we obtain
$$|I_2^1(\tau)| \leq K \bigg| \int_0^\tau \exp K |s| \ \| Z_\varepsilon(\cdot, s) \|_{H^k(\R)} ds \bigg|,$$
and
$$|I_2^2(\tau)| \leq K \varepsilon^2 \exp K |\tau|.$$
By \eqref{delgado}, we are led to
\begin{align*}
& \bigg| \int_0^\tau \int_\R \partial_\tau \partial_x^k Z_\varepsilon \partial_x^{k + 2} V_\varepsilon - \int_0^\tau \int_\R \partial_x^{k + 2} V_\varepsilon \partial_x^k f_\varepsilon \bigg|\\
\leq K \bigg( & \varepsilon^2 \exp K |\tau| + \bigg| \int_0^\tau \exp K |s| \ \| Z_\varepsilon(\cdot, s) \|_{H^k(\R)} ds \bigg| \bigg).
\end{align*}
Estimate \eqref{estim-I2} follows combining with \eqref{benzema}, \eqref{govou} and \eqref{fred}.

We finally consider the finite sum $I_3(\tau)$.

\begin{step}
\label{Estim-I3}
Under the assumptions of Lemma \ref{Jesuisborne}, there exists a positive constant $K$ depending only on $K_0$ and $k$, such that
\begin{equation}
\label{estim-I3}
\begin{split}
& \bigg| I_3(\tau) - \frac{\varepsilon^2}{24} \int_0^\tau \int_\R \partial_x^k \big( U_\varepsilon V_\varepsilon \big) \partial_x^k f_\varepsilon - \frac{\varepsilon^2}{24} \sum_{j = 1}^k \binom{k}{j - 1} \int_0^\tau \int_\R \partial_x^k f_\varepsilon \partial_x^{k - j} V_\varepsilon \partial_x^j U_\varepsilon \bigg|\\
\leq & K \varepsilon^2 \bigg( \big( \varepsilon^2 + \| \partial_x^k Z_\varepsilon(\cdot, \tau) \|_{L^2(\R)} \big) \exp K |\tau| + \bigg| \int_0^\tau \exp K |s| \ \| Z_\varepsilon(\cdot, s) \|_{H^k(\R)} ds \bigg| \bigg).
\end{split}
\end{equation}
\end{step}

The proof is similar to the proof of Steps \ref{Estim-I1} and \ref{Estim-I2}. In view of \eqref{slow2}, we have integrating by parts in time,
\begin{equation}
\label{grosso}
\begin{split}
I_3(\tau) & = \frac{\varepsilon^2}{24} \sum_{j = 0}^k \binom{k + 1}{j} \Bigg( - \int_\R \partial_x^k Z_\varepsilon(x, \tau) \partial_x^j U_\varepsilon(x, \tau) \partial_x^{k - j} V_\varepsilon(x, \tau) dx + \int_0^\tau \int_\R \partial_\tau \partial_x^k Z_\varepsilon \partial_x^j U_\varepsilon \partial_x^{k - j} V_\varepsilon\\
& + \int_0^\tau \int_\R \partial_x^k Z_\varepsilon \partial_\tau \partial_x^j U_\varepsilon \partial_x^{k - j} V_\varepsilon + \int_0^\tau \int_\R \partial_x^k Z_\varepsilon \partial_x^j U_\varepsilon \partial_x^{k - j} g_\varepsilon + \varepsilon^2 \int_0^\tau \int_\R \partial_x^k Z_\varepsilon \partial_x^j U_\varepsilon \partial_x^{k - j} r_\varepsilon \Bigg).
\end{split}
\end{equation}
The estimates of the first, fourth and fifth integrals on the right-hand side of \eqref{grosso} are similar to \eqref{cris} and \eqref{govou} for the first one, \eqref{mensah} and \eqref{fred}, for the other ones. More precisely, we obtain
\begin{equation}
\label{mounier}
\bigg| \int_\R \partial_x^k Z_\varepsilon(x, \tau) \partial_x^j U_\varepsilon(x, \tau) \partial_x^{k - j} V_\varepsilon(x, \tau) dx \bigg| \leq K \| \partial_x^k Z_\varepsilon(\cdot, \tau) \|_{L^2(\R)} \exp K |\tau|,
\end{equation}
respectively
\begin{equation}
\label{keita}
\bigg| \int_0^\tau \int_\R \partial_x^k Z_\varepsilon \partial_x^j U_\varepsilon \partial_x^{k - j} g_\varepsilon \bigg| + \bigg| \int_0^\tau \int_\R \partial_x^k Z_\varepsilon \partial_x^j U_\varepsilon \partial_x^{k - j} r_\varepsilon \bigg| \leq K \bigg| \int_0^\tau \exp K |s| \ \| \partial_x^k Z_\varepsilon(\cdot, s) \|_{L^2(\R)} ds \bigg|.
\end{equation}
Similarly, for the third integral, we have
\begin{equation}
\label{makoun}
\begin{split}
\bigg| \int_0^\tau \int_\R \partial_x^k Z_\varepsilon \partial_\tau \partial_x^j U_\varepsilon \partial_x^{k - j} V_\varepsilon \bigg| \leq & \bigg| \int_0^\tau \| \partial_x^k Z_\varepsilon \|_{L^2(\R)} \| \partial_\tau U_\varepsilon \|_{H^k(\R)} \| V_\varepsilon \|_{H^{k + 1}(\R)} \bigg|\\
\leq & K \bigg| \int_0^\tau \exp K |s| \ \| \partial_x^k Z_\varepsilon(\cdot, s) \|_{L^2(\R)} ds \bigg|,
\end{split}
\end{equation}
where we invoke \eqref{allfolks1} to bound the time derivative $\partial_\tau \partial_x^j U_\varepsilon$.

The second term on the right-hand side of \eqref{grosso} is more involved to estimate. We first introduce as above the expression of $\partial_\tau \partial_x^k Z_\varepsilon$ given by \eqref{dede}, so that
\begin{equation}
\label{clerc}
\sum_{j = 0}^k \binom{k + 1}{j} \int_0^\tau \int_\R \partial_\tau \partial_x^k Z_\varepsilon \partial_x^j U_\varepsilon \partial_x^{k - j} V_\varepsilon = I_3^1(\tau) + I_3^2(\tau) + I_3^3(\tau),
\end{equation}
where
$$I_3^1(\tau) \equiv - \sum_{j = 0}^k \binom{k + 1}{j} \int_0^\tau \int_\R \partial_x^j U_\varepsilon \partial_x^{k - j} V_\varepsilon \Big( \partial_x^{3 + k} Z_\varepsilon + \partial_x^{k + 1} \big( \boU Z_\varepsilon \big) + \partial_x^k \big( Z_\varepsilon \partial_x Z_\varepsilon \big) \Big),$$
$$I_3^2(\tau) \equiv \sum_{j = 0}^k \binom{k + 1}{j} \int_0^\tau \int_\R \partial_x^j U_\varepsilon \partial_x^{k - j} V_\varepsilon \partial_x^k f_\varepsilon,$$
and
$$I_3^3(\tau) \equiv - \varepsilon^2 \sum_{j = 0}^k \binom{k + 1}{j} \int_0^\tau \int_\R \partial_x^j U_\varepsilon \partial_x^{k - j} V_\varepsilon \partial_x^k r_\varepsilon.$$
We then estimate each of the above sum $I_3^j(\tau)$. Integrating by parts in space the integrals in $I_3^1(\tau)$, we obtain
$$I_3^1(\tau) = \sum_{j = 0}^k \binom{k + 1}{j} \int_0^\tau \int_\R \bigg( \partial_x^3 \Big( \partial_x^j U_\varepsilon \partial_x^{k - j} V_\varepsilon \Big) \partial_x^k Z_\varepsilon + \partial_x \Big( \partial_x^j U_\varepsilon \partial_x^{k - j} V_\varepsilon \Big) \Big( \partial_x^k \big( \boU Z_\varepsilon \big) + \frac{1}{2} \partial_x^k \big( Z_\varepsilon^2 \big) \Big) \bigg).$$
Following the lines of the proof of \eqref{juninho}, we are led to
\begin{equation}
\label{reveillere}
|I_3^1(\tau)| \leq K \bigg| \int_0^\tau \int_\R \exp K |s| \ \| Z_\varepsilon(\cdot, s) \|_{H^k(\R)} ds \bigg|.
\end{equation}
Concerning the integral $I_3^3(\tau)$, we have as in the proof of \eqref{ederson},
\begin{equation}
\label{vercoutre}
|I_3^3(\tau)| \leq K \varepsilon^2 \exp K |\tau|.
\end{equation}
Finally, applying the Pascal rule and the Leibniz formula to $I_3^2(\tau)$, we obtain
$$I_3^2(\tau) = \int_0^\tau \int_\R \partial_x^k \big( U_\varepsilon V_\varepsilon \big) \partial_x^k f_\varepsilon + \sum_{j = 1}^k \binom{k}{j- 1} \int_0^\tau \int_\R \partial_x^j U_\varepsilon \partial_x^{k - j} V_\varepsilon \partial_x^k f_\varepsilon.$$
Combining with \eqref{grosso}, \eqref{mounier}, \eqref{keita}, \eqref{makoun}, \eqref{clerc}, \eqref{reveillere} and \eqref{vercoutre}, estimate \eqref{estim-I3} follows.

We are now in position to complete the proof of Lemma \ref{Itisbounded}.

\begin{proof}[End of the proof of Lemma \ref{Itisbounded}]
In view of identity \eqref{lloris}, and estimates \eqref{estim-I1}, \eqref{estim-I2} and \eqref{estim-I3}, we have
\begin{equation}
\label{kallstrom}
\begin{split}
\bigg| I_\varepsilon(\tau) & - \frac{\varepsilon^2}{8} \int_0^\tau \int_\R \partial_x^k f_\varepsilon \partial_x^k \Big( \frac{1}{6} V_\varepsilon^2 - \partial_x^2 V_\varepsilon + \frac{1}{3} U_\varepsilon V_\varepsilon \Big) - \frac{\varepsilon^2}{24} \sum_{j = 1}^k \binom{k}{j - 1} \int_0^\tau \int_\R \partial_x^k f_\varepsilon \partial_x^{k - j} V_\varepsilon \partial_x^j U_\varepsilon \bigg|\\
& \leq K \varepsilon^2 \bigg( \big( \varepsilon^2 + \| \partial_x^k Z_\varepsilon(\cdot, \tau) \|_{L^2(\R)} \big) \exp K |\tau| + \bigg| \int_0^\tau \exp K |s| \ \| Z_\varepsilon(\cdot, s) \|_{H^k(\R)} ds \bigg| \bigg).
\end{split}
\end{equation}
On the other hand, in view of the definition of $f_\varepsilon$,
$$\int_0^\tau \int_\R \partial_x^k f_\varepsilon \partial_x^k \Big( \frac{1}{6} V_\varepsilon^2 - \partial_x^2 V_\varepsilon + \frac{1}{3} U_\varepsilon V_\varepsilon \Big) = \int_0^\tau \int_\R \partial_x^{k + 1} F_\varepsilon \partial_x^k F_\varepsilon = 0.$$
Combining with \eqref{kallstrom}, this completes the proof of \eqref{ohele}.
\end{proof}

\subsection{Proof of Lemma \ref{Jesuisborne}}

The proof is similar to the proof of Lemma \ref{Itisbounded}. Given any $\tau \in \R$, we introduce the function $\Upsilon_\varepsilon$ in the expression of the integral $J_\varepsilon(\tau)$. In view of \eqref{equpsilon} and \eqref{zubar}, this yields
\begin{equation}
\label{taiwo}
J_\varepsilon(\tau) = \frac{1}{3} \int_0^\tau \int_\R \partial_x^{k + 1} U_\varepsilon \partial_x^k Z_\varepsilon \partial_x \Upsilon_\varepsilon = \frac{1}{3} \Big( J_1(\tau) + J_2(\tau) + J_3(\tau) \Big),
\end{equation}
where
$$J_1(\tau) \equiv - \frac{\varepsilon^2}{8} \int_0^\tau \int_\R \partial_x^{k + 1} U_\varepsilon \partial_x^k Z_\varepsilon \partial_\tau \Upsilon_\varepsilon,$$
\begin{equation}
\label{J2}
J_2(\tau) \equiv \frac{\varepsilon^2}{8} \int_0^\tau \int_\R \partial_x^{k + 1} U_\varepsilon \partial_x^k Z_\varepsilon \big( G_\varepsilon + \varepsilon^2 R_\varepsilon \big),
\end{equation}
and
\begin{equation}
\label{J3}
J_3(\tau) \equiv \int_0^\tau \Big( V_\varepsilon(- R, s) - \frac{\varepsilon^2}{8} \Big( G_\varepsilon(- R, s) + \varepsilon^2 R_\varepsilon(- R, s) \Big) \Big) \bigg( \int_\R \partial_x^{k + 1} U_\varepsilon(x, s) \partial_x^k Z_\varepsilon(x, s) dx \bigg) ds.
\end{equation}
We now consider each of the above integrals $J_k(\tau)$. For the first one, we have

\setcounter{step}{0}
\begin{step}
\label{Estim-J1}
Under the assumptions of Lemma \ref{Jesuisborne}, there exists a positive constant $K$ depending only on $K_0$ and $k$, such that
\begin{equation}
\label{estim-J1}
\begin{split}
& \bigg| J_1(\tau) - \frac{\varepsilon^2}{8} \int_0^\tau \int_\R \partial_x^{k + 1} U_\varepsilon \Upsilon_\varepsilon \partial_x^k f_\varepsilon \bigg|\\
\leq K \varepsilon^2 \bigg( \Big( \varepsilon^2 + \| \partial_x^k Z_\varepsilon(\cdot, \tau) & \|_{L^2(\R)} \Big) \exp K |\tau| + \bigg| \int_0^\tau \exp K |s| \ \| \partial_x^k Z_\varepsilon(\cdot, s) \|_{L^2(\R)} ds \bigg| \bigg).
\end{split}
\end{equation}
\end{step}

We first integrate by parts in time the integral $J_1(\tau)$, and use the fact that $Z_\varepsilon^0 = 0$ to derive
\begin{equation}
\label{bonnart}
\begin{split}
J_1(\tau) = - & \frac{\varepsilon^2}{8} \int_\R \partial_x^{k + 1} U_\varepsilon(x, \tau) \partial_x^k Z_\varepsilon(x, \tau) \Upsilon_\varepsilon(x, \tau) dx + \frac{\varepsilon^2}{8} \int_0^\tau \int_\R \partial_x^{k + 1} \partial_\tau U_\varepsilon \partial_x^k Z_\varepsilon \Upsilon_\varepsilon\\
+ & \frac{\varepsilon^2}{8} \int_0^\tau \int_\R \partial_x^{k + 1} U_\varepsilon \partial_\tau \partial_x^k Z_\varepsilon \Upsilon_\varepsilon.
\end{split}
\end{equation}
We then invoke Proposition \ref{Bounded} and Lemmas \ref{Wanted} and \ref{Nonbounded} to bound each term on the right-hand side of \eqref{bonnart}. Notice in particular that in view of assumption \eqref{grinzing1}, Lemma \ref{Wanted} provides
\begin{equation}
\label{cana}
\| \Upsilon_\varepsilon(\cdot, s) \|_{L^\infty(\R)} \leq \| V_\varepsilon(\cdot, s) \|_{\boM(\R)} \leq K \big( 1 + |s| \big),
\end{equation}
where $K$ denotes some positive constant, possibly depending on $K_0$, but neither on $\varepsilon$ nor $s$.

Hence, applying the H\"older inequality to the first term on the right-hand side of \eqref{bonnart}, we are led to
$$\bigg| \int_\R \partial_x^{k + 1} U_\varepsilon(x, \tau) \partial_x^k Z_\varepsilon(x, \tau) \Upsilon_\varepsilon(x, \tau) dx \bigg| \leq \| \partial_x^{k + 1} U_\varepsilon(\cdot, \tau) \|_{L^2(\R)} \| \partial_x^k Z_\varepsilon(\cdot, \tau) \|_{L^2(\R)} \| \Upsilon_\varepsilon(\cdot, \tau) \|_{L^\infty(\R)},$$
so that by \eqref{frites} and \eqref{cana},
\begin{equation}
\label{ziani}
\bigg| \int_\R \partial_x^{k + 1} U_\varepsilon(x, \tau) \partial_x^k Z_\varepsilon(x, \tau) \Upsilon_\varepsilon(x, \tau) dx \bigg| \leq K \| \partial_x^k Z_\varepsilon(\cdot, \tau) \|_{L^2(\R)} \exp K |\tau|.
\end{equation}
Similarly, for the second term on the right-hand side of \eqref{bonnart}, we compute
$$\bigg| \int_0^\tau \int_\R \partial_x^{k + 1} \partial_\tau U_\varepsilon \partial_x^k Z_\varepsilon \Upsilon_\varepsilon \bigg| \leq K \bigg| \int_0^\tau \big( 1 + |s| \big) \| \partial_x^{k + 1} \partial_\tau U_\varepsilon(\cdot, s) \|_{L^2(\R)} \| \partial_x^k Z_\varepsilon(\cdot, s) \|_{L^2(\R)} ds \bigg|.$$
so that by \eqref{allfolks1},
\begin{equation}
\label{cheyrou}
\bigg| \int_0^\tau \int_\R \partial_x^{k + 1} \partial_\tau U_\varepsilon \partial_x^k Z_\varepsilon \Upsilon_\varepsilon \bigg| \leq K \bigg| \int_0^\tau \exp K |s| \ \| \partial_x^k Z_\varepsilon(\cdot, s) \|_{L^2(\R)} ds \bigg|.
\end{equation}
We finally turn to the last term on the right-hand side of \eqref{bonnart}. In view of \eqref{dede}, we can write
\begin{equation}
\label{kone}
\int_0^\tau \int_\R \partial_x^{k + 1} U_\varepsilon \partial_\tau \partial_x^k Z_\varepsilon \Upsilon_\varepsilon = J_1^1(\tau) + J_1^2(\tau) + \int_0^\tau \int_\R \partial_x^{k + 1} U_\varepsilon \Upsilon_\varepsilon \partial_x^k f_\varepsilon,
\end{equation}
where
$$J_1^1(\tau) \equiv - \int_0^\tau \int_\R \partial_x^{k + 1} U_\varepsilon \Upsilon_\varepsilon \Big( \partial_x^{k + 3} Z_\varepsilon + \partial_x^{k + 1} \big( Z_\varepsilon \boU \big) + \frac{1}{2} \partial_x^{k + 1} \big( Z_\varepsilon^2 \big) \Big).$$
and
$$J_1^2(\tau) \equiv - \varepsilon^2 \int_0^\tau \int_\R \partial_x^{k + 1} U_\varepsilon \Upsilon_\varepsilon \partial_x^k r_{\varepsilon}.$$
Concerning the integral $J_1^1(\tau)$, we integrate by parts in space to obtain
$$J_1^1(\tau) = \int_0^\tau \int_\R \bigg( \partial_x^k Z_\varepsilon \partial_x^3 \big( \partial_x^{k + 1} U_\varepsilon \Upsilon_\varepsilon \big) + \partial_x^k \big( Z_\varepsilon \boU \big) \partial_x \big( \partial_x^{k + 1} U_\varepsilon \Upsilon_\varepsilon \big) + \frac{1}{2} \partial_x^k \big( Z_\varepsilon^2 \big) \partial_x \big( \partial_x^{k + 1} U_\varepsilon \Upsilon_\varepsilon \big) \bigg),$$
so that
\begin{align*}
|J_1^1(\tau)| \leq \bigg| \int_0^\tau \| Z_\varepsilon \|_{H^k(\R)} \| U_\varepsilon \|_{H^{k + 4}(\R)} \big( \| \Upsilon_\varepsilon \|_{L^\infty(\R)} + \| V_\varepsilon \|_{H^3(\R)} \big) \big( 1 + \| \boU \|_{H^{k + 1}(\R)} + \| Z_\varepsilon \|_{H^{k + 1}(\R)} \big) \bigg|.
\end{align*}
Hence, by \eqref{frites}, \eqref{allfolks2}, \eqref{allfolks3} and \eqref{cana},
\begin{equation}
\label{kabore}
|J_1^1(\tau)| \leq K \bigg| \int_0^\tau \exp K |s| \ \| Z_\varepsilon(\cdot, s) \|_{H^k(\R)} ds \bigg|.
\end{equation}
Similarly, the integral $J_1^2(\tau)$ is bounded by
\begin{equation}
\label{samassa}
|J_1^2(\tau)| \leq K \varepsilon^2 \bigg| \int_0^\tau \| \partial_x^{k + 1} U_\varepsilon \|_{L^2(\R)} \| \Upsilon_\varepsilon \|_{L^\infty(\R)} \| \partial_x^k r_{\varepsilon} \|_{L^2(\R)} \bigg| \leq K \varepsilon^2 \exp K |\tau|.
\end{equation}
Combining \eqref{kone} with \eqref{kabore} and \eqref{samassa}, we have
\begin{align*}
& \bigg| \int_0^\tau \int_\R \partial_x^{k + 1} U_\varepsilon \partial_\tau \partial_x^k Z_\varepsilon \Upsilon_\varepsilon - \int_0^\tau \int_\R \partial_x^{k + 1} U_\varepsilon \Upsilon_\varepsilon \partial_x^k f_\varepsilon \bigg|\\
\leq & K \bigg( \varepsilon^2 \exp K |\tau| + \bigg| \int_0^\tau \exp K |s| \ \| Z_\varepsilon(\cdot, s) \|_{H^k(\R)} ds \bigg| \bigg),
\end{align*}
so that \eqref{estim-J1} follows from \eqref{bonnart}, \eqref{ziani} and \eqref{cheyrou}.

We now turn to the second integral $J_2(\tau)$.

\begin{step}
\label{Estim-J2}
Under the assumptions of Lemma \ref{Jesuisborne}, there exists a positive constant $K$ depending only on $K_0$ and $k$, such that
\begin{equation}
\label{estim-J2}
|J_2(\tau)| \leq K \varepsilon^2 \bigg| \int_0^\tau \exp K |s| \ \| \partial_x^k Z_\varepsilon(\cdot, s) \|_{L^2(\R)} ds \bigg|,
\end{equation}
for any $\tau \in \R$.
\end{step}

Applying the H\"older inequality and the Sobolev embedding theorem to definition \eqref{J2}, we obtain
$$|J_2(\tau)| \leq \frac{\varepsilon^2}{8} \bigg| \int_0^\tau \| \partial_x^{k + 1} U_\varepsilon \|_{H^1(\R)} \| \partial_x^k Z_\varepsilon \|_{L^2(\R)} \big( \| G_\varepsilon \|_{L^2(\R)} + \varepsilon^2 \| R_\varepsilon \|_{L^2(\R)} \big) \bigg|,$$
so that \eqref{estim-J2} follows from \eqref{frites} and \eqref{allfolks1}.

Concerning the last integral $J_3(\tau)$, we finally derive

\begin{step}
\label{Estim-J3}
Under the assumptions of Lemma \ref{Jesuisborne}, given any $\tau \in \R$, there exist a positive number $R_1$, depending on $\varepsilon$ and $V_\varepsilon$, and a positive constant $K$, depending only on $K_0$ and $k$, such that
\begin{equation}
\label{estim-J3}
|J_3(\tau)| \leq K \varepsilon^2 \bigg( \varepsilon^2 \exp K |\tau| + \bigg| \int_0^\tau \exp K |s| \ \| \partial_x^k Z_\varepsilon(\cdot, s) \|_{L^2(\R)} ds \bigg| \bigg),
\end{equation}
for any choice of the number $R$ of definition \eqref{defupsilon} in $(R_1, + \infty)$.
\end{step}

Combining estimates \eqref{frites} and \eqref{allfolks1} with definition \eqref{J3}, we have
\begin{equation}
\label{erbate}
\begin{split}
|J_3(\tau)| \leq & \bigg| \int_0^\tau \Big( |V_\varepsilon(- R, \cdot)| + \varepsilon^2 \| G_\varepsilon \|_{L^\infty(\R)} + \varepsilon^4 \| R_\varepsilon \|_{L^\infty(\R)} \Big) \| \partial_x^{k + 1} U_\varepsilon \|_{L^2(\R)} \| \partial_x^k Z_\varepsilon \|_{L^2(\R)} \bigg|\\
\leq & K \bigg| \int_0^\tau \Big( |V_\varepsilon(- R, s)| + \varepsilon^2 \Big) \exp K |s| \ \| \partial_x^k Z_\varepsilon(\cdot, s) \|_{L^2(\R)} ds \bigg|.
\end{split}
\end{equation}
We then invoke \eqref{allfolks3} to write
$$\bigg| \int_0^\tau |V_\varepsilon(- R, s)| \exp K |s| \ \| \partial_x^k Z_\varepsilon(\cdot, s) \|_{L^2(\R)} ds \bigg| \leq K' \exp K' |\tau| \ \bigg| \int_0^\tau |V_\varepsilon(- R, s)|^2 ds \bigg|^\frac{1}{2},$$
where $K'$ is some further constant depending only on $K_0$ and $k$. Next, we have
$$|V_\varepsilon(- R, s)|^2 \leq \int_{- \infty}^{- R} \Big( \big( V_\varepsilon(x, s) \big)^2 + \big( \partial_x V_\varepsilon(x, s) \big)^2 \Big) dx,$$
so that
\begin{equation}
\label{givet}
\begin{split}
& \bigg| \int_0^\tau |V_\varepsilon(- R, s)| \exp K |s| \ \| \partial_x^k Z_\varepsilon(\cdot, s) \|_{L^2(\R)} ds \bigg|\\
\leq K' \exp & K' |\tau| \ \bigg| \int_0^\tau \int_{- \infty}^{- R} \Big( \big( V_\varepsilon(x, s) \big)^2 + \big( \partial_x V_\varepsilon(x, s) \big)^2 \Big) dx ds \bigg|^\frac{1}{2}.
\end{split}
\end{equation}
By the dominated convergence theorem, the integral on the right-hand side of \eqref{givet} decays to $0$ as $R \to + \infty$, so that there exists a positive number $R_1$, depending on $\varepsilon$ and $V_\varepsilon$, such that
\begin{equation}
\label{pape}
\bigg| \int_0^\tau \int_{- \infty}^{- R} \Big( \big( V_\varepsilon(x, s) \big)^2 + \big( \partial_x V_\varepsilon(x, s) \big)^2 \Big) dx ds \bigg| \leq \varepsilon^8,
\end{equation}
for any $R > R_1$. Combining with \eqref{erbate} and \eqref{givet}, this provides \eqref{estim-J3}.

We are now in position to complete the proof of Lemma \ref{Jesuisborne}.

\begin{proof}[End of the proof of Lemma \ref{Jesuisborne}]
Estimate \eqref{oheme} follows from applying bounds \eqref{estim-J1}, \eqref{estim-J2} and \eqref{estim-J3} to definition \eqref{taiwo}.
\end{proof}

\subsection{Proof of Lemma \ref{Ture}}

Though the proof is similar to the proofs of Lemmas \ref{Itisbounded} and \ref{Jesuisborne}, the computations are somewhat more technical. In view of definition \eqref{fF} and the fact that $\partial_x \Upsilon_\varepsilon = V_\varepsilon$, we first decompose the quantity $K_\varepsilon(\tau)$ as
\begin{equation}
\label{rame}
K_\varepsilon(\tau) = \frac{\varepsilon^2}{24} \sum_{j = 0}^k \binom{k}{j} \Big( \frac{1}{6} K_1^j(\tau) - K_2^j(\tau) + \frac{1}{3} K_3^j(\tau) \Big),
\end{equation}
where
\begin{equation}
\label{K1j}
K_1^j(\tau) = \int_0^\tau \int_\R \partial_x^{k + 1} \big( V_\varepsilon^2 \big) \partial_x^j \Upsilon_\varepsilon \partial_x^{k + 1 - j} U_\varepsilon,
\end{equation}
\begin{equation}
\label{K2j}
K_2^j(\tau) = \int_0^\tau \int_\R \partial_x^{k + 3} V_\varepsilon \partial_x^j \Upsilon_\varepsilon \partial_x^{k + 1 - j} U_\varepsilon,
\end{equation}
and
\begin{equation}
\label{K3j}
K_3^j(\tau) = \int_0^\tau \int_\R \partial_x^{k + 1} \big( U_\varepsilon V_\varepsilon \big) \partial_x^j \Upsilon_\varepsilon \partial_x^{k + 1 - j} U_\varepsilon.
\end{equation}
We then estimate each of the above integrals $K_i^j(\tau)$.

\begin{claim}
\label{Estim-K}
Let $i \in \{ 1, 2, 3 \}$ and $0 \leq j \leq k$. Under the assumptions of Lemma \ref{Ture}, given any $\tau \in \R$, there exists a positive constant $K$, depending only on $K_0$ and $k$, such that
\begin{equation}
\label{estim-K}
|K_i^j(\tau)| \leq K \varepsilon^2 \exp K |\tau|,
\end{equation}
for any choice of the number $R$ of definition \eqref{defupsilon} in $(R_1, + \infty)$, where $R_1$ denotes the positive number given by Lemma \ref{Jesuisborne}.
\end{claim}

Estimate \eqref{moisi} follows from combining decomposition \eqref{rame} with bounds \eqref{estim-K}, so that Lemma \ref{Ture} is a direct consequence of Claim \ref{Estim-K}, and it only remains to show Claim \ref{Estim-K}.

\begin{proof}[Proof of Claim \ref{Estim-K}]
We split the proof in six cases according to the values of $i$ and $j$.

\begin{case}
\label{Estim-K1j}
$i = 1$ and $1 \leq j \leq k$.
\end{case}

Using the fact that $\partial_x \Upsilon_\varepsilon = V_\varepsilon$, we have in view of definition \eqref{K1j},
$$K_1^j(\tau) = \int_0^\tau \int_\R \partial_x^{k + 1} \big( V_\varepsilon^2 \big) \partial_x^{j - 1} V_\varepsilon \partial_x^{k + 1 - j} U_\varepsilon,$$
so that integrating by parts in space
\begin{equation}
\label{K1jbis}
K_1^j(\tau) = - \int_0^\tau \int_\R \partial_x \Big( \partial_x^{k + 1} \big( V_\varepsilon^2 \big) \partial_x^{j - 1} V_\varepsilon \Big) \partial_x^{k - j} U_\varepsilon.
\end{equation}
In view of \eqref{slow2} and \eqref{eq:Vcarre}, we now have
\begin{align*}
\partial_x \Big( \partial_x^{k + 1} \big( V_\varepsilon^2 \big) \partial_x^{j - 1} V_\varepsilon \Big)& = - \frac{\varepsilon^2}{8} \partial_\tau \Big( \partial_x^{k + 1} \big( V_\varepsilon^2 \big) \partial_x^{j - 1} V_\varepsilon \Big) + \frac{\varepsilon^2}{4} \partial_x^{j - 1} V_\varepsilon \Big( \partial_x^{k + 1} \big( g_\varepsilon V_\varepsilon \big) + \varepsilon^2 \partial_x^{k + 1} \big( r_\varepsilon V_\varepsilon \big) \Big)\\
& + \frac{\varepsilon^2}{8} \partial_x^{k + 1} \big( V_\varepsilon^2 \big) \Big( \partial_x^{j - 1} g_\varepsilon + \varepsilon^2 \partial_x^{j - 1} r_\varepsilon \Big),
\end{align*}
so that \eqref{K1jbis} becomes after an integration by parts in time,
\begin{equation}
\label{henrique}
\begin{split}
K_1^j(\tau) & = \frac{\varepsilon^2}{8} \bigg[ \int_\R \partial_x^{k - j} U_\varepsilon \partial_x^{k + 1} \big( V_\varepsilon^2 \big) \partial_x^{j - 1} V_\varepsilon \bigg]_0^\tau - \frac{\varepsilon^2}{8} \int_0^\tau \int_\R \bigg( \partial_\tau \partial_x^{k - j} U_\varepsilon \partial_x^{k + 1} \big( V_\varepsilon^2 \big) \partial_x^{j - 1} V_\varepsilon\\
+ & 2 \partial_x^{k - j} U_\varepsilon \partial_x^{j - 1} V_\varepsilon \Big( \partial_x^{k + 1} \big( g_\varepsilon V_\varepsilon \big) + \varepsilon^2 \partial_x^{k + 1} \big( r_\varepsilon V_\varepsilon \big) \Big) + \partial_x^{k - j} U_\varepsilon \partial_x^{k + 1} \big( V_\varepsilon^2 \big) \Big( \partial_x^{j - 1} g_\varepsilon + \varepsilon^2 \partial_x^{j - 1} r_\varepsilon \Big) \bigg).
\end{split}
\end{equation}
We finally argue as in the proof of Lemma \ref{Itisbounded} using bounds \eqref{frites} and \eqref{allfolks1} to bound any term on the right-hand side of \eqref{henrique}. This provides estimate \eqref{estim-K} in case $i = 1$ and $1 \leq j \leq k$.

\begin{case}
\label{K1k}
$i = 1$ and $j = 0$.
\end{case}

When $j = 0$, we also integrate by parts in space to obtain
\begin{equation}
\label{K10}
K_1^0(\tau) = - \int_0^\tau \int_\R \partial_x \big( \partial_x^{k + 1} \big( V_\varepsilon^2 \big) \Upsilon_\varepsilon \big) \partial_x^{k - j} U_\varepsilon.
\end{equation}
We then combine \eqref{eq:Vcarre} with \eqref{equpsilon} to establish
\begin{align*}
\partial_x & \Big( \partial_x^{k + 1} \big( V_\varepsilon^2 \big) \Upsilon_\varepsilon \Big) = - \frac{\varepsilon^2}{8} \partial_\tau \Big( \partial_x^{k + 1} \big( V_\varepsilon^2 \big) \Upsilon_\varepsilon \Big) + \frac{\varepsilon^2}{4} \Upsilon_\varepsilon \Big( \partial_x^{k + 1} \big( g_\varepsilon V_\varepsilon \big) + \varepsilon^2 \partial_x^{k + 1} \big( r_\varepsilon V_\varepsilon \big) \Big)\\
& + \frac{\varepsilon^2}{8} \partial_x^{k + 1} \big( V_\varepsilon^2 \big) \Big( G_\varepsilon + \varepsilon^2 R_\varepsilon \Big) - \frac{\varepsilon^2}{8} \partial_x^{k + 1} \big( V_\varepsilon^2 \big) \Big( G_\varepsilon(- R, \cdot) + \varepsilon^2 R_\varepsilon(- R, \cdot) \Big) + V_\varepsilon(- R, \cdot) \partial_x^{k + 1} \big( V_\varepsilon^2 \big),
\end{align*}
so that \eqref{K10} becomes after an integration by parts in time,
\begin{equation}
\label{diarra}
\begin{split}
K_1^0(\tau) = & \frac{\varepsilon^2}{8} \bigg[ \int_\R \partial_x^{k - j} U_\varepsilon \partial_x^{k + 1} \big( V_\varepsilon^2 \big) \Upsilon_\varepsilon \bigg]_0^\tau - \frac{\varepsilon^2}{8} \int_0^\tau \int_\R \partial_\tau \partial_x^{k - j} U_\varepsilon \partial_x^{k + 1} \big( V_\varepsilon^2 \big) \Upsilon_\varepsilon\\
& - \frac{\varepsilon^2}{8} \int_0^\tau \int_\R \partial_x^{k - j} U_\varepsilon \bigg( 2 \Upsilon_\varepsilon \Big( \partial_x^{k + 1} \big( g_\varepsilon V_\varepsilon \big) + \varepsilon^2 \partial_x^{k + 1} \big( r_\varepsilon V_\varepsilon \big) \Big) + \partial_x^{k + 1} \big( V_\varepsilon^2 \big) \Big( G_\varepsilon + \varepsilon^2 R_\varepsilon \Big) \bigg)\\
& - \int_0^\tau \Big( V_\varepsilon(- R, \cdot) - \frac{\varepsilon^2}{8} G_\varepsilon(- R, \cdot) -\frac{\varepsilon^4}{8} R_\varepsilon(- R, \cdot) \Big) \int_\R \partial_x^{k - j} U_\varepsilon \partial_x^{k + 1} \big( V_\varepsilon^2 \big).
\end{split}
\end{equation}
We then estimate any terms on the right-hand side of \eqref{diarra} similarly to the terms on the right-hand side of \eqref{henrique}. Using bounds \eqref{frites}, \eqref{allfolks1} and \eqref{cana}, we are led to
$$|K_1^0(\tau)| \leq K \bigg( \varepsilon^2 \exp K |\tau| + \bigg| \int_0^\tau \Big( |V_\varepsilon(- R, s)| + \varepsilon^2 \| G_\varepsilon(\cdot, s) \|_{L^\infty} + \varepsilon^4 \| R_\varepsilon(\cdot, s) \|_{L^\infty} \Big) \exp K |s| \ ds \bigg| \bigg),$$
so that by the Sobolev embedding theorem and \eqref{allfolks1},
\begin{equation}
\label{cavenaghi}
|K_1^0(\tau)| \leq K \bigg( \varepsilon^2 \exp K |\tau| + \bigg| \int_0^\tau |V_\varepsilon(- R, s)| \exp K |s| \ ds \bigg| \bigg). 
\end{equation}
We then argue as in the proof of Step \ref{Estim-J3} of Lemma \ref{Jesuisborne}. Given any $R > R_1$, we obtain
\begin{equation}
\label{wendel}
\bigg| \int_0^\tau |V_\varepsilon(- R, s)| \exp K |s| \ ds \bigg| \leq K' \exp K' |\tau| \ \bigg| \int_0^\tau \int_{- \infty}^{- R} \Big( V_\varepsilon^2 + \big( \partial_x V_\varepsilon \big)^2 \Big) \bigg|^\frac{1}{2} \leq K' \varepsilon^4 \exp K' |\tau|,
\end{equation}
where $K'$ is some further positive constant depending only on $K_0$ and $k$. Combining with \eqref{cavenaghi}, this completes the proof of \eqref{estim-K} in case $i = 1$ and $j = 0$.

\begin{case}
\label{Estim-K2j}
$i = 2$ and $1 \leq j \leq k$.
\end{case}

The proof is similar to Case \ref{Estim-K1j}. Using the fact that $\partial_x \Upsilon_\varepsilon = V_\varepsilon$, and integrating by parts the derivative $\partial_x^{k + 3} V_\varepsilon$, we have in view of definition \eqref{K2j},
\begin{align*}
K_2^j(\tau) = & \int_0^\tau \int_\R \partial_x^{k + 1} V_\varepsilon \partial_x^{j + 1} V_\varepsilon \partial_x^{k + 1 - j} U_\varepsilon + 2 \int_0^\tau \int_\R \partial_x^{k + 1} V_\varepsilon \partial_x^j V_\varepsilon \partial_x^{k + 2 - j} U_\varepsilon\\
+ & \int_0^\tau \int_\R \partial_x^{k + 1} V_\varepsilon \partial_x^{j - 1} V_\varepsilon \partial_x^{k + 3 - j} U_\varepsilon,
\end{align*}
so that integrating by parts the derivatives of $U_\varepsilon$,
\begin{equation}
\label{K2jbis}
\begin{split}
K_2^j(\tau) = & - \int_0^\tau \int_\R \partial_x \Big( \partial_x^{k + 1} V_\varepsilon \partial_x^{j + 1} V_\varepsilon \Big) \partial_x^{k - j} U_\varepsilon - 2 \int_0^\tau \int_\R \partial_x \Big( \partial_x^{k + 1} V_\varepsilon \partial_x^j V_\varepsilon \Big) \partial_x^{k + 1 - j} U_\varepsilon\\
- & \int_0^\tau \int_\R \partial_x \Big( \partial_x^{k + 1} V_\varepsilon \partial_x^{j - 1} V_\varepsilon \Big) \partial_x^{k + 2 - j} U_\varepsilon.
\end{split}
\end{equation}
In view of \eqref{slow2}, we now remark that
\begin{equation}
\label{planus}
\begin{split}
\partial_x \Big( \partial_x^l V_\varepsilon \partial_x^m V_\varepsilon \Big) & = - \frac{\varepsilon^2}{8} \partial_\tau \Big( \partial_x^l V_\varepsilon \partial_x^m V_\varepsilon \Big) + \frac{\varepsilon^2}{8} \Big( \partial_x^l V_\varepsilon \partial_x^m g_\varepsilon + \partial_x^m V_\varepsilon \partial_x^l g_\varepsilon \Big)\\
& + \frac{\varepsilon^4}{8} \Big( \partial_x^l V_\varepsilon \partial_x^m r_\varepsilon + \partial_x^m V_\varepsilon \partial_x^l r_\varepsilon \Big),
\end{split}
\end{equation}
for any $(l, m) \in \N^2$. We then introduce the expression of the functions $\partial_x (\partial_x^l V_\varepsilon \partial_x^m V_\varepsilon)$ given by \eqref{planus} in the three integrals on the right-hand side of \eqref{K2jbis}, integrate by parts in time as in the proof of Case \ref{Estim-K1j}, and bound the resulting terms using estimates \eqref{frites} and \eqref{allfolks1}. This provides inequality \eqref{estim-K} in case $i = 2$ and $1 \leq j \leq k$.

\begin{case}
\label{K2k}
$i = 2$ and $j = 0$.
\end{case}

The proof is similar to Cases \ref{K1k} and \ref{Estim-K2j}. We first integrate by parts as in \eqref{K2jbis} to obtain 
\begin{equation}
\label{K20}
\begin{split}
K_2^0(\tau) = & - \int_0^\tau \int_\R \partial_x \big( \partial_x^{k + 1} V_\varepsilon \partial_x V_\varepsilon \big) \partial_x^{k - j} U_\varepsilon - 2 \int_0^\tau \int_\R \partial_x \big( \partial_x^{k + 1} V_\varepsilon V_\varepsilon \big) \partial_x^{k + 1 - j} U_\varepsilon\\
& - \int_0^\tau \int_\R \partial_x \big( \partial_x^{k + 1} V_\varepsilon \Upsilon_\varepsilon \big) \partial_x^{k + 2 - j} U_\varepsilon.
\end{split}
\end{equation}
Concerning the first and the second integrals on the right-hand side of \eqref{K20}, we then replace the functions $\partial_x (\partial_x^{k + 1} V_\varepsilon V_\varepsilon)$ and $\partial_x (\partial_x^{k + 1} V_\varepsilon \partial_x V_\varepsilon)$ by their expression given by \eqref{planus}, integrate by parts in time as in the proof of Case \ref{Estim-K1j}, and bound the resulting terms using estimates \eqref{frites} and \eqref{allfolks1}. This provides
\begin{equation}
\label{bellion}
\bigg| \int_0^\tau \int_\R \partial_x \big( \partial_x^{k + 1} V_\varepsilon \partial_x V_\varepsilon \big) \partial_x^{k - j} U_\varepsilon \bigg| + \bigg| \int_0^\tau \int_\R \partial_x \big( \partial_x^{k + 1} V_\varepsilon V_\varepsilon \big) \partial_x^{k + 1 - j} U_\varepsilon \bigg| \leq K \varepsilon^2 \exp K |\tau|.
\end{equation}
In contrast, for the last integral on the right-hand side of \eqref{K20}, we combine \eqref{slow2} with \eqref{equpsilon} to establish
\begin{align*}
\partial_x \Big( \partial_x^{k + 1} V_\varepsilon \Upsilon_\varepsilon \Big) = - & \frac{\varepsilon^2}{8} \partial_\tau \Big( \partial_x^{k + 1} V_\varepsilon \Upsilon_\varepsilon \Big) + \frac{\varepsilon^2}{8} \Upsilon_\varepsilon \Big( \partial_x^{k + 1} g_\varepsilon + \varepsilon^2 \partial_x^{k + 1} r_\varepsilon \Big) + \frac{\varepsilon^2}{8} \partial_x^{k + 1} V_\varepsilon \Big( G_\varepsilon + \varepsilon^2 R_\varepsilon \Big)\\
- & \frac{\varepsilon^2}{8} \partial_x^{k + 1} V_\varepsilon \Big( G_\varepsilon(- R, \cdot) + \varepsilon^2 R_\varepsilon(- R, \cdot) \Big) + V_\varepsilon(- R, \cdot) \partial_x^{k + 1} V_\varepsilon,
\end{align*}
so that after an integration by parts in time,
\begin{align*}
& \int_0^\tau \int_\R \partial_x \big( \partial_x^{k + 1} V_\varepsilon \Upsilon_\varepsilon \big) \partial_x^{k + 2 - j} U_\varepsilon = - \frac{\varepsilon^2}{8} \bigg[ \int_\R \partial_x^{k + 2 - j} U_\varepsilon \partial_x^{k + 1} V_\varepsilon \Upsilon_\varepsilon \bigg]_0^\tau\\
+ \frac{\varepsilon^2}{8} \int_0^\tau \int_\R & \bigg( \partial_\tau \partial_x^{k + 2 - j} U_\varepsilon \partial_x^{k + 1} V_\varepsilon \Upsilon_\varepsilon + \partial_x^{k + 2 - j} U_\varepsilon \Big( \Upsilon_\varepsilon \big( \partial_x^{k + 1} g_\varepsilon + \varepsilon^2 \partial_x^{k + 1} r_\varepsilon \big)\\
+ \partial_x^{k + 1} V_\varepsilon \big( G_\varepsilon & + \varepsilon^2 R_\varepsilon \big) \Big) \bigg) + \int_0^\tau \Big( V_\varepsilon(- R, \cdot) - \frac{\varepsilon^2}{8} G_\varepsilon(- R, \cdot) -\frac{\varepsilon^4}{8} R_\varepsilon(- R, \cdot) \Big) \int_\R \partial_x^{k + 2 - j} U_\varepsilon \partial_x^{k + 1} V_\varepsilon.
\end{align*}
We finally argue as in the proof of Case \ref{K1k}. Using bounds \eqref{frites}, \eqref{allfolks1}, \eqref{cana} and \eqref{wendel}, we are led to
$$\bigg| \int_0^\tau \int_\R \partial_x \big( \partial_x^{k + 1} V_\varepsilon \Upsilon_\varepsilon \big) \partial_x^{k + 2 - j} U_\varepsilon \bigg| \leq K \varepsilon^2 \exp K |\tau|,$$
so that \eqref{estim-K} follows for $i = 2$ and $j = 0$ from \eqref{K20} and \eqref{bellion}.

\begin{case}
\label{Estim-K3j}
$i = 3$ and $1 \leq j \leq k$.
\end{case}

The proof is similar to Cases \ref{Estim-K1j} and \ref{Estim-K2j}. Using the fact that $\partial_x \Upsilon_\varepsilon = V_\varepsilon$ and applying the Newton formula to \eqref{K3j}, we are led to
\begin{equation}
\label{K3jl}
K_3^j(\tau) = \sum_{l = 0}^{k + 1} \binom{k + 1}{l} \boK_l^j.
\end{equation}
where the integrals $\boK_l^j$ are given by
\begin{equation}
\label{bok}
\boK_l^j \equiv \int_0^\tau \int_\R \partial_x^l V_\varepsilon \partial_x^{j - 1} V_\varepsilon \partial_x^{k + 1 - l} U_\varepsilon \partial_x^{k + 1 - j} U_\varepsilon.
\end{equation}
Assuming first that $l \geq j$, our arguments to estimate the integrals $\boK_l^j$ then depend on the parity of the difference $l - j$. When $l -j = 2 m$ is even, we can write
$$\partial_x^l V_\varepsilon \partial_x^{j - 1} V_\varepsilon = \frac{(- 1)^m}{2} \partial_x \Big( \big( \partial_x^{j + m - 1} V_\varepsilon \big)^2 \Big) + \sum_{p = 0}^{m - 1} (- 1)^p \partial_x \Big( \partial_x^{l - 1 - p} V_\varepsilon \partial_x^{j - 1 + p} V_\varepsilon \Big),$$
so that \eqref{bok} becomes
\begin{align*}
\boK_l^j = & \frac{(- 1)^m}{2} \int_0^\tau \int_\R \partial_x \Big( \big( \partial_x^{j + m - 1} V_\varepsilon \big)^2 \Big) \partial_x^{k + 1 - l} U_\varepsilon \partial_x^{k + 1 - j} U_\varepsilon\\
+ & \sum_{p = 0}^{m - 1} (- 1)^p \int_0^\tau \int_\R \partial_x \Big( \partial_x^{l - 1 - p} V_\varepsilon \partial_x^{j - 1 + p} V_\varepsilon \Big) \partial_x^{k + 1 - l} U_\varepsilon \partial_x^{k + 1 - j} U_\varepsilon.
\end{align*}
We then replace the functions $\partial_x \Big( \big( \partial_x^{j + m - 1} V_\varepsilon \big)^2 \Big)$ and $\partial_x \Big( \partial_x^{l - 1 - p} V_\varepsilon \partial_x^{j - 1 + p} V_\varepsilon \Big)$ by their expressions given by \eqref{planus} and argue as in the proof of Cases \ref{Estim-K1j} and \ref{Estim-K2j} to obtain that
\begin{equation}
\label{jurietti}
|\boK_l^j| \leq K \varepsilon^2 \exp K |\tau|.
\end{equation}

In contrast, when $l -j = 2 m + 1$ is odd, we can write
\begin{equation}
\label{riera}
\partial_x^{k + 1 - l} U_\varepsilon \partial_x^{k + 1 - j} U_\varepsilon = \frac{(- 1)^m}{2} \partial_x \Big( \big( \partial_x^{k - j - m} U_\varepsilon \big)^2 \Big) + \sum_{p = 0}^{m - 1} (- 1)^p \partial_x \Big( \partial_x^{k - j - p} U_\varepsilon \partial_x^{k - l + 1 + p} U_\varepsilon \Big),
\end{equation}
so that \eqref{bok} becomes after an integration by parts in space,
\begin{align*}
\boK_l^j = & \frac{(- 1)^{m + 1}}{2} \int_0^\tau \int_\R \partial_x \Big( \big( \partial_x^{j - 1} V_\varepsilon \partial_x^l V_\varepsilon \Big) \big( \partial_x^{k - j - m} U_\varepsilon \big)^2\\
- & \sum_{p = 0}^{m - 1} (- 1)^p \int_0^\tau \int_\R \partial_x \Big( \partial_x^{j - 1} V_\varepsilon \partial_x^l V_\varepsilon \Big) \partial_x^{k - j - p} U_\varepsilon \partial_x^{k - l + 1 + p} U_\varepsilon.
\end{align*}
We then complete the proof of \eqref{jurietti} as in the case $l - j$ was even. The proof is similar when $j > l$ permuting the roles of $j$ and $l$. Hence, \eqref{jurietti} holds for any choice of $j$ and $k$. In view of \eqref{K3jl}, this concludes the proof of \eqref{estim-K} when $i = 3$ and $1 \leq j \leq k$.

\begin{case}
\label{K3k}
$i = 3$ and $j = 0$.
\end{case}

The proof is similar to Case \ref{Estim-K3j}. Applying the Leibniz formula to \eqref{K3j}, we are led to
\begin{equation}
\label{K30l}
K_3^0(\tau) = \sum_{l = 0}^{k + 1} \binom{k + 1}{l} \boK_l^0.
\end{equation}
where the integrals $\boK_l^0$ are given by
\begin{equation}
\label{bok0}
\boK_l^0 \equiv \int_0^\tau \int_\R \partial_x^l V_\varepsilon \Upsilon_\varepsilon \partial_x^{k + 1 - l} U_\varepsilon \partial_x^{k + 1} U_\varepsilon.
\end{equation}
When $l = 2 m$ is even and positive, we can write using the fact that $\partial_x \Upsilon_\varepsilon = V_\varepsilon$,
\begin{equation}
\label{mavuba}
\partial_x^l V_\varepsilon \Upsilon_\varepsilon = \frac{(- 1)^m}{2} \partial_x \Big( \big( \partial_x^{m - 1} V_\varepsilon \big)^2 \Big) + \partial_x \Big( \partial_x^{l - 1} V_\varepsilon \Upsilon_\varepsilon \Big) + \sum_{p = 1}^{m - 1} (- 1)^p \partial_x \Big( \partial_x^{l - 1 - p} V_\varepsilon \partial_x^{p - 1} V_\varepsilon \Big),
\end{equation}
so that \eqref{bok0} becomes
\begin{equation}
\label{micoud}
\begin{split}
\boK_l^0 = & \frac{(- 1)^m}{2} \int_0^\tau \int_\R \partial_x \Big( \big( \partial_x^{m - 1} V_\varepsilon \big)^2 \Big) \partial_x^{k + 1 - l} U_\varepsilon \partial_x^{k + 1} U_\varepsilon + \int_0^\tau \int_\R \partial_x \Big( \partial_x^{l - 1} V_\varepsilon \Upsilon_\varepsilon \Big) \partial_x^{k + 1 - l} U_\varepsilon \partial_x^{k + 1} U_\varepsilon\\
+ & \sum_{p = 1}^{m - 1} (- 1)^p \int_0^\tau \int_\R \partial_x \Big( \partial_x^{l - 1 - p} V_\varepsilon \partial_x^{p - 1} V_\varepsilon \Big) \partial_x^{k + 1 - l} U_\varepsilon \partial_x^{k + 1} U_\varepsilon.
\end{split}
\end{equation}
We then rely on the arguments of the proofs of Cases \ref{K1k} and \ref{K2k} to bound any integrals on the right-hand side of \eqref{micoud} in order to obtain
\begin{equation}
\label{blanc}
|\boK_l^0| \leq K \varepsilon^2 \exp K |\tau|.
\end{equation}
When $l = 2 m + 1$, we invoke \eqref{riera} as in the proof of Case \ref{Estim-K3j}, instead of \eqref{mavuba}, and complete the proof of \eqref{blanc} similarly.

Finally, for $l = 0$, we write
\begin{equation}
\label{K00}
\boK_0^0 = \frac{1}{2} \int_0^\tau \int_\R \partial_x \big( \Upsilon_\varepsilon^2 \big) \big( \partial_x^{k + 1} U_\varepsilon \big)^2.
\end{equation}
In view of \eqref{equpsilon}, we have
$$\partial_x \big( \Upsilon_\varepsilon^2 \big) = - \frac{\varepsilon^2}{8} \partial_\tau \big( \Upsilon_\varepsilon^2 \big) + \frac{\varepsilon^2}{4} \Upsilon_\varepsilon \Big( G_\varepsilon + \varepsilon^2 R_\varepsilon - G_\varepsilon(- R, \cdot) - \varepsilon^2 R_\varepsilon(- R, \cdot) \Big) + 2 V_\varepsilon(- R, \cdot) \Upsilon_\varepsilon,$$
so that integrating by parts in time, \eqref{K00} becomes
\begin{equation}
\label{meriem}
\begin{split}
\boK_0^0 = & - \frac{\varepsilon^2}{16} \bigg[ \int_\R \Upsilon_\varepsilon^2 \big( \partial_x^{k + 1} U_\varepsilon \big)^2 \bigg]_0^\tau + \frac{\varepsilon^2}{8} \int_0^\tau \int_\R \Upsilon_\varepsilon^2 \partial_x^{k + 1} U_\varepsilon \partial_\tau \partial_x^{k + 1} U_\varepsilon + \int_0^\tau V_\varepsilon(- R, \cdot) \int_\R \Upsilon_\varepsilon \big( \partial_x^{k + 1} U_\varepsilon \big)^2\\
+ & \frac{\varepsilon^2}{8} \int_0^\tau \int_\R \Upsilon_\varepsilon \Big( G_\varepsilon + \varepsilon^2 R_\varepsilon - G_\varepsilon(- R, \cdot) - \varepsilon^2 R_\varepsilon(- R, \cdot) \Big) \big( \partial_x^{k + 1} U_\varepsilon \big)^2.
\end{split}
\end{equation}
Using bounds \eqref{frites}, \eqref{allfolks1}, \eqref{cana} and \eqref{wendel} to estimate the right-hand side of \eqref{meriem}, we deduce \eqref{blanc} for $l = 0$ provided we choose $R > R_1$ as above. Combining with \eqref{K30l}, this completes the proof of \eqref{estim-K} in case $i = 3$ and $j = 0$, and therefore of Claim \ref{Estim-K}.
\end{proof}

\section{Proof of Theorem \ref{cochondore}}
\label{quatre}

We now turn to the proof of Theorem \ref{cochondore}. As mentioned in the introduction, we focus on the coordinate $x = x^- = \varepsilon (\x + \sqrt{2} t)$ and the associated functions $U_\varepsilon = U_\varepsilon^-$ and $\boU = \boU^-$. Given any $k \in \N$, we first recall identity \eqref{gege}, which may be written as
\begin{align*}
\partial_\tau \boZ_\varepsilon^k(\tau) = & - 2 \int_0^\tau \int_\R \partial_x^{k+1} \big( \boU Z_\varepsilon \big) \partial_x^k Z_\varepsilon - \int_0^\tau \int_\R \partial_x^{k+1} \big( Z_\varepsilon^2 \big) \partial_x^k Z_\varepsilon\\
+ & 2 \int_0^\tau \int_\R \partial_x^k f_\varepsilon \partial_x^k Z_\varepsilon - 2 \varepsilon^2 \int_0^\tau \int_\R \partial_x^k r_\varepsilon \partial_x^k Z_\varepsilon,
\end{align*}
with $Z_\varepsilon = U_\varepsilon - \boU$. We then bound inductively any term on the right-hand side of \eqref{gege} in order to apply a Gronwall lemma to the quantity $\boZ_\varepsilon^k(\tau)$ defined by \eqref{beauzire}, and derive inequality \eqref{ineq1}.

For $k = 0$, identity \eqref{gege} may be recast as
\begin{equation}
\label{krilin}
\partial_\tau \boZ_\varepsilon^0(\tau) = - \int_0^\tau \int_\R \partial_x \boU Z_\varepsilon^2 + 2 \int_0^\tau \int_\R f_\varepsilon Z_\varepsilon - 2 \varepsilon^2 \int_0^\tau \int_\R r_\varepsilon Z_\varepsilon,
\end{equation}
so that in view of Proposition \ref{Estim-f}, and bounds \eqref{allfolks1} and \eqref{allfolks2}, we have
$$\partial_\tau \boZ_\varepsilon^0(\tau) \leq K \bigg( |\boZ_\varepsilon^0(\tau)| + \varepsilon^2 \bigg| \int_0^\tau \exp K |s| \ \| Z_\varepsilon(\cdot, s) \|_{L^2(\R)} ds \bigg| + \varepsilon^2 \big( \| Z_\varepsilon(\cdot, \tau) \|_{L^2(\R)} + \varepsilon^2 \big) \exp K |\tau| \bigg),$$
where $K$ is a positive constant depending only on $K_0$. Using the inequality $2 |a b| \leq a^2 + b^2$ and the identity
\begin{equation}
\label{derivz0}
\partial_\tau \boZ_\varepsilon^0(\tau) = \| Z_\varepsilon(\cdot, \tau) \|_{L^2(\R)}^2,
\end{equation}
we are led to the differential inequality
\begin{equation}
\label{gronwall0}
\partial_\tau \boZ_\varepsilon^0(\tau) \leq K \Big( \sign(\tau) \boZ_\varepsilon^0(\tau) + \varepsilon^4 \exp K |\tau| \Big),
\end{equation}
so that
$$|\boZ_\varepsilon^0(\tau)| \leq K' \varepsilon^4 \exp K' |\tau|,$$
where $K'$ is some further positive constant depending only on $K_0$. Combining with \eqref{derivz0} and \eqref{gronwall0}, this provides \eqref{ineq1} for $k = 0$ (and the functions $U_\varepsilon^-$ and $\boU^-$).

We now assume that \eqref{ineq1} holds for any $0 \leq j \leq k - 1$, i.e. that there exists a positive constant $K$ depending only on $K_0$ and $k$, such that
\begin{equation}
\label{quirecure}
\| Z_\varepsilon(\cdot, \tau) \|_{H^{k - 1}(\R)} \leq K \varepsilon^2 \exp K |\tau|,
\end{equation}
for any $\tau \in \R$. We then bound any integral on the right-hand side of \eqref{gege}. For the first one, we compute by the Leibniz formula,
$$\int_0^\tau \int_\R \partial_x^{k+1} \big( \boU Z_\varepsilon \big) \partial_x^k Z_\varepsilon = \sum_{j = 0}^{k - 1} \binom{k + 1}{j} \int_0^\tau \int_\R \partial_x^{k + 1 - j} \boU \partial_x^j Z_\varepsilon \partial_x^k Z_\varepsilon + \Big( k + \frac{1}{2} \Big) \int_0^\tau \int_\R \partial_x \boU \big( \partial_x^k Z_\varepsilon \big)^2,$$
so that by \eqref{allfolks2} and \eqref{quirecure}, we are led to
\begin{equation}
\label{foncia}
\bigg| \int_0^\tau \int_\R \partial_x^{k+1} \big( \boU Z_\varepsilon \big) \partial_x^k Z_\varepsilon \bigg| \leq K \bigg( \bigg| \int_0^\tau \int_\R \big( \partial_x^k Z_\varepsilon \big)^2 \bigg| + \varepsilon^2 \bigg| \int_0^\tau \exp K |s| \ \| \partial_x^k Z_\varepsilon(\cdot, s) \|_{L^2(\R)} ds \bigg| \bigg).
\end{equation}
The estimates of the second integral on the right-hand side of \eqref{gege} are similar. The Leibniz formula yields
$$\int_0^\tau \int_\R \partial_x^{k+1} \big( Z_\varepsilon^2 \big) \partial_x^k Z_\varepsilon = \sum_{j = 2}^{k - 1} \binom{k + 1}{j} \int_0^\tau \int_\R \partial_x^{k + 1 - j} Z_\varepsilon \partial_x^j Z_\varepsilon \partial_x^k Z_\varepsilon + (2 k + 1) \int_0^\tau \int_\R \partial_x Z_\varepsilon \big( \partial_x^k Z_\varepsilon \big)^2,$$
so that by the H\"older inequality and bounds \eqref{allfolks3} and \eqref{quirecure}, we are led to
\begin{equation}
\label{veolia}
\bigg| \int_0^\tau \int_\R \partial_x^{k+1} \big( Z_\varepsilon^2 \big) \partial_x^k Z_\varepsilon \bigg| \leq K \bigg( \bigg| \int_0^\tau \| \partial_x Z_\varepsilon \|_{L^\infty} \| \partial_x^k Z_\varepsilon \|_{L^2}^2 \bigg| + \varepsilon^2 \bigg| \int_0^\tau \exp K |s| \ \| \partial_x^k Z_\varepsilon(\cdot, s) \|_{L^2} ds \bigg| \bigg).
\end{equation}
We now recall that $\partial_x Z_\varepsilon = \frac{1}{2} \big( \partial_x N_\varepsilon + \partial_x^2 \Theta_\varepsilon \big) - \partial_x \boU$, so that, combining bounds \eqref{praterbis} and \eqref{allfolks2} with the Sobolev embedding theorem,
$$\| \partial_x Z_\varepsilon(\cdot, \tau) \|_{L^\infty(\R)} \leq K,$$
for any $\tau \in \R$. Hence, by \eqref{veolia},
\begin{equation}
\label{virbac}
\bigg| \int_0^\tau \int_\R \partial_x^{k+1} \big( Z_\varepsilon^2 \big) \partial_x^k Z_\varepsilon \bigg| \leq K \bigg( \bigg| \int_0^\tau \int_\R \big( \partial_x^k Z_\varepsilon \big)^2 \bigg| + \varepsilon^2 \bigg| \int_0^\tau \exp K |s| \ \| \partial_x^k Z_\varepsilon(\cdot, s) \|_{L^2(\R)} ds \bigg| \bigg).
\end{equation}
Concerning the last integral on the right-hand side of \eqref{gege}, we invoke \eqref{allfolks1} to obtain
$$\bigg| \varepsilon^2 \int_0^\tau \int_\R \partial_x^k r_\varepsilon \partial_x^k Z_\varepsilon \bigg| \leq K \varepsilon^2 \bigg| \int_0^\tau \exp K |s| \ \| \partial_x^k Z_\varepsilon(\cdot, s) \|_{L^2(\R)} ds \bigg|.$$
Therefore, in view of \eqref{estim-f}, \eqref{foncia} and \eqref{virbac}, we are led to the differential inequality
$$\partial_\tau \boZ_\varepsilon^k(\tau) \leq K \bigg( \sign(\tau) \boZ_\varepsilon^k(\tau) + \varepsilon^4 \exp K |\tau| \bigg),$$
so that by the Gronwall lemma, inequality \eqref{ineq1} also holds for the integer $k$. By induction, this completes the proof of \eqref{ineq1} for the functions $U_\varepsilon^-$ and $\boU^-$.

We next say a few words of the proof for the functions $U_\varepsilon^+$ and $\boU^+$. Setting
$$V_\varepsilon^+(x, \tau) = \frac{1}{2} \big( N_\varepsilon^+(x, \tau) + \partial_x \Theta_\varepsilon^+(x, \tau) \big),$$
the functions $U_\varepsilon^+$ and $V_\varepsilon^+$ satisfy the system of equations
\begin{equation}
\label{slow1+}
- \partial_\tau U_\varepsilon^+ + \partial_x^3 U_\varepsilon^+ + U_\varepsilon^+ \partial_x U_\varepsilon^+ = f_\varepsilon^+ - \varepsilon^2 r_\varepsilon^+,
\end{equation}
where
$$f_\varepsilon^+ \equiv \partial_x \Big( \frac{1}{6} (V_\varepsilon^+)^2 - \partial_x^2 V_\varepsilon^+ + \frac{1}{3} U_\varepsilon^+ V_\varepsilon^+ \Big),$$
and
\begin{equation}
\label{slow2+}
- \partial_\tau V_\varepsilon^+ + \frac{8}{\varepsilon^2} \partial_x V_\varepsilon^+ = g_\varepsilon^+ + \varepsilon^2 r_\varepsilon^+,
\end{equation}
where
$$g_\varepsilon^+ \equiv \partial_x \Big( \partial_x^2 N_\varepsilon^+ + \frac{1}{2} (V_\varepsilon^+)^2 - \frac{1}{6} (U_\varepsilon^+)^2 - \frac{1}{3} U_\varepsilon^+ V_\varepsilon^+ \Big),$$
and the remainder term $r_\varepsilon^+$ is given by the formula
$$r_\varepsilon^+ = \frac{N_\varepsilon^+ \partial_x^3 N_\varepsilon^+}{6 (1 - \frac{\varepsilon^2}{6} N_\varepsilon^+)} + \frac{(\partial_x N_\varepsilon^+) (\partial_x^2 N_\varepsilon^+)}{3 (1 - \frac{\varepsilon^2}{6} N_\varepsilon^+)^2} + \frac{\varepsilon^2}{36} \frac{(\partial_x N_\varepsilon^+)^3}{(1 - \frac{\varepsilon^2}{6} N_\varepsilon^+)^3}.$$
Up to a reverse orientation of time, equations \eqref{slow1+} and \eqref{slow2+} are identical to equations \eqref{slow1} and \eqref{slow2}. In particular, we can apply to the functions $\tau \mapsto U_\varepsilon^+(\cdot, - \tau)$ and $\tau \mapsto V_\varepsilon^+(\cdot, - \tau)$, the analysis developed above to prove \eqref{ineq1} for the functions $U_\varepsilon = U_\varepsilon^-$ and $V_\varepsilon$. Given any $k \in \N$, the associated initial datum
$$\big( U_\varepsilon^+(\cdot, 0), V_\varepsilon^+(\cdot, 0) \big) = \Big( \frac{1}{2} \Big( N_\varepsilon^0 - \partial_x \Theta_\varepsilon^0 \Big), \frac{1}{2} \Big( N_\varepsilon^0 + \partial_x \Theta_\varepsilon^0 \Big) \Big),$$
also satisfies assumption \eqref{grinzing1}, so that there exists some constant $K$ depending only on $K_0$ and $k$, such that
\begin{equation}
\label{ineq1+}
\| U_\varepsilon^+(\cdot, - \tau) - U(\cdot, \tau) \|_{H^k(\R)} \leq K \varepsilon^2 \exp K |\tau|,
\end{equation}
for any $\tau \in \R$. Here, the function $U$ denotes the solution to \eqref{KdV} with initial datum
$$U(\cdot, 0) = U_\varepsilon^+(\cdot, 0).$$
By the uniqueness of the solution to \eqref{KdV-} for any fixed initial datum in $H^k(\R)$, we notice that $U(\cdot, - \tau) = \boU^+(\cdot, \tau)$. Reverting the orientation of time in \eqref{ineq1+}, this completes the proof of \eqref{ineq1} for the functions $U_\varepsilon^+$ and $\boU^+$.

\section{Proof of Theorem \ref{cochon2}}
\label{cinq}

In this section, we provide the proof of Theorem \ref{cochon2}. This first requires to show Proposition \ref{Estim-V}, Lemma \ref{Initier} and Proposition \ref{Estim-f2}.

\subsection{Proof of Proposition \ref{Estim-V}}

In order to estimate the $H^k$-norm of $V_\varepsilon$, we apply the differential operator $\partial_x^k$ to \eqref{slow2}, multiply the resulting equation by $\partial_x^k V_\varepsilon$ and integrate by parts on $\R \times (0, \tau)$. In view of definition \eqref{gG}, this yields
\begin{equation}
\label{sangohan}
\int_\R \big( \partial_x^k V_\varepsilon(\cdot, \tau) \big)^2 = \int_\R \big( \partial_x^k V_\varepsilon^0 \big)^2 + \int_0^\tau \int_\R \partial_x^{k + 1} \big( V_\varepsilon^2 \big) \partial_x^k V_\varepsilon - 2 L_1(\tau) + \frac{1}{3} L_2(\tau) - \frac{2}{3} L_3(\tau) + 2 L_4(\tau),
\end{equation}
where, in view of the fact that $N_\varepsilon = U_\varepsilon + V_\varepsilon$,
\begin{align*}
L_1(\tau) \equiv & \int_0^\tau \int_\R \partial_x^{k + 2} N_\varepsilon \partial_x^{k + 1} V_\varepsilon = \int_0^\tau \int_\R \partial_x^{k + 2} U_\varepsilon \partial_x^{k + 1} V_\varepsilon,\\
L_2(\tau) \equiv & \int_0^\tau \int_\R \partial_x^k \big( U_\varepsilon^2 \big) \partial_x^{k + 1} V_\varepsilon,\\
L_3(\tau) \equiv & \int_0^\tau \int_\R \partial_x^{k + 1} \big( U_\varepsilon V_\varepsilon \big) \partial_x^k V_\varepsilon,
\end{align*}
and
$$L_4(\tau) \equiv \varepsilon^2 \int_0^\tau \int_\R \partial_x^k r_\varepsilon \partial_x^k V_\varepsilon.$$
We now estimate each integral $L_j(\tau)$ as in the proof of Theorem \ref{cochondore}. For the first one, we have

\setcounter{step}{0}
\begin{step}
\label{Estim-L1}
Under the assumptions of Proposition \ref{Estim-V}, there exists a positive constant $K$ depending only on $K_0$ and $k$, such that
\begin{equation}
\label{estim-L1}
\begin{split}
\bigg| L_1(\tau) + \frac{\varepsilon^2}{48} \int_0^\tau \int_\R \partial_x^{k + 2} U_\varepsilon \partial_x^{k + 1} \big( U_\varepsilon^2 \big) \bigg| \leq & K \varepsilon^2 \bigg( \Big( \varepsilon^2 + \| \partial_x^k V_\varepsilon^0 \|_{L^2(\R)} + \| \partial_x^k V_\varepsilon(\cdot, \tau) \|_{L^2(\R)} \Big) \exp K |\tau|\\
& + \bigg| \int_0^\tau \exp K |s| \ \| V_\varepsilon(\cdot, s) \|_{H^k(\R)} ds \bigg| \bigg).
\end{split}
\end{equation}
\end{step}

In view of \eqref{slow2}, we compute
\begin{align*}
\int_0^\tau \int_\R \partial_x^{k + 2} U_\varepsilon \partial_x^{k + 1} V_\varepsilon = & - \frac{\varepsilon^2}{8} \bigg[ \int_\R \partial_x^{k + 2} U_\varepsilon \partial_x^k V_\varepsilon \bigg]_0^\tau + \frac{\varepsilon^2}{8} \int_0^\tau \int_\R \partial_\tau \partial_x^{k + 2} U_\varepsilon \partial_x^k V_\varepsilon\\
+ & \frac{\varepsilon^2}{8} \int_0^\tau \int_\R \partial_x^{k + 2} U_\varepsilon \big( \partial_x^k g_\varepsilon + \varepsilon^2 \partial_x^k r_\varepsilon),
\end{align*}
so that by \eqref{frites} and \eqref{allfolks1},
\begin{equation}
\label{garlic}
\begin{split}
\bigg| L_1(\tau) - \frac{\varepsilon^2}{8} \int_0^\tau \int_\R \partial_x^{k + 2} U_\varepsilon \partial_x^k g_\varepsilon \bigg| \leq & K \varepsilon^2 \bigg( \Big( \varepsilon^2 + \| \partial_x^k V_\varepsilon^0 \|_{L^2(\R)} + \| \partial_x^k V_\varepsilon(\cdot, \tau) \|_{L^2(\R)} \Big) \exp K |\tau|\\
& + \bigg| \int_0^\tau \exp K |s| \ \| \partial_x^k V_\varepsilon(\cdot, s) \|_{L^2(\R)} ds \bigg| \bigg).
\end{split}
\end{equation}
In view of \eqref{gG}, we next write
$$\int_0^\tau \int_\R \partial_x^{k + 2} U_\varepsilon \partial_x^k g_\varepsilon = \int_0^\tau \int_\R \partial_x^{k + 2} U_\varepsilon \Big( \partial_x^{k + 3} V_\varepsilon + \frac{1}{2} \partial_x^{k + 1} \big( V_\varepsilon^2 \big) - \frac{1}{6} \partial_x^{k + 1} \big( U_\varepsilon^2 \big) - \frac{1}{3} \partial_x^{k + 1} \big( U_\varepsilon V_\varepsilon \big) \Big),$$
so that integrating by parts, we are led to
\begin{align*}
& \int_0^\tau \int_\R \partial_x^{k + 2} U_\varepsilon \partial_x^k g_\varepsilon + \frac{1}{6} \int_0^\tau \int_\R \partial_x^{k + 2} U_\varepsilon \partial_x^{k + 1} \big( U_\varepsilon^2 \big)\\
= - \int_0^\tau \int_\R & \partial_x^{k + 5} U_\varepsilon \partial_x^k V_\varepsilon - \int_0^\tau \int_\R \partial_x^{k + 3} U_\varepsilon \Big( \frac{1}{2} \partial_x^k \big( V_\varepsilon^2 \big) - \frac{1}{3} \partial_x^k \big( U_\varepsilon V_\varepsilon \big) \Big).
\end{align*}
Combining \eqref{frites} with the Leibniz formula, the H\"older inequality and the Sobolev embedding theorem, we deduce that
$$\bigg| \int_0^\tau \int_\R \partial_x^{k + 2} U_\varepsilon \partial_x^k g_\varepsilon + \frac{1}{6} \int_0^\tau \int_\R \partial_x^{k + 2} U_\varepsilon \partial_x^{k + 1} \big( U_\varepsilon^2 \big) \bigg| \leq K \bigg| \int_0^\tau \exp K |s| \ \| V_\varepsilon(\cdot, s) \|_{H^k(\R)} ds \bigg|$$
Invoking \eqref{garlic}, this completes the proof of \eqref{estim-L1}.

We similarly derive for the fourth integral on the right-hand side of \eqref{sangohan}.

\begin{step}
\label{Estim-L2}
Under the assumptions of Proposition \ref{Estim-V}, there exists a positive constant $K$ depending only on $K_0$ and $k$, such that
\begin{align*}
\bigg| L_2(\tau) - \frac{\varepsilon^2}{8} \int_0^\tau \int_\R \partial_x^k \big( U_\varepsilon^2 \big) \partial_x^{k + 3} U_\varepsilon \bigg| \leq & K \varepsilon^2 \bigg( \Big( \varepsilon^2 + \| \partial_x^k V_\varepsilon^0 \|_{L^2(\R)} + \| \partial_x^k V_\varepsilon(\cdot, \tau) \|_{L^2(\R)} \Big) \exp K |\tau|\\
& + \bigg| \int_0^\tau \exp K |s| \ \| V_\varepsilon(\cdot, s) \|_{H^k(\R)} ds \bigg| \bigg).
\end{align*}
\end{step}

The proof is identical to the proof of Step \ref{Estim-L1}, so that we omit it, and instead we turn to the fifth integral on the right-hand side of \eqref{sangohan}.

\begin{step}
\label{Estim-L3}
Under the assumptions of Proposition \ref{Estim-V}, there exists a positive constant $K$ depending only on $K_0$ and $k$, such that
\begin{equation}
\label{estim-L3}
\big| L_3(\tau) \big| \leq K \varepsilon^2 \bigg( \Big( \varepsilon^2 + \| \partial_x^k V_\varepsilon^0 \|_{L^2} + \| \partial_x^k V_\varepsilon(\cdot, \tau) \|_{L^2} \Big) \exp K |\tau| + \bigg| \int_0^\tau \exp K |s| \ \| \partial_x^k V_\varepsilon(\cdot, s) \|_{L^2} ds \bigg| \bigg).
\end{equation}
\end{step}

Using the Leibniz formula, we are led to
\begin{equation}
\label{tortuegeniale}
L_3(\tau) = \sum_{j = 0}^{k + 1} \binom{k + 1}{j} \boL_3^j(\tau),
\end{equation}
where
$$\boL_3^j(\tau) \equiv \int_0^\tau \int_\R \partial_x^{k + 1 - j} U_\varepsilon \partial_x^j V_\varepsilon \partial_x^k V_\varepsilon.$$
In view of \eqref{planus}, given any $0 \leq j \leq k$, we have integrating by parts in space, then in time,
\begin{align*}
\boL_3^j(\tau) = & \frac{\varepsilon^2}{8} \bigg[ \int_\R \partial_x^{k - j} U_\varepsilon \partial_x^j V_\varepsilon \partial_x^k V_\varepsilon \bigg]_0^\tau - \frac{\varepsilon^2}{8} \int_0^\tau \int_\R \partial_\tau \partial_x^{k - j} U_\varepsilon \partial_x^j V_\varepsilon \partial_x^k V_\varepsilon\\
- & \frac{\varepsilon^2}{8} \int_0^\tau \int_\R \partial_x^{k - j} U_\varepsilon \big( \partial_x^j V_\varepsilon \partial_x^k g_\varepsilon + \partial_x^k V_\varepsilon \partial_x^j g _\varepsilon \big) - \frac{\varepsilon^4}{8} \int_0^\tau \int_\R \partial_x^{k - j} U_\varepsilon \big( \partial_x^j V_\varepsilon \partial_x^k r_\varepsilon + \partial_x^k V_\varepsilon \partial_x^j r_\varepsilon \big),
\end{align*}
so that by \eqref{frites} and \eqref{allfolks1},
\begin{equation}
\label{vegeta}
\big| \boL_3^j(\tau) \big| \leq K \varepsilon^2 \bigg( \Big( \varepsilon^2 + \| \partial_x^k V_\varepsilon^0 \|_{L^2} + \| \partial_x^k V_\varepsilon(\cdot, \tau) \|_{L^2} \Big) \exp K |\tau| + \bigg| \int_0^\tau \exp K |s| \ \| \partial_x^k V_\varepsilon(\cdot, s) \|_{L^2} ds \bigg| \bigg).
\end{equation}
For $j = k + 1$, we can also invoke \eqref{planus} to establish that
\begin{align*}
\boL_3^{k + 1}(\tau) = \frac{1}{2} \int_0^\tau \int_\R U_\varepsilon \partial_x \big( (\partial_x^k V_\varepsilon)^2 \big) = & - \frac{\varepsilon^2}{16} \bigg[ \int_\R U_\varepsilon \big( \partial_x^k V_\varepsilon \big)^2 \bigg]_0^\tau + \frac{\varepsilon^2}{16} \int_0^\tau \int_\R \partial_\tau U_\varepsilon \big( \partial_x^k V_\varepsilon \big)^2\\
& + \frac{\varepsilon^2}{8} \int_0^\tau \int_\R U_\varepsilon \partial_x^k V_\varepsilon \big( \partial_x^k g_\varepsilon + \varepsilon^2 \partial_x^k r_\varepsilon \big),
\end{align*}
so that \eqref{vegeta} follows similarly. In view of \eqref{tortuegeniale}, this completes the proof of \eqref{estim-L3}.

Using \eqref{allfolks1}, we now compute directly the next estimate of the integral $L_4(\tau)$.

\begin{step}
\label{Estim-L4}
Under the assumptions of Proposition \ref{Estim-V}, there exists a positive constant $K$ depending only on $K_0$ and $k$, such that
$$\big| L_4(\tau) \big| \leq K \varepsilon^2 \bigg| \int_0^\tau \exp K |s| \ \| \partial_x^k V_\varepsilon(\cdot, s) \|_{H^k(\R)} ds \bigg|.$$
\end{step}

We finally complete the proof of \eqref{estim-V} by induction.

\begin{proof}[End of the proof of Proposition \ref{Estim-V}]
Let us denote
$$S_\varepsilon^k(\tau) = \sup_{s \in (0, \tau)} \| \partial_x^k V_\varepsilon(\cdot, s) \|_{L^2(\R)}.$$
For $k = 0$, we have in view of \eqref{sangohan}, and Steps \ref{Estim-L1}, \ref{Estim-L2}, \ref{Estim-L3} and \ref{Estim-L4},
\begin{align*}
\int_\R V_\varepsilon(x, \tau)^2 dx \leq \| V_\varepsilon^0 \|_{L^2(\R)}^2 + K \varepsilon^2 \bigg( & \Big( \varepsilon^2 + \| V_\varepsilon^0 \|_{L^2(\R)} + \| V_\varepsilon(\cdot, \tau) \|_{L^2(\R)} \Big) \exp K |\tau|\\
& + \bigg| \int_0^\tau \exp K |s| \ \| V_\varepsilon(\cdot, s) \|_{L^2(\R)} ds \bigg| \bigg),
\end{align*}
so that
$$S_\varepsilon^0(\tau)^2 \leq \Big( \| V_\varepsilon^0 \|_{L^2(\R)}^2 + K \varepsilon^4 \Big) \exp K |\tau| \ + K \varepsilon^2 \exp K |\tau| \ S_\varepsilon^0(\tau).$$
Using the inequality $2 a b \leq a^2 + b^2$, this completes the proof of \eqref{estim-V} for $k = 0$.

We now assume that \eqref{estim-V} holds for any $0 \leq k \leq \kappa - 1$ and establish it for $k = \kappa \geq 1$. Invoking the Leibniz formula, we compute
$$\int_0^\tau \int_\R \partial_x^{\kappa + 1} \big( V_\varepsilon^2 \big) \partial_x^\kappa V_\varepsilon = \sum_{j = 2}^{\kappa - 1} \binom{\kappa + 1}{j} \int_0^\tau \int_\R \partial_x^{\kappa + 1 - j} V_\varepsilon \partial_x^j V_\varepsilon \partial_x^\kappa V_\varepsilon - (4 k + 2) \int_0^\tau \int_\R V_\varepsilon \partial_x^\kappa V_\varepsilon \partial_x^{\kappa + 1} V_\varepsilon,$$
so that by the inductive assumption and bound \eqref{frites},
\begin{align*}
\bigg| \int_0^\tau \int_\R \partial_x^{\kappa + 1} \big( V_\varepsilon^2 \big) \partial_x^\kappa V_\varepsilon \bigg| \leq & K \Big( \| V_\varepsilon^0 \|_{H^{\kappa - 1}(\R)} + \varepsilon^2 \Big) \bigg( \Big( \| V_\varepsilon^0 \|_{H^{\kappa - 1}(\R)} + \varepsilon^2 \Big) \exp K |\tau|\\
& + \bigg| \int_0^\tau \exp K |s| \ \| \partial_x^\kappa V_\varepsilon(\cdot, s) \|_{L^2(\R)} ds \bigg| \bigg).
\end{align*}
Combining with \eqref{sangohan}, the inductive assumption and Steps \ref{Estim-L1}, \ref{Estim-L2}, \ref{Estim-L3} and \ref{Estim-L4}, we are led to
$$S_\varepsilon^\kappa(\tau)^2 \leq \| \partial_x^\kappa V_\varepsilon^0 \|_{L^2(\R)}^2 + K \Big( \varepsilon^2 + \| V_\varepsilon^0 \|_{H^{\kappa - 1}(\R)} \Big) \Big( \varepsilon^2 + \| V_\varepsilon^0 \|_{H^\kappa(\R)} + S_\varepsilon^\kappa(\tau) \Big) \exp K |\tau|.$$
Using again the inequality $2 a b \leq a^2 + b^2$, this completes the proof of \eqref{estim-V} for $k = \kappa$. By induction, this concludes the proof of Proposition \ref{Estim-V}.
\end{proof}

\subsection{Proof of Lemma \ref{Initier}}

Since $F$ and $G$ are solutions to \eqref{KdV}, their difference $H \equiv F - G$ is solution to
\begin{equation}
\label{eq:hache}
\partial_\tau H + \partial_x^3 H + F \partial_x H + H \partial_x G = 0.
\end{equation}
In order to prove \eqref{init}, we now compute inductively energy estimates on \eqref{eq:hache}.

For $k = 0$, we multiply \eqref{eq:hache} by $H$ and integrate by parts on $\R$ to obtain
\begin{equation}
\label{hachisparmentier}
\partial_\tau \bigg( \int_\R H^2 \bigg) = \int_\R \Big( \partial_x F - 2 \partial_x G \Big) H^2.
\end{equation}
Since $F^0$ and $G^0$ are in $H^2(\R)$, we recall (see also the proof of \eqref{allfolks2}) that in view of the integrability properties of \eqref{KdV}, there exists a constant $K$ depending only on the $H^2$-norms of $F^0$ and $G^0$, such that
$$\| F(\cdot, \tau) \|_{H^2(\R)} + \| G(\cdot, \tau) \|_{H^2(\R)} \leq K,$$
for any $\tau \in \R$. Applying the H\"older inequality and the Sobolev embedding theorem to \eqref{hachisparmentier}, we are led to
$$\bigg| \partial_\tau \bigg( \int_\R H^2 \bigg) \bigg| \leq K \int_\R H^2,$$
so that \eqref{init} follows from the Gronwall lemma.

We now assume that \eqref{init} holds for any $0 \leq k \leq \kappa - 1$ and derive it for $k = \kappa \geq 1$. For this purpose, we apply the differential operator $\partial_x^\kappa$ to \eqref{eq:hache}, multiply the resulting equation by $\partial_x^\kappa H$ and integrate by parts on $\R$. This provides
\begin{equation}
\label{steackhache}
\frac{1}{2} \partial_\tau \bigg( \int_\R \big( \partial_x^\kappa H \big)^2 \bigg) = - \int_\R \partial_x^\kappa \big( F \partial_x H \big) \partial_x^\kappa H - \int_\R \partial_x^\kappa \big( H \partial_x G \big) \partial_x^\kappa H.
\end{equation}
By the Leibniz formula, the first integral on the right-hand side of \eqref{steackhache} reduces to
$$\int_\R \partial_x^\kappa \big( F \partial_x H \big) \partial_x^\kappa H = \sum_{j = 1}^\kappa \binom{\kappa}{j} \int_\R \partial_x^j F \partial_x^{\kappa + 1 -j} H \partial_x^\kappa H - \frac{1}{2} \int_\R \partial_x F \big( \partial_x^\kappa H \big)^2.$$
In view of the integrability properties of \eqref{KdV}, there again exists some constant $K$ depending only on the $H^{\kappa + 2}$-norms of $F^0$ and $G^0$, such that
$$\| F(\cdot, \tau) \|_{H^{\kappa + 2}(\R)} + \| G(\cdot, \tau) \|_{H^{\kappa + 2}(\R)} \leq K,$$
for any $\tau \in \R$. Hence, we are led to
$$\bigg| \int_\R \partial_x^\kappa \big( F \partial_x H \big) \partial_x^\kappa H \bigg| \leq K \| H \|_{H^\kappa(\R)}^2.$$
The same estimate holds for the second integral on the right-hand side of \eqref{steackhache}. Invoking the inductive assumption, we deduce that
$$\bigg| \partial_\tau \bigg( \int_\R \big( \partial_x^\kappa H(x, \tau) \big)^2 dx \bigg) \bigg| \leq K \bigg( \int_\R \big( \partial_x^\kappa H(x, \tau) \big)^2 dx + \| F^0 - G^0 \|_{H^{\kappa - 1}(\R)}^2 \exp K |\tau| \bigg).$$
Inequality \eqref{init} follows applying the Gronwall lemma. By induction, this concludes the proof of Lemma \ref{Initier}.

\subsection{Proof of Proposition \ref{Estim-f2}}

The proof is a direct adaptation of the proof of Proposition \ref{Estim-f} using only assumption \eqref{grinzing2}. In view of \eqref{mandanda}, and Lemma \ref{Itisbounded} which remains valid under assumption \eqref{grinzing2}, we have
\begin{equation}
\label{trunk}
\begin{split}
\bigg| \int_0^\tau \int_\R \partial_x^k & f_\varepsilon \partial_x^k Z_\varepsilon \bigg| \leq K\bigg( \bigg| \int_0^\tau \int_\R \partial_x^{k + 1} U_\varepsilon V_\varepsilon \partial_x^k Z_\varepsilon \bigg| + \varepsilon^2 \sum_{j = 1}^k \bigg| \int_0^\tau \int_\R \partial_x^k f_\varepsilon \partial_x^{k - j} V_\varepsilon \partial_x^j U_\varepsilon \bigg|\\
+ & \varepsilon^2 \Big( \varepsilon^2 + \| \partial_x^k Z_\varepsilon(\cdot, \tau) \|_{L^2(\R)} \Big) \exp K |\tau| + \varepsilon^2 \bigg| \int_0^\tau \exp K |s| \ \| Z_\varepsilon(\cdot, s) \|_{H^k(\R)} ds \bigg| \bigg),
\end{split}
\end{equation}
where $K$ refers to some positive constant depending only on $K_0$ and $k$. We then invoke Proposition \ref{Estim-V} to bound the first and the second integrals on the right-hand side of \eqref{trunk}. Combining with \eqref{frites} and \eqref{allfolks1}, this leads to
\begin{align*}
& \bigg| \int_0^\tau \int_\R \partial_x^{k + 1} U_\varepsilon V_\varepsilon \partial_x^k Z_\varepsilon \bigg| + \varepsilon^2 \sum_{j = 1}^k \bigg| \int_0^\tau \int_\R \partial_x^k f_\varepsilon \partial_x^{k - j} V_\varepsilon \partial_x^j U_\varepsilon \bigg|\\
\leq K \Big( \| V_\varepsilon^0 & \|_{H^k(\R)} + \varepsilon^2 \Big) \bigg( \varepsilon^2 \exp K |\tau| + \bigg| \int_0^\tau \exp K |s| \ \| \partial_x^k Z_\varepsilon(\cdot, s) \|_{L^2(\R)} ds \bigg| \bigg).
\end{align*}
In view of \eqref{trunk}, this completes the proof of \eqref{estim-f2}, and of Proposition \ref{Estim-f2}.

\subsection{Proof of Theorem \ref{cochon2} completed}
\label{Uhu}

As mentioned in the introduction, we first focus on the coordinate $x^-$ and the associated functions $N_\varepsilon \equiv N_\varepsilon^-$ and $\Theta_\varepsilon \equiv \Theta_\varepsilon^-$. Theorem \ref{cochon2} then follows from combining decompositions \eqref{trianguler1} and \eqref{trianguler2} with estimates \eqref{estim-V}, \eqref{initfin} and \eqref{ineq12}, once \eqref{ineq12} is established, which we do next.

\begin{proof}[Proof of inequality \eqref{ineq12}]
The proof is an adaptation of the proof of Theorem \ref{cochondore}. For $k = 0$, coming back to \eqref{krilin}, we deduce from \eqref{estim-f2}, \eqref{allfolks1} and \eqref{allfolks2} that
\begin{align*}
\partial_\tau \boZ_\varepsilon^0(\tau) \leq & K \bigg( |\boZ_\varepsilon^0(\tau)| + \big( \varepsilon^2 + \| V_\varepsilon^0 \|_{L^2(\R)} \big) \bigg| \int_0^\tau \exp K |s| \ \| Z_\varepsilon(\cdot, s) \|_{L^2(\R)} ds \bigg|\\
+ & \varepsilon^2 \Big( \varepsilon^2 + \| V_\varepsilon^0 \|_{L^2(\R)} + \| Z_\varepsilon(\cdot, \tau) \|_{L^2(\R)} \Big) \exp K |\tau| \bigg),
\end{align*}
where $K$ is some positive constant depending only on $K_0$. Using again the inequality $2 |a b| \leq a^2 + b^2$ and identity \eqref{derivz0}, we are led to the differential inequality
$$\partial_\tau \boZ_\varepsilon^0(\tau) \leq K \Big( \sign(\tau) \boZ_\varepsilon^0(\tau) + \big( \varepsilon^2 + \| V_\varepsilon^0 \|_{L^2(\R)} \big)^2 \exp K |\tau| \Big),$$
so that
$$|\boZ_\varepsilon^0(\tau)| \leq K' \big( \varepsilon^2 + \| V_\varepsilon^0 \|_{L^2(\R)} \big)^2 \exp K' |\tau|,$$
where $K'$ is some further positive constant depending only on $K_0$. Combining with \eqref{derivz0} and \eqref{gronwall0}, this provides \eqref{ineq12} for $k = 0$.

We now assume that \eqref{ineq12} holds for any $0 \leq j \leq k - 1$, i.e. that there exists a positive constant depending only on $K_0$ and $k$, such that
$$\| Z_\varepsilon(\cdot, \tau) \|_{H^{k - 1}(\R)} \leq K \big( \varepsilon^2 + \| V_\varepsilon^0 \|_{H^{k - 1}(\R)} \big) \exp K |\tau|,$$
for any $\tau \in \R$. We then bound any integral on the right-hand side of \eqref{gege}. For the first and second ones, we have following the lines of the proofs of \eqref{foncia} and \eqref{virbac},
\begin{equation}
\label{virbac2}
\begin{split}
& \bigg| \int_0^\tau \int_\R \partial_x^{k+1} \big( \boU Z_\varepsilon \big) \partial_x^k Z_\varepsilon \bigg| + \bigg| \int_0^\tau \int_\R \partial_x^{k+1} \big( Z_\varepsilon^2 \big) \partial_x^k Z_\varepsilon \bigg|\\
\leq K \bigg( \bigg| \int_0^\tau \int_\R \big( \partial_x^k & Z_\varepsilon \big)^2 \bigg| + \big( \varepsilon^2 + \| V_\varepsilon^0 \|_{H^{k - 1}(\R)} \big) \bigg| \int_0^\tau \exp K |s| \ \| \partial_x^k Z_\varepsilon(\cdot, s) \|_{L^2(\R)} ds \bigg| \bigg).
\end{split}
\end{equation}
Concerning the last integral on the right-hand side of \eqref{gege}, we invoke \eqref{allfolks1} to obtain
$$\bigg| \varepsilon^2 \int_0^\tau \int_\R \partial_x^k r_\varepsilon \partial_x^k Z_\varepsilon \bigg| \leq K \varepsilon^2 \bigg| \int_0^\tau \exp K |s| \ \| \partial_x^k Z_\varepsilon(\cdot, s) \|_{L^2(\R)} ds \bigg|.$$
Therefore, in view of \eqref{estim-f2} and \eqref{virbac2}, we are led to the differential inequality
$$\partial_\tau \boZ_\varepsilon^k(\tau) \leq K \bigg( \sign(\tau) \boZ_\varepsilon^k(\tau) + \big( \varepsilon^2 + \| V_\varepsilon^0 \|_{H^k(\R)} \big)^2 \exp K |\tau| \bigg),$$
so that by the Gronwall lemma, inequality \eqref{ineq12} also holds for the integer $k$. By induction, this completes the proof of \eqref{ineq12}.
\end{proof}

We are now in position to end the proof of Theorem \ref{cochon2}.

\begin{proof}[End of the proof of Theorem \ref{cochon2}]
As mentioned above, Theorem \ref{cochon2} is a direct consequence of decompositions \eqref{trianguler1} and \eqref{trianguler2} and estimates \eqref{estim-V}, \eqref{initfin} and \eqref{ineq12}, when the coordinate $x^-$ is considered. For the functions $N_\varepsilon^+$ and $\partial_x \Theta_\varepsilon^+$, the proof reduces as in Theorem \ref{cochondore} to consider the system of equations \eqref{slow1+}-\eqref{slow2+} instead of the system \eqref{slow1}-\eqref{slow2}. Since the functions $\tau \mapsto U_\varepsilon^+(\cdot, - \tau)$ and $\tau \mapsto V_\varepsilon^+(\cdot, - \tau)$ are solutions to \eqref{slow1} and \eqref{slow2}, we can apply Propositions \ref{Estim-V} and \ref{Estim-f2} to them in order to obtain inequalities \eqref{estim-V} and \eqref{ineq12} in the coordinate $x^+$. Combining with the versions of \eqref{trianguler1} and \eqref{trianguler2} in the coordinate $x^+$, and Lemma \ref{Initier}, this provides \eqref{ineq2} in the coordinate $x^+$, and concludes the proof of Theorem \ref{cochon2}.
\end{proof}

\appendix
\section{Defining a notion of the mass for \eqref{GP}}
\label{Alamasse}

The purpose of this appendix is to provide a framework where the notion of mass for the one-dimensional Gross-Pitaevskii equation may be rigorously handled. At least on a formal level, the mass
\footnote{It would be more appropriate to call it a relative mass (with respect to the vacuum) since it may be negative as such.}
may be defined by
\begin{equation}
\label{alamasse}
m(\Psi) = \frac{1}{2} \int_\R \Big( 1 - |\Psi|^2 \Big),
\end{equation}
and it is a conserved quantity along the Gross-Pitaevskii flow. Indeed, a solution $\Psi$ to \eqref{GP} satisfies the conservation law
\begin{equation}
\label{eqmass}
\partial_t \eta = 2 \partial_x \big( \langle i \Psi, \partial_x \Psi \rangle \big),
\end{equation}
where we denote as above $\eta \equiv 1 - |\Psi|^2$. Hence, we have
$$\partial_t m(\Psi) = \frac{1}{2} \int_\R \partial_t \eta = \int_\R \partial_x \big( \langle i \Psi, \partial_x \Psi \rangle \big) = 0,$$
provided that the functions $\Psi$ and $\eta$ are sufficiently smooth and decay suitably at infinity.

The quantity $m(\psi)$ is however not well-defined in general for an arbitrary function $\psi$ in the energy space $X^1(\R)$. Consider for instance, the function $\psi$ defined by
$$\psi(x) = \frac{|x|}{|x| + 1}, \ \forall x \in \R,$$
which belongs to $X^1(\R)$, but for which $m(\psi) = + \infty$. In order to circumvent this difficulty, we introduce the set
$$X_\boM(\R) \equiv \{ \psi \in X^1(\R), \ {\rm s.t.} \ 1 - |\psi|^2 \in \boM(\R)\}.$$
We first claim that the Gross-Pitaevskii equation is well-posed in this new functional setting.

\begin{lemma}
\label{Wellposed}
Given any function $\Psi_0 \in X_\boM(\R)$, there exists a unique solution $\Psi(\cdot, t)$ to \eqref{GP} in $\boC^0(\R, X_\boM(\R))$ with initial datum $\Psi_0$. Moreover, there exists a universal constant $K$ such that
\begin{equation}
\label{metalips}
\| \eta(t) - \eta(s) \|_{\boM(\R)} \leq K \Big( E(\Psi_0)^\frac{1}{2} \big( 1 + E(\Psi_0)^\frac{1}{2} \big) |t - s| + \| \eta(t) - \eta(s) \|_{L^2(\R)} \Big),
\end{equation}
for any $(s, t) \in \R^2$.
\end{lemma}

\begin{remark}
In view of Proposition \ref{thm:existe} and Lemma \ref{Wellposed}, the Gross-Pitaevskii equation is also globally well-posed in the space
$$X_\boM^k(\R) \equiv \{ \psi \in X^k(\R), \ {\rm s.t.} \ 1 - |\psi|^2 \in \boM(\R)\},$$
for any $k \geq 2$.
\end{remark}

\begin{proof}
We recall that in view of Proposition \ref{thm:existe}, the Gross-Pitaevskii equation is well-posed in $X^1(\R)$ with conservation of the energy $E$, i.e.
\begin{equation}
\label{consE}
E \big( \Psi(\cdot, t) \big) = E \big( \Psi_0 \big),
\end{equation}
for any $t \in \R$. Therefore, the proof of Lemma \ref{Wellposed} reduces to show that the function $\eta(\cdot, t) \equiv 1 - |\Psi(\cdot, t)|^2$ associated to the unique solution $\Psi(\cdot, t)$ in the space $X^1(\R)$ is continuous with values in $\boM(\R)$. This fact is a direct consequence of \eqref{metalips} which we show next.

For the proof of \eqref{metalips}, we introduce a cut-off function $\chi \in \boC^\infty(\R)$ such that $0 \leq \chi \leq 1$,
$$\chi(x) = 1, \ {\rm for} \ x \leq 0, \ {\rm and} \ \chi(x) = 0, \ {\rm for} \ x \geq 1,$$
and denote
$$\chi_{a, b}(x) = \left\{ \begin{array}{lll} \chi(a - x), \ {\rm for} \ x \leq a,\\ 1, \ {\rm for} \ a \leq x \leq b,\\ \chi(x - b), \ {\rm for} \ x \geq b, \end{array} \right.$$
for any given numbers $a < b$. When $\Psi$ is a solution to \eqref{GP} in $X^1(\R)$, identity \eqref{eqmass} holds in the sense of distributions and involves quantities which are in $H^{- 1}(\R)$, so that we may multiply \eqref{eqmass} by the test function $\chi_{a, b}$ and integrate by parts to obtain
$$\partial_t \bigg( \int_\R \eta \chi_{a, b} \bigg) = \int_\R \partial_t \eta \chi_{a, b} = 2 \int_\R \partial_x \big( \langle i \Psi, \partial_x \Psi \rangle \big) \chi_{a, b} = - 2 \int_\R \langle i \Psi, \partial_x \Psi \rangle \partial_x \chi_{a, b}.$$
By the Cauchy-Schwarz inequality, we are led to
\begin{equation}
\label{aubry}
\bigg| \partial_t \bigg( \int_\R \eta \chi_{a, b} \bigg) \bigg| \leq 2 \big\| |\Psi| \partial_x \chi_{a, b} \big\|_{L^2(I(a, b))} \big\| \partial_x \Psi \big\|_{L^2(I(a, b))},
\end{equation}
where we denote $I(a, b) = (a - 1, a) \cup (b, b + 1)$. We now recall that it is proved in \cite{Gerard1} that there exists some universal constant $K$ such that
\begin{equation}
\label{psiinfini}
\| \psi \|_{L^\infty} \leq K \big( 1 + E(\psi)^\frac{1}{2} \big),
\end{equation}
for any $\psi \in L^2(\R)$, so that \eqref{aubry} may be recast as
$$\bigg| \partial_t \bigg( \int_\R \eta \chi_{a, b} \bigg) \bigg| \leq K \Big( 1 + E(\Psi)^\frac{1}{2} \Big) E(\Psi)^\frac{1}{2},$$
where $K$ denotes a further positive constant, depending only on our choice of the function $\chi$. Integrating in time and invoking the conservation of the energy provided by \eqref{consE}, we are led to
\begin{equation}
\label{hollande}
\bigg| \int_\R \big( \eta(\cdot, t) - \eta(\cdot, s) \big) \chi_{a, b} \bigg| \leq K \big( 1 + E(\Psi_0)^\frac{1}{2} \Big) E(\Psi_0)^\frac{1}{2} |t - s|.
\end{equation}
Notice finally that
\begin{equation}
\label{dray}
\bigg| \int_\R f \chi_{a, b} - \int_a^b f \bigg| = \bigg| \int_{I(a, b)} f \chi_{a, b} \bigg| \leq K \| f \|_{L^2(I(a, b))},
\end{equation}
for any function $f \in X^1(\R)$, so that \eqref{metalips} is a consequence of \eqref{hollande} and \eqref{dray} (for $f = \eta(t) - \eta(s)$). This completes the proof of Lemma \ref{Wellposed}.
\end{proof}

We now turn to the notion of mass and define
\begin{equation}
\label{weakmass1}
m^+(\psi) = \frac{1}{2} \bigg( \underset{x \to + \infty}{\limsup} \int_0^x \big( 1 - |\psi|^2 \big) + \underset{y \to - \infty}{\limsup} \int_y^0 \big( 1 - |\psi|^2 \big) \bigg),
\end{equation}
and
\begin{equation}
\label{weakmass2}
m^-(\psi) = \frac{1}{2} \bigg( \underset{x \to + \infty}{\liminf} \int_0^x \big( 1 - |\psi|^2 \big) + \underset{y \to - \infty}{\liminf} \int_y^0 \big( 1 - |\psi|^2 \big) \bigg),
\end{equation}
for any function $\psi \in X_\boM(\R)$. Recall that the above integrals are bounded and continuous functions of $x$ and $y$, when $\psi$ belongs to $X_\boM(\R)$, so that $m^+(\psi)$ and $m^-(\psi)$ are well-defined
\footnote{In our definitions of $m^+(\psi)$ and $m^-(\psi)$, the number $0$ may be replaced by any arbitrary other real number.}.

We next show that both the quantities $m^+$ and $m^-$ are conserved along the Gross-Pitaevskii flow provided that the initial datum $\Psi_0$ belongs to $X_\boM(\R)$.

\begin{lemma}
\label{Conspm}
Given any function $\Psi_0 \in X_\boM(\R)$, we have
\begin{equation}
\label{plusoumoins}
m^+(\Psi(\cdot, t)) = m^+(\Psi_0), \ {\rm and} \ m^-(\Psi(\cdot, t)) = m^-(\Psi_0),
\end{equation}
for any $t \in \R$.
\end{lemma}

\begin{proof}
Given any numbers $a < b$, we deduce from \eqref{consE}, \eqref{aubry} and \eqref{psiinfini} that
$$\bigg| \partial_t \bigg( \int_\R \eta \chi_{a, b} \bigg) \bigg| \leq K \big( 1 + E(\Psi_0)^\frac{1}{2} \Big) \| \partial_x \Psi \|_{L^2(I(a, b))},$$
for any fixed $t \in \R$. Integrating this relation in time, combining with \eqref{dray} and applying the Cauchy-Schwarz inequality in time, we are led to
\begin{equation}
\label{peillon}
\begin{split}
& \bigg| \int_a^b \eta(x, t) dx - \int_a^b \eta(x, 0) dx \bigg|
\\ \leq K \bigg( \Big( 1 + E(\Psi_0)^\frac{1}{2} \Big)\Big( \int_0^t \| & \partial_x \Psi(\cdot, s) \|_{L^2(I(a, b))} ds\Big) + \| \eta(\cdot, t) - \eta(\cdot, 0) \|_{L^2(I(a, b))} \bigg)
\\ \leq K \bigg( \Big( 1 + E(\Psi_0)^\frac{1}{2} \Big) |t|^\frac{1}{2} \bigg| \int_0^t & \int_{I(a, b)} \frac{1}{2} |\partial_x \Psi(x, s)|^2 dx ds \bigg|^\frac{1}{2} + \| \eta(\cdot, 0) \|_{L^2(I(a, b))} + \|\eta(\cdot, t) \|_{L^2(I(a, b))} \bigg).
\end{split}
\end{equation}
On the other hand, it follows from the conservation of the energy that
$$\bigg| \int_0^t \int_\R \Big( \frac{1}{2} |\partial_x \Psi(x, s)|^2 + \frac{1}{4} |\eta(x, s)|^2 \Big) dx ds \bigg| = |t| \ E(\Psi_0) < + \infty,$$
so that, by the dominated convergence theorem,
$$\int_0^t \int_{I(a, b)} \Big( \frac{1}{2} |\partial_x \Psi(x, s)|^2 + \frac{1}{4} |\eta(x, s)|^2 \Big) dx ds \to 0, \ {\rm as} \ a \to - \infty \ {\rm and} \ b \to + \infty.$$
Similarly, since $\eta(\cdot, 0)$ and $\eta(\cdot, t)$ belong to $L^2(\R)$,
$$\| \eta(\cdot, 0) \|_{L^2(I(a, b))} + \| \eta(\cdot, t) \|_{L^2(I(a, b))} \to 0, \ {\rm as} \ a \to - \infty \ {\rm and} \ b \to + \infty.$$
The conclusion then follows from \eqref{peillon}, and definitions \eqref{weakmass1} and \eqref{weakmass2}.
\end{proof}

When the function $1 - |\psi|^2$ belongs to $L^1(\R)$, the quantities $m^+(\psi)$ and $m^-(\psi)$ are equal to the mass of $\psi$ defined by \eqref{alamasse}. However, for an arbitrary map in $X_\boM(\R)$, the quantities $m^+(\psi)$ and $m^-(\psi)$, which are preserved by the flow, may be different. In order to define a generalized notion of mass, we are led to restrict ourselves to an even smaller class of maps. More precisely, we consider the subset of $X_\boM(\R)$ defined by
$$\boX_\boM(\R) = \{ \psi \in X_\boM(\R), \ {\rm s.t} \ m^+(\psi) = m^-(\psi) \},$$
and define the generalized mass of an arbitrary function $\psi \in \boX_\boM(\R)$ as the quantity
\begin{equation}
\label{goodmass}
m(\psi) = m^+(\psi) = m^-(\psi) = \frac{1}{2} \bigg( \underset{x \to + \infty}{\lim} \int_0^x \big( 1 - |\psi|^2 \big) + \underset{y \to - \infty}{\lim} \int_y^0 \big( 1 - |\psi|^2 \big) \bigg).
\end{equation}
We then have

\begin{prop}
\label{Consm}
Given any function $\Psi_0 \in \boX_\boM(\R)$, there exists a unique solution $\Psi(\cdot, t)$ to \eqref{GP} in $\boC^0(\R, \boX_\boM(\R))$ with initial datum $\Psi_0$. Moreover, we have
\begin{equation}
\label{consm}
m(\Psi(\cdot, t)) = m (\Psi_0),
\end{equation}
for any $t \in \R$.
\end{prop}

\begin{proof}
Proposition \ref{Wellposed} is a direct consequence of Lemmas \ref{Wellposed} and \ref{Conspm}. Given any function $\Psi_0 \in \boX_\boM(\R)$, there exists a unique solution $\Psi$ to \eqref{GP} in $\boC^0(\R, X_\boM(\R))$ with initial datum $\Psi_0$. Since $m^+(\Psi_0) = m^-(\Psi_0)$, it follows from Lemma \ref{Conspm} that
$$m^+(\Psi(\cdot, t)) = m^-(\Psi(\cdot, t)),$$
so that $\Psi(\cdot, t)$ belongs to $\boX_\boM(\R)$. Equality \eqref{consm} then follows from \eqref{plusoumoins} and \eqref{goodmass}.
\end{proof}

\begin{remark}
As already mentioned, if the function $1 - |\Psi_0|^2$ belongs to $L^1(\R)$, then the function $\Psi_0$ belongs to $\boX_\boM(\R)$, and it follows from Proposition \ref{Consm} that the generalized mass of the solution $\Psi(\cdot, t)$ to \eqref{GP} with initial datum $\Psi_0$ is well-defined for any time $t \in \R$ and conserved by the flow. Notice however that we do not claim that the function $1 - |\Psi|^2$ remains in $L^1(\R)$.
\end{remark}

\begin{remark}
In our proofs, we use several estimates involving control on the norm $\| \cdot \|_\boM$ which are closely related to the conservation of the generalized mass. The conservation of the generalized mass itself does actually not provide any bound on the solution $\Psi(\cdot, t)$. We believe that this fact is of independent interest. In particular, it might be relevant for the physical phenomena the equation was designed to describe.
\end{remark}

\begin{merci}
The authors are grateful to the referees for their careful reading of the paper, and their valuable remarks and comments which helped to improve the manuscript.\\
A large part of this work was completed while the four authors were visiting the Wolfgang Pauli Institute in Vienna. We wish to thank warmly this institution, as well as Prof. Norbert Mauser for the hospitality and support. We are also thankful to Dr. Martin Sepp for substantial digressions.\\
F.B., P.G. and D.S. are partially sponsored by project JC05-51279 of the Agence Nationale de la Recherche. J.-C. S. acknowledges support from project ANR-07-BLAN-0250 of the Agence Nationale de la Recherche.
\end{merci}

\bibliographystyle{plain}
\bibliography{Bibliogr}

\end{document}